\documentclass[10pt,a4paper]{article}
\usepackage{amsmath,amssymb,amsthm,graphics,graphicx,epic,eepic,multicol,ascmac}
\usepackage[mathscr]{euscript}
\usepackage[all]{xy}
\usepackage{fancybox}
\usepackage{url}
\usepackage{ytableau}

\numberwithin{equation}{section}
\theoremstyle{theorem}
\newtheorem{thm}{Theorem}[section]
\newtheorem{prop}[thm]{Proposition}
\newtheorem{lem}[thm]{Lemma}
\newtheorem{rem}[thm]{Remark}
\newtheorem{cor}[thm]{Corollary}

\theoremstyle{definition}
\newtheorem{defn}[thm]{Definition}
\newtheorem{ex}[thm]{Example}

\def\al{\alpha}

\def\wht(#1){\widehat{\ #1\ }}

\newcommand{\frg}{\mathfrak g}

\newcommand{\frt}{\mathfrak t}

\newcommand{\bbQ}{\mathbb Q}

\newcommand{\bbZ}{\mathbb Z}

\newcommand{\ch}{\mathrm{ch}}

\newcommand{\lbr}{\begin{bmatrix}}
\newcommand{\rbr}{\end{bmatrix}}

\newcommand{\cd}{commutative diagram }

\def\ge{\frg}

\def\al{\alpha}

\def\beneme{\begin{enumerate}}
\def\beq{\begin{equation}}
\def\beqn{\begin{eqnarray}}
\def\beqnn{\begin{eqnarray*}}

\def\bfii0{{\bf i_0}}

\def\bbra#1,#2,#3{\left\{\begin{array}{c}\hspace{-5pt}
#1;#2\\ \hspace{-5pt}#3\end{array}\hspace{-5pt}\right\}}
\def\cd{\cdots}
\def\ci(#1,#2){c_{#1}^{(#2)}}
\def\Ci(#1,#2){C_{#1}^{(#2)}}
\def\mpp(#1,#2,#3){#1^{(#2)}_{#3}}
\def\bCi(#1,#2){\ovl C_{#1}^{(#2)}}
\def\ch(#1,#2){c_{#2,#1}^{-h_{#1}}}
\def\cc(#1,#2){c_{#2,#1}}

\def\di(#1,#2){D_{#1}^{(#2)}}
\def\dbi(#1,#2){\ovl D_{#1}^{(#2)}}

\def\eit{\tilde{e}_i}
\def\eneme{\end{enumerate}}

\def\eeq{\end{equation}}
\def\eeqn{\end{eqnarray}}
\def\eeqnn{\end{eqnarray*}}
\def\fit{\tilde{f}_i}

\def\gau#1,#2{\left[\begin{array}{c}\hspace{-5pt}#1\\
\hspace{-5pt}#2\end{array}\hspace{-5pt}\right]}

\def\ify{\infty}
\def\io{\iota}
\def\ji(#1,#2){j_{#1}^{(#2)}}
\def\kp{k^{(+)}}
\def\km{k^{(-)}}

\def\lan{\langle}
\def\lar{\longrightarrow}
\def\lm{\lambda}

\def\nd{\noindent}

\def\ovl{\overline}

\def\qq{\qquad}

\def\qed{\hfill\framebox[2mm]{}}
\def\QQ{\mathbb Q}

\def\ran{\rangle}

\def\tt{\frt}
\def\TY(#1,#2,#3){#1^{(#2)}_{#3}}

\def\uq{U_q(\ge)}

\def\uqm{U^-_q(\ge)}

\def\vp{\varphi}

\def\xxi(#1,#2,#3){\displaystyle {}^{#1}\Xi^{(#2)}_{#3}}
\def\xsi(#1,#2,#3){\displaystyle {}^{#1}\Sigma^{(#2)}_{#3}}
\def\xE(#1,#2,#3){\displaystyle {}^{#1}E_{#2}[#3]}
\def\xF(#1,#2){\displaystyle {}^{#1}F_{#2}}
\def\xx(#1,#2){\displaystyle {}^{#1}\Xi_{#2}}
\def\W1{W(\varpi_1)}

\def\what{\widehat}

\def\ZZ{\mathbb Z}

\def\m@th{\mathsurround=0pt}
\def\fsquare(#1,#2){
\hbox{\vrule$\hskip-0.4pt\vcenter to #1{\normalbaselines\m@th
\hrule\vfil\hbox to #1{\hfill$\scriptstyle #2$\hfill}\vfil\hrule}$\hskip-0.4pt
\vrule}}

\newcommand{\ba}{\begin{array}}
\newcommand{\ea}{\end{array}}

\newcommand{\eq}{\begin{eqnarray}}
\newcommand{\eneq}{\end{eqnarray}}

\title{\textbf{\large{Adapted Sequences and Polyhedral Realizations of Crystal Bases for highest weight modules}}}
\author{\normalsize{YUKI KANAKUBO\thanks{Faculty of Pure and Applied Sciences, University of Tsukuba,
1-1-1 Tennodai, Tsukuba, Ibaraki 305-8577,
Japan: {y-kanakubo@math.tsukuba.ac.jp}.}
\ and\ 
TOSHIKI NAKASHIMA\thanks{Division of Mathematics, 
Sophia University, Kioicho 7-1, Chiyoda-ku, Tokyo 102-8554,
Japan: {toshiki@sophia.ac.jp}.}
}}
\date{}

\begin{document}

\maketitle
\vspace{-10pt}

\begin{abstract}
The polyhedral realizations for crystal bases of
the integrable highest weight modules of $U_q(\mathfrak{g})$ have been introduced in
([T.Nakashima, J. Algebra, vol.219, no. 2, (1999)]),
which describe the crystal bases as sets of lattice points in the infinite $\bbZ$-lattice $\mathbb{Z}^{\infty}$
given by some system of linear inequalities, where $\mathfrak{g}$ is a symmetrizable Kac-Moody Lie algebra.
To construct the polyhedral realization, we need to fix
an infinite sequence $\iota$ from the indices of the simple roots.
If the pair ($\iota$,$\lambda$) ($\lambda$: a dominant integral weight) satisfies the `ample' condition then
there are some procedure to calculate the sets of linear inequalities.

In this article, 
we show that if $\iota$ is an adapted sequence (defined in our paper [Y.Kanakubo, T.Nakashima, arXiv:1904.10919]) then
the pair ($\iota$, $\lambda$) satisfies the ample condition for any dominant integral weight $\lambda$ in the case $\mathfrak{g}$ is
a classical Lie algebra.
Furthermore, we reveal the explicit
forms of the polyhedral realizations of the crystal bases $B(\lambda)$
 associated with arbitrary 
adapted sequences $\iota$ in terms of column tableaux. As an application,
we will give a combinatorial description of the function $\varepsilon_i^*$ on the crystal base $B(\infty)$.
\end{abstract}

\section{Introduction}

The invention of crystal bases (\cite{K1,L}) developed the combinatorial
study of the quantum group $U_q(\mathfrak{g})$ and its representations,
where $\mathfrak{g}$ is a symmetrizable Kac-Moody Lie algebra with an index set $I=\{1,2,\cdots,n\}$.
The crystal bases $B(\lambda)$ tell us the skeleton structures
of the irreducible integrable highest weight modules $V(\lambda)$ and
are realized via combinatorial objects like as Young tableaux, LS paths, Laurent monomials, etc..

In \cite{NZ}, the polyhedral realization of crystal base $B(\infty)$ for the negative part 
$U^-_q(\mathfrak{g})$ has been introduced as an image of  `Kashiwara embedding'
$\Psi_{\iota}:B(\infty)\hookrightarrow \mathbb{Z}^{\infty}_{\iota}$, where $\iota$ is an infinite sequence of entries in $I$
and $\mathbb{Z}^{\infty}_{\iota}=\{(\cd,a_k,\cd,a_2,a_1)| a_k\in\ZZ
\,\,{\rm and}\,\,a_k=0\,\,{\rm for}\,\,k\gg 0\}$ is an infinite $\mathbb{Z}$-lattice with certain crystal structure associated with $\iota$.
Under the `positivity condition' on $\iota$, a procedure to describe
an explicit form of the image ${\rm Im}(\Psi_{\iota})$
is presented. 
If $\mathfrak{g}$ is a finite dimensional simple Lie algebra and $\iota=(\cdots,i_{N+1},i_N,\cdots,i_2,i_1)$ is a sequence such that
$(i_N,\cdots,i_2,i_1)$ is a reduced word of the longest element in the Weyl group $W$ then
the image ${\rm Im}(\Psi_{\iota})$ coincides with a set of lattice points in the string cone associated to the reduced word 
$(i_1,i_2\cdots,i_N)$ \cite{Lit}, which is a polyhedral convex cone \cite{BZ}.

The polyhedral realization for crystal bases $B(\lambda)$ is introduced 
as the image of embedding of crystals $\Psi^{\lambda}_{\iota}:B(\lambda)\hookrightarrow \mathbb{Z}^{\infty}_{\iota}[\lambda]$
(see \ref{poly-lm})
and
under the `ample condition' on the pair $(\iota,\lambda)$, an algorithm to calculate
${\rm Im}(\Psi^{\lambda}_{\iota})$ is presented in \cite{N99}.
In \cite{FN}, it is proved that ${\rm Im}(\Psi^{\lambda}_{\iota})$ is the set of lattice points
in a rational convex polytope (Newton-Okounkov convex body) for a finite dimensional semisimple Lie algebra $\mathfrak{g}$.
In \cite{H1, N99}, for the specific sequence $\iota=(\cdots,2,1,n,\cdots,2,1,n\cdots,2,1)$ and
simple Lie algebras $\mathfrak{g}$,
it is shown that the pair $(\iota,\lambda)$ ($\lambda$ is any dominant integral weight) satisfies
the ample condition and explicit forms of ${\rm Im}(\Psi^{\lambda}_{\iota})$ are given.

In \cite{KaN}, we considered a condition called `adaptedness' on sequences $\iota$ and proved that
if $\iota$ satisfies the adaptedness condition then the positivity condition holds in the case $\mathfrak{g}$ is a classical Lie algebra.
Using the method in \cite{NZ}, we also found explicit forms of polyhedral realizations
${\rm Im}(\Psi_{\iota})$ for $B(\infty)$ in terms of column tableaux. 
One defined a set ${\rm Tab}_{{\rm X},\iota}$ of column tableaux which lie in $(\mathbb{Q}^{\infty})^*$
and expressed ${\rm Im}(\Psi_{\iota})$ as ${\rm Im}(\Psi_{\iota})=\{\textbf{a}\in \mathbb{Z}^{\infty} |
\varphi(\textbf{a})\geq0\ {\rm for\ any\ }\varphi\in {\rm Tab}_{{\rm X},\iota} \}$.

In this article, for classical Lie algebras $\mathfrak{g}$,
we will prove that if $\iota$ satisfies the adaptedness condition then the pair
$(\iota,\lambda)$ is ample for any dominant integral weight $\lambda$.
One also give explicit forms of ${\rm Im}(\Psi^{\lambda}_{\iota})$ 
in terms of column tableaux. 
More precisely, defining the set ${\rm Tab}_{{\rm X},\iota}[\lambda]$ of column tableaux
which describe linear inequalities, we describe the polyhedral realization as
\[
{\rm Im}(\Psi^{\lambda}_{\iota})=\{\textbf{a}\in \mathbb{Z}^{\infty} |
\varphi(\textbf{a})\geq0\ {\rm for\ any\ }\varphi\in {\rm Tab}_{{\rm X},\iota}[\lambda]\cup {\rm Tab}_{{\rm X},\iota} \}.
\]
As an application, we get a tableaux description of $\varepsilon_i^*$, which is the composition of $\varepsilon_i$ on $B(\infty)$
and the operator $*:B(\infty)\rightarrow B(\infty)$ (see \ref{cry-sec}). 
As an example, let us consider the case $\mathfrak{g}$ is of type ${\rm A}_2$ and $\iota=(\cdots,2,1,2,1,2,1)$. 
We rewrite each element $(\cdots,a_3,a_2,a_1)$ in $\mathbb{Z}^{\infty}$ as $(\cdots,a_{3,2},a_{3,1},a_{2,2},a_{2,1},a_{1,2},a_{1,1})$.
Then
$\iota$ is adapted and hence the pair $(\iota,\lambda)$ satisfies the ample condition for any dominant integral weight $\lambda$.
We get
\begin{eqnarray*}
{\rm Tab}_{{\rm A},\iota}
&=&\{\begin{ytableau}
 i 
\end{ytableau}^{\rm A}_{s} | 1\leq i\leq 3,\ s\in\mathbb{Z}_{\geq1}\}\cup
\{\begin{ytableau}
 i \\
j
\end{ytableau}^{\rm A}_{s} | 1\leq i<j\leq 3,\ s\in\mathbb{Z}_{\geq1}\},\\
{\rm Tab}_{{\rm A},\iota}[\lambda]&=&
\{-x_{1,1}+\lan h_1,\lambda \ran\}\cup
\{\begin{ytableau}
2 \\
3
\end{ytableau}^{\rm A}_{0}+\lan h_2,\lambda\ran,
\begin{ytableau}
1 \\
3
\end{ytableau}^{\rm A}_{0}+\lan h_2,\lambda\ran
\},
\end{eqnarray*}
where the tableaux mean
\[
\begin{ytableau}
i 
\end{ytableau}^{\rm A}_{s}=x_{s,i}-x_{s+1,i-1},
\quad
\begin{ytableau}
 i \\
j
\end{ytableau}^{\rm A}_{s}=x_{s+1,i}-x_{s+2,i-1} +x_{s,j}-x_{s+1,j-1},
\]
and each $x_{s,i}\in (\mathbb{Q}^{\infty})^*$ is defined as
$x_{s,i}(\cdots,a_{3,2},a_{3,1},a_{2,2},a_{2,1},a_{1,2},a_{1,1})=a_{s,i}$ for $s\in \mathbb{Z}_{\geq1}$, $i\in\{1,2\}$
and $x_{s,i}=0$ if $i\notin\{1,2\}$.
Thus
\[
{\rm Im}(\Psi^{\lambda}_{\iota})=\left\{ (\cdots,a_{2,2},a_{2,1},a_{1,2},a_{1,1})\in\mathbb{Z}^{\infty}
\left|\begin{array}{l}
 a_{s,i}-a_{s+1,i-1}\geq0\ (1\leq i\leq3,\ s\in\mathbb{Z}_{\geq1}),\\
 a_{s+1,i}-a_{s+2,i-1} +a_{s,j}-a_{s+1,j-1}\geq0\ (1\leq i<j\leq 3,\ s\in\mathbb{Z}_{\geq1}),\\
 -a_{1,1}+\lan h_1,\lambda \ran\geq0,\\
 -a_{1,2}+a_{1,1}+\lan h_2,\lambda \ran\geq0,\ 
  -a_{2,1}+\lan h_2,\lambda \ran\geq0.
\end{array}\right.
\right\}.
\]
Simplifying the inequalities, we get 
\[
{\rm Im}(\Psi^{\lambda}_{\iota})=
\left\{(\cdots,a_{2,2},a_{2,1},a_{1,2},a_{1,1})\in\mathbb{Z}^{\infty} 
\left|\begin{array}{l}  
a_{1,2}\geq a_{2,1}\geq0,\ a_{1,1}\geq0,\ a_{m+1,2}=a_{m+2,1}=0,\ \forall m\in\mathbb{Z}_{\geq1}\\
\lan h_1,\lambda \ran\geq a_{1,1},\ 
-a_{1,2}+a_{1,1}+\lan h_2,\lambda \ran\geq0,\ 
 \lan h_2,\lambda \ran\geq a_{2,1}.
\end{array}\right.
\right\}.
\]

We explain our motivation
to describe inequalities defining polyhedral realizations
as column tableaux.
In \cite{Nj,NN}, $q$-characters are expressed in terms of
tableaux. In \cite{Nj}, monomial realizations of crystal bases
are introduced and it is proved that the sets of monomials appearing in ($t$-analogs of)
$q$-characters of standard modules are isomorphic to crystal bases of highest weight
$\mathfrak{g}$-modules in the case $\mathfrak{g}$ is simply laced.
On the other hand,
in \cite{N13,N14}, a relation between
monomial realizations
and inequalities
defining polyhedral realizations
are referred. Thus, there seems to be an indirect linkage
between the inequalities and $q$-characters, 
which gives us a motivation of 
column tableaux expression of the inequalities.

The organization of this article is as follows. In Sect.2, after a concise reminder on 
crystals, we review the polyhedral realizations for $B(\infty)$ and $B(\lambda)$. 
In Sect.3, we recall the column tableaux descriptions of the
polyhedral realizations for $B(\infty)$, which is shown in our previous article \cite{KaN}.
Sect.4 is devoted to present our main results, which provide a
column tableaux descriptions of the polyhedral realizations for $B(\lambda)$.
In Sect.5, 6, 
we will prove our main theorem.

\vspace{2mm}

\nd
{\bf Acknowledgements} 
Y.K. was supported by JSPS KAKENHI Grant Number JP20J00186.
T.N. was supported in part by JSPS KAKENHI Grant Number JP15K04794.

\section{Crystal and its polyhedral realization}

Let us recall the definition of {\it crystals} \cite{K3}. 

\subsection{Notations}

We list the notations used in this paper.
Let $\ge$ be
a symmetrizable Kac-Moody Lie algebra over $\bbQ$
with a Cartan subalgebra
$\tt$, a weight lattice $P \subset \tt^*$, the set of simple roots
$\{\al_i| i\in I\} \subset \tt^*$,
and the set of coroots $\{h_i| i\in I\} \subset \tt$,
where $I=\{1,2,\cdots,n\}$ is a finite index set.
Let $\lan h,\lm\ran=\lm(h)$ be the pairing between $\tt$ and $\tt^*$,
and $(\al, \beta)$ be an inner product on
$\tt^*$ such that $(\al_i,\al_i)\in 2\bbZ_{\geq 0}$ and
$\lan h_i,\lm\ran={{2(\al_i,\lm)}\over{(\al_i,\al_i)}}$
for $\lm\in\tt^*$ and $A:=(\lan h_i,\al_j\ran)_{i,j\in I}=(a_{i,j})_{i,j\in I}$ be the associated symmetrizable generalized Cartan matrix.
Let $P^*=\{h\in \tt| \lan h,P\ran\subset\ZZ\}$ and
$P_+:=\{\lm\in P|\lan h_i,\lm\ran\in\ZZ_{\geq 0}\}$.
We call an element in $P_+$ a {\it dominant integral weight}.
The quantum algebra $\uq$
is an associative
$\QQ(q)$-algebra generated by $e_i$, $f_i \,\, (i\in I)$,
and $q^h \,\, (h\in P^*)$
satisfying the usual relations.
The algebra $\uqm$ is the subalgebra of $\uq$ generated 
by $f_i$ $(i\in I)$.

For the irreducible integrable highest weight module of $\uq$
with the highest weight $\lm\in P_+$, we denote it by $V(\lm)$
and its {\it crystal base} we denote $(L(\lm),B(\lm))$.
Similarly, for the crystal base of the algebra $\uqm$ we denote 
$(L(\ify),B(\ify))$ (see \cite{K0,K1}).
For positive integers $l$ and $m$ with $l\leq m$, we set $[l,m]:=\{l,l+1,\cdots,m-1,m\}$.

\subsection{Crystals}\label{cry-sec}

By the terminology {\it crystal } we mean some combinatorial object 
obtained by abstracting the properties of crystal bases:

\begin{defn}
A {\it crystal} is a set $\mathcal{B}$ together with the maps
${\rm wt}:\mathcal{B}\rightarrow P$,
$\varepsilon_i,\varphi_i:\mathcal{B}\rightarrow \mathbb{Z}\sqcup \{-\infty\}$
and $\tilde{e}_i$,$\tilde{f}_i:\mathcal{B}\rightarrow \mathcal{B}\sqcup\{0\}$
($i\in I$) satisfying the following: For $b,b'\in\mathcal{B}$, $i\in I$,
\begin{enumerate}
\item[$(1)$] $\varphi_i(b)=\varepsilon_i(b)+\langle h_i,{\rm wt}(b)\rangle$,
\item[$(2)$] ${\rm wt}(\tilde{e}_ib)={\rm wt}(b)+\alpha_i$ if $\tilde{e}_i(b)\in\mathcal{B}$,
\quad ${\rm wt}(\tilde{f}_ib)={\rm wt}(b)-\alpha_i$ if $\tilde{f}_i(b)\in\mathcal{B}$,
\item[$(3)$] $\varepsilon_i(\tilde{e}_i(b))=\varepsilon_i(b)-1,\ \ 
\varphi_i(\tilde{e}_i(b))=\varphi_i(b)+1$\ if $\tilde{e}_i(b)\in\mathcal{B}$, 
\item[$(4)$] $\varepsilon_i(\tilde{f}_i(b))=\varepsilon_i(b)+1,\ \ 
\varphi_i(\tilde{f}_i(b))=\varphi_i(b)-1$\ if $\tilde{f}_i(b)\in\mathcal{B}$, 
\item[$(5)$] $\tilde{f}_i(b)=b'$ if and only if $b=\tilde{e}_i(b')$,
\item[$(6)$] if $\varphi_i(b)=-\infty$ then $\tilde{e}_i(b)=\tilde{f}_i(b)=0$.
\end{enumerate}
We call $\tilde{e}_i$,$\tilde{f}_i$ {\it Kashiwara operators}.
\end{defn}

\begin{defn}
A {\it strict morphism} $\psi:\mathcal{B}_1\rightarrow\mathcal{B}_2$
of crystals $\mathcal{B}_1$, $\mathcal{B}_2$ is a map 
$\mathcal{B}_1\sqcup\{0\}\rightarrow\mathcal{B}_2\sqcup\{0\}$
satisfying the following conditions:
 $\psi(0)=0$, 
${\rm wt}(\psi(b))={\rm wt}(b)$,
 $\varepsilon_i(\psi(b))=\varepsilon_i(b)$,
 $\varphi_i(\psi(b))=\varphi_i(b)$,
 if $b\in \mathcal{B}_1$, $\psi(b)\in \mathcal{B}_2,$ $i\in I$,
 and
$\psi: \mathcal{B}_1\sqcup\{0\} \lar \mathcal{B}_2\sqcup\{0\}$
commutes with all $\eit$ and $\fit$, where $\eit(0)=\fit(0)=0$.
An injective strict morphism is said to be {\it strict embedding} of crystals.
\end{defn}

Let ($B(\infty),\{\tilde{e}_i\}_{i\in I},\{\tilde{f}_i\}_{i\in I},\{\varepsilon_i\}_{i\in I},\{\varphi_i\}_{i\in I}, {\rm wt})$ be the crystal
structure of $B(\infty)$ and
$*:U_q(\frak{g})\rightarrow U_q(\frak{g})$ the antiautomorphism
such that $e_i^*=e_i$, $f_i^*=f_i$ and $(q^h)^*=q^{-h}$ in \cite{K1}.
It is known that the map $*$ induces a bijection $*:B(\infty)\rightarrow B(\infty)$ satisfying $*\circ *=id$.
Let $B(\infty)^{*}$ be the crystal as follows : $B(\infty)^{*}$ is equal to $B(\infty)$ as sets,
and maps are defined as $\tilde{e}_i^*:=*\circ \tilde{e}_i\circ *$, 
$\tilde{f}_i^*:=*\circ \tilde{f}_i\circ *$,
$\varepsilon_i^*:=\varepsilon_i\circ *$,
$\varphi_i^*:=\varphi_i\circ *$ and ${\rm wt}^*:={\rm wt}$.

\subsection{Polyhedral realizations of $B(\ify)$}
\label{poly-uqm}
Let us recall the results in \cite{NZ}.

First, we consider the infinite $\bbZ$-lattice
\[
\ZZ^{\ify}
:=\{(\cd,a_k,\cd,a_2,a_1)| a_k\in\ZZ
\,\,{\rm and}\,\,a_k=0\,\,{\rm for}\,\,k\gg 0\}
\]
and denote by $\ZZ^{\ify}_{\geq 0} \subset \ZZ^{\ify}$
the subsemigroup of nonnegative sequences.
For the rest of this section, we fix an infinite sequence of indices
$\io=(\cd,i_k,\cd,i_2,i_1)$ from $I$ such that
\begin{equation}
{\hbox{
$i_k\ne i_{k+1}$ and $\sharp\{k\in\mathbb{Z}_{\geq1}| i_k=i\}=\ify$ for any $i\in I$.}}
\label{seq-con}
\end{equation}
We can define a crystal structure
on $\ZZ^{\ify}$ associated with $\iota$ and denote it by $\ZZ^{\ify}_{\io}$ 
(\cite[2.4]{NZ}).

\begin{prop}[\cite{K3}, See also \cite{NZ}]
\label{emb}
There is a unique strict embedding of crystals
$($called Kashiwara embedding$)$
\begin{equation}
\Psi_{\io}:B(\ify)\hookrightarrow \ZZ^{\ify}_{\geq 0}
\subset \ZZ^{\ify}_{\io},
\label{psi}
\end{equation}
such that
$\Psi_{\io} (u_{\ify}) = (\cd,0,\cd,0,0)$, where 
$u_{\ify}\in B(\ify)$ is the vector corresponding to $1\in \uqm$.
\end{prop}

\begin{defn}
The image ${\rm Im} \Psi_{\iota} (\cong B(\infty))$ is called a \it{polyhedral realization} of $B(\infty)$.
\end{defn}

Let us consider the infinite dimensional vector space
$$
\QQ^{\ify}:=\{\textbf{a}=
(\cd,a_k,\cd,a_2,a_1)| a_k \in \QQ\,\,{\rm and }\,\,
a_k = 0\,\,{\rm for}\,\, k \gg 0\},
$$
and its dual space $(\QQ^{\ify})^*:={\rm Hom}(\QQ^{\ify},\QQ)$.
Let $x_k\in (\QQ^{\ify})^*$ be the linear function defined as $x_k((\cd,a_k,\cd,a_2,a_1)):=a_k$
for $k\in\mathbb{Z}_{\geq1}$.
We will also write a linear form $\vp \in (\QQ^{\ify})^*$ as
$\vp=\sum_{k \geq 1} \vp_k x_k$ ($\vp_j\in \QQ$).

For the fixed infinite sequence
$\io=(i_k)_{k\in\mathbb{Z}_{\geq1}}$ and $k\geq1$ we set $\kp:={\rm min}\{l:l>k\,\,{\rm and }\,\,i_k=i_l\}$ and
$\km:={\rm max}\{l:l<k\,\,{\rm and }\,\,i_k=i_l\}$ if it exists,
or $\km=0$  otherwise.
We set $\beta_0:=0$ and
\begin{equation}
\beta_k:=x_k+\sum_{k<j<\kp}\lan h_{i_k},\al_{i_j}\ran x_j+x_{\kp}\in (\QQ^{\ify})^*
\qq(k\geq1).
\label{betak}
\end{equation}
We define the piecewise-linear operator 
$S_k=S_{k,\io}$ on $(\QQ^{\ify})^*$ as
\begin{equation}
S_k(\vp):=
\begin{cases}
\vp-\vp_k\beta_k & {\mbox{ if }}\vp_k>0,\\
 \vp-\vp_k\beta_{\km} & {\mbox{ if }}\vp_k\leq 0.
\end{cases}
\label{Sk}
\end{equation}
We set
\begin{eqnarray}
\Xi_{\io} &:=  &\{S_{j_l}\cd S_{j_2}S_{j_1}x_{j_0}\,|\,
l\geq0,j_0,j_1,\cd,j_l\geq1\},
\label{Xi_io} \label{xiiodef}\\
\Sigma_{\io} & := &
\{\textbf{x}\in \ZZ^{\ify}\subset \QQ^{\ify}\,|\,\vp(\textbf{x})\geq0\,\,{\rm for}\,\,
{\rm any}\,\,\vp\in \Xi_{\io}\}.
\end{eqnarray}
We impose on $\io$ the following {\it positivity condition}:
\begin{equation}
{\hbox{if $\km=0$ then $\vp_k\geq0$ for any 
$\vp=\sum_k\vp_kx_k\in \Xi_{\io}$}}.
\label{posi}
\end{equation}
\begin{thm}[\cite{NZ}]\label{polyhthm}
Let $\io$ be a sequence of indices satisfying $(\ref{seq-con})$ 
and $(\ref{posi})$. Then it holds 
${\rm Im}(\Psi_{\io})=\Sigma_{\io}$.
\end{thm}

\subsection{Polyhedral realizations of $B(\lambda)$}
\label{poly-lm}

For $\lambda\in P_+$, let $B(\lambda)$ be the crystal base 
of the irreducible integrable highest weight module $V(\lambda)$ with 
the highest weight vector $u_{\lambda}$. 
We fix an infinite sequence of indices
$\io=(\cd,i_k,\cd,i_2,i_1)$ from $I$ satisfying (\ref{seq-con}).
For $\lambda\in P_+$,
let $\mathbb{Z}^{\infty}_{\iota}[\lambda]$ be the crystal 
whose underlying set is $\mathbb{Z}^{\infty}$ and whose crystal structure associated with $\iota$
is defined as in Sect.4.1 of \cite{N99}.

\begin{thm}\cite{N99}\label{emb2}
There exists the unique strict embedding of crystals
\[
\Psi^{(\lambda)}_{\io}:B(\lambda)\hookrightarrow \ZZ^{\ify}_{\iota}[\lambda]
\]
such that $\Psi^{(\lambda)}_{\io}(u_{\lambda})=(\cdots,0,0,0)$.
\end{thm}

\begin{defn}
The image ${\rm Im} \Psi^{(\lambda)}_{\iota} (\cong B(\lambda))$ is called a \it{polyhedral realization} of $B(\lambda)$.
\end{defn}

Let $\beta^{(\pm)}_k$ be the linear functions given by $\beta^{(+)}_k=\beta_k$ and
\begin{equation}\label{betadef}
\beta^{(-)}_k = 
\begin{cases}
x_{k^{(-)}}+\sum_{k^{(-)}<j<k} \lan h_{i_k},\alpha_{i_j}\ran x_j +x_k & {\rm if}\ k^{(-)}>0,\\
- \lan h_{i_k},\lambda\ran + \sum_{1\leq j<k} \lan h_{i_k},\alpha_{i_j}\ran x_j + x_k & {\rm if}\ k^{(-)}=0.
\end{cases}
\end{equation}
Note that $\beta^{(-)}_k= \beta_{k^{(-)}}$ if  $k^{(-)}>0$.

Using this notation, for each $k\in\mathbb{Z}_{\geq1}$, we define an operator $\what{S}_k=\what{S}_{k,\io}$ as follows:
For a linear function $\varphi=c+\sum_{k\geq1} \varphi_kx_k$ ($c$, $\varphi_k\in\mathbb{Q}$),
we set
\begin{equation}\label{Shat}
\what{S}_k(\varphi):=
\begin{cases}
\varphi - \varphi_k\beta^{(+)}_k & {\rm if}\ \varphi_k>0, \\
\varphi - \varphi_k\beta^{(-)}_k & {\rm if}\ \varphi_k\leq0.
\end{cases}
\end{equation}
One can easily check $(\what{S}_k)^2=\what{S}_k$. For $i\in I$, we also set 
\begin{equation}\label{lmi-def}
\lambda^{(i)} :=\lan h_i, \lambda \ran - \sum_{1\leq j<\iota^{(i)}} \lan h_i, \alpha_{i_j} \ran x_j - x_{\iota^{(i)}},
\end{equation}
where $\iota^{(i)}:={\rm min}\{k\in\mathbb{Z}_{\geq1} | i_k=i\}$.
For $\iota$ and a dominant integral weight $\lambda$, let $\Xi_{\iota}[\lambda]$ be the set of all linear
functions generated by applying $\what{S}_k$ on the functions $x_j$ ($j \geq 1$) and $\lambda^{(i)}$
($i\in I$), namely,
\begin{equation}\label{xilamdef}
\begin{array}{l}
\Xi_{\iota}[\lambda]: =\{\what{S}_{j_l}\cdots \what{S}_{j_1}x_{j_0} | l\in\mathbb{Z}_{\geq0}, j_0,\cdots,j_l\in \mathbb{Z}_{\geq1} \} \\
\qquad\qquad \cup \{\what{S}_{j_k}\cdots \what{S}_{j_1} \lambda^{(i)} | k\in\mathbb{Z}_{\geq0},\ i\in I,\ j_1,\cdots,j_k\in \mathbb{Z}_{\geq1} \}.
\end{array}
\end{equation}
Now we set
\begin{equation}\label{siglamdef}
\Sigma_{\iota}[\lambda]:=
\{
x\in\mathbb{Z}^{\infty}_{\iota}[\lambda] | \varphi(x) \geq0\ \ {\rm for\ any\ }\varphi\in \Xi_{\iota}[\lambda] 
\}.
\end{equation}
\begin{defn}
We say the pair $(\iota,\lambda)$ is {\it ample}
if $\textbf{0}:=(\cdots,0,0,0)\in \Sigma_{\iota}[\lambda]$.
\end{defn}

\begin{thm}\label{Nthm1}\cite{N99}
We suppose that $(\iota,\lambda)$ is ample. 
Let $\Psi^{(\lambda)}_{\iota}:B(\lambda)\hookrightarrow \mathbb{Z}^{\infty}_{\iota}[\lambda]$
be the embedding as in Theorem \ref{emb2}. Then
the image ${\rm Im}(\Psi^{(\lambda)}_{\iota}) (\cong B(\lambda))$ is equal to 
$\Sigma_{\iota}[\lambda]$.
\end{thm}

\begin{ex}

Let $\mathfrak{g}$ be of type ${\rm A}_2$, $\iota=(\cdots,2,1,2,1,2,1)$ and $\lambda\in P_+$.
It follows
\[
1^{(-)}=2^{(-)}=0,\ k^{(-)}>0\ (k>2).\]
We rewrite a vector $(\cdots,x_6,x_5,x_4,x_3,x_2,x_1)$ as 
\[
(\cdots,x_{3,2},x_{3,1},x_{2,2},x_{2,1},x_{1,2},x_{1,1}), \]
that is, $x_{2l-1}=x_{l,1}$, $x_{2l}=x_{l,2}$ for $l\in\mathbb{Z}_{\geq1}$.
Similarly, we rewrite $\what{S}_{2l-1}=\what{S}_{l,1}$, $\what{S}_{2l}=\what{S}_{l,2}$.
For $k\in\mathbb{Z}_{\geq1}$, the action of the operators are the following:
\[
x_{k,1}\overset{\what{S}_{k,1}}{\rightleftarrows}x_{k,2}-x_{k+1,1}\overset{\what{S}_{k,2}}{\rightleftarrows} -x_{k+1,2},
\]
\vspace{-7mm}
\[\hspace{-6mm}
 {\scriptstyle \what{S}_{k+1,1}}\quad \quad \quad \quad \quad {\scriptstyle \what{S}_{k+1,2}}
\]

\[
x_{k,2}\overset{\what{S}_{k,2}}{\rightleftarrows}x_{k+1,1}-x_{k+1,2}\overset{\what{S}_{k+1,1}}{\rightleftarrows} -x_{k+2,1},
\]
\vspace{-6mm}
\[\hspace{-7mm}
 {\scriptstyle \what{S}_{k+1,2}}\quad \quad \quad \quad \quad \quad {\scriptstyle \what{S}_{k+2,1}}
\]
and the other actions are trivial. Thus we obtain
\[
\{\what{S}_{j_l}\cdots \what{S}_{j_1}x_{j_0} | l\in\mathbb{Z}_{\geq0}, j_0,\cdots,j_l\in \mathbb{Z}_{\geq1} \} 
=\{x_{k,1},\ x_{k,2}-x_{k+1,1},\ -x_{k+1,2},\ x_{k,2},\ x_{k+1,1}-x_{k+1,2}, -x_{k+2,1} | k\geq1\}.
\]
The definition (\ref{lmi-def}) of $\lambda^{(i)}$ means that
$\lambda^{(1)}=\lan h_1,\lambda \ran- x_{1,1}$ and $\lambda^{(2)}=\lan h_2,\lambda \ran+x_{1,1}- x_{1,2}$.
We also obtain
\[
\lambda^{(1)} \overset{\what{S}_{1,1}}{\longrightarrow} 0,\ \ 
\lambda^{(2)} \overset{\what{S}_{1,2}}{\longrightarrow} 0,\ \ 
\lambda^{(2)} \overset{\what{S}_{1,1}}{\longrightarrow} \lan h_2,\lambda \ran-x_{2,1}
\overset{\what{S}_{2,1}}{\longrightarrow}\lambda^{(2)},
\]
which means that
\[
\{\what{S}_{j_k}\cdots \what{S}_{j_1} \lambda^{(1)} | k\in\mathbb{Z}_{\geq0},\ j_1,\cdots,j_k\in \mathbb{Z}_{\geq1} \}
=\{0,\lambda^{(1)}\},
\]
\[
\{\what{S}_{j_k}\cdots \what{S}_{j_1} \lambda^{(2)} | k\in\mathbb{Z}_{\geq0},\ j_1,\cdots,j_k\in \mathbb{Z}_{\geq1} \}
=\{0,\lambda^{(2)},\ \lan h_2,\lambda \ran-x_{2,1}\}.
\]
Thus,
\[
\Xi_{\iota}[\lambda] =\{x_{k,1},\ x_{k,2}-x_{k+1,1},\ -x_{k+1,2},\ x_{k,2},\ x_{k+1,1}-x_{k+1,2}, -x_{k+2,1} | k\geq1\}
\cup 
\{0,\lambda^{(1)},\lambda^{(2)},\ \lan h_2,\lambda \ran-x_{2,1}\}
\]
and it is easy to see $\textbf{0}\in \Sigma_{\iota}[\lambda]$. Hence $(\iota,\lambda)$ is ample and
$\Sigma_{\iota}[\lambda]={\rm Im}(\Psi_{\iota}^{(\lambda)})$.
For $\textbf{x}=(\cdots,x_{3,2},x_{3,1},x_{2,2},x_{2,1},x_{1,2},x_{1,1})\in {\rm Im}(\Psi^{(\lambda)}_{\iota})$,
combining the inequalities $x_{k,1}\geq0$, $-x_{k+2,1}\geq0$ $(k\geq1)$ in $\Sigma_{\iota}[\lambda]$, we obtain
$x_{k,1}=0$ $(k\geq3)$. Similarly, by $x_{k,2}\geq0$, $-x_{k+1,2}\geq0$ $(k\geq1)$,
we get $x_{k,2}=0$ $(k\geq2)$. Hence, we obtain
\[
 {\rm Im}(\Psi^{(\lambda)}_{\iota})=\Sigma_{\io}[\lambda]=
\left\{\textbf{x}\in \mathbb{Z}^{\infty}_{\iota} \left|
\begin{array}{l}
x_{k+1,1}=x_{k,2}=0\ {\rm for}\ k\in\mathbb{Z}_{\geq2},\\ 
x_{1,2}\geq x_{2,1}\geq 0,\ x_{1,1}\geq0,\ \lan h_1,\lambda \ran\geq x_{1,1},\\
\lan h_2,\lambda \ran\geq x_{1,2}-x_{1,1},\ \lan h_2,\lambda \ran\geq x_{2,1}\end{array}\right.
\right\}.
\]

\end{ex}

\begin{ex}\cite{N99}\label{counterofpos}
Let $\mathfrak{g}$ be of type ${\rm A}_3$, $\iota=(\cdots,2,1,2,3,2,1)$ and $\lambda\in P_+$ such that $\lan h_2,\lambda\ran>0$.
We obtain
\[
x_{1} \overset{\what{S}_1}{\rightarrow} -x_5+x_4+x_2 \overset{\what{S}_2}{\rightarrow} -x_5 +x_3  
\overset{\what{S}_5}{\rightarrow} -x_4+x_3-x_2+x_1
\overset{\what{S}_2}{\rightarrow} -x_4+x_3-\lan h_2,\lambda\ran.
\]
Thus $\varphi:=-x_4+x_3-\lan h_2,\lambda\ran\in \Xi_{\iota}[\lambda]$, $\varphi(\textbf{0})=-\lan h_2,\lambda\ran<0$
and $(\iota,\lambda)$ is not ample.
\end{ex}

\subsection{Strict positivity condition}

Let $\xi^{(i)}$ ($i\in I$) be the linear function on $\mathbb{Q}^{\infty}$ defined by
\begin{equation}
\xi^{(i)}:=
- \sum_{1\leq j<\iota^{(i)}} \lan h_i, \alpha_{i_j} \ran x_j - x_{\iota^{(i)}}.
\end{equation}
Note that for any $\lambda\in P_+$, it follows
$\xi^{(i)}= -\lan h_i,\lambda \ran+\lambda^{(i)}$.
We define the set of linear forms 
\begin{equation}
\Xi_{\iota}^{(i)}:= 
\{
S_{j_l}\cdots S_{j_1} \xi^{(i)} | l\geq0,\ j_1,\cdots,j_l\geq1
\},
\end{equation}
and $\Xi^{(\infty)}_{\iota}:= \Xi_{\iota}$ in (\ref{xiiodef}).

\begin{defn}\cite{N99}
We say $\iota$ satisfies the {\it strict positivity condition} if the following condition holds:
\[
{\rm if}\ k^{(-)} = 0\ {\rm then}\ \varphi_k\geq0 \ {\rm for\ any}\ \varphi=\sum_{k}\varphi_kx_k
\in (\bigcup_{j\in I\cup\{\infty\}} \Xi^{(j)}_{\iota}) \setminus \{\xi^{(i)} | i\in I \}.
\]
\end{defn}

\begin{thm}\cite{N99}
Let $\iota$ be a sequence of indices satisfying (\ref{seq-con}) and the strict positivity condition, and $\lambda$ be a dominant
integral weight.
Then for $i\in I$ and $x\in \Sigma_{\iota}$, we get
\begin{equation}
\varepsilon^*_i(x) = {\rm max}\{-\varphi(x) | \varphi\in\Xi^{(i)}_{\iota}\}.
\end{equation} 
\end{thm}

\subsection{Infinite sequences adapted to $A$}

\begin{defn}\label{adapt}\cite{KaN}
Let $A=(a_{i,j})$ be the symmetrizable generalized Cartan matrix of $\mathfrak{g}$ and $\io$
a sequence of indices satisfying $(\ref{seq-con})$.
If $\iota$ satisfies the following condition then we say $\iota$ is {\it adapted} to $A$ : 
For $i,j\in I$ with $i\neq j$ and $a_{i,j}\neq0$, the subsequence of $\iota$ consisting of all $i$, $j$ is
\[
(\cdots,i,j,i,j,i,j,i,j)\quad {\rm or}\quad (\cdots,j,i,j,i,j,i,j,i).
\]
If the Cartan matrix is fixed then the sequence $\iota$ is shortly said to be {\it adapted}.
\end{defn}

\begin{ex}
We consider the case
 $\mathfrak{g}$ is of type ${\rm A}_3$, $\iota=(\cdots,2,1,3,2,1,3,2,1,3)$.
\begin{itemize}
\item
The subsequence consisting of $1$, $2$ is $(\cdots,2,1,2,1,2,1)$.
\item
The subsequence consisting of $2$, $3$ is $(\cdots,2,3,2,3,2,3)$.
\item
Since $a_{1,3}=0$ we do not need consider the pair $1$, $3$.
\end{itemize}

Hence $\iota$ is an adapted sequence.
\end{ex}

\begin{ex}
We consider the case
 $\mathfrak{g}$ is  of type ${\rm A}_3$, $\iota=(\cdots,2,1,2,3,2,1)$.
The subsequence consisting of $1$, $2$ is $(\cdots,2,1,2,2,1)$.
Thus $\iota$ is not an adapted sequence.
\end{ex}

\section{Tableaux descriptions of polyhedral realizations of $B(\infty)$}

In this section, we take $\mathfrak{g}$ as a finite dimensional simple Lie algebra of type ${\rm A}_n$,
${\rm B}_n$, ${\rm C}_n$ or ${\rm D}_n$. We shall review the column tableaux descriptions of the
polyhedral realizations for $B(\infty)$ shown in our previous article \cite{KaN}.
In the rest of article, we follow Kac's notation \cite{Kac} and
suppose that $\iota=(\cdots,i_3,i_2,i_1)$ satisfies
$(\ref{seq-con})$ and is adapted to the Cartan matrix $A$ of $\mathfrak{g}$.
Let $(p_{i,j})_{i\neq j,\ a_{i,j}\neq0}$ be the set of integers such that
\begin{equation}\label{pij}
p_{i,j}=\begin{cases}
1 & {\rm if}\ {\rm the\ subsequence\ of\ }\iota{\rm\ consisting\ of}\ i,j\ {\rm is}\ (\cdots,j,i,j,i,j,i),\\
0 & {\rm if}\ {\rm the\ subsequence\ of\ }\iota{\rm\ consisting\ of}\ i,j\ {\rm is}\ (\cdots,i,j,i,j,i,j).
\end{cases}
\end{equation}
For $k$ $(2\leq k\leq n)$, we set
\begin{equation}\label{LaP}
P(k):=\begin{cases}
p_{2,1}+p_{3,2}+\cdots+p_{n-2,n-3}+p_{n,n-2} & {\rm if}\ k=n\ {\rm and}\ \mathfrak{g}\ {\rm is\ of}\ {\rm type\ D}_n, \\
p_{2,1}+p_{3,2}+p_{4,3}+\cdots+p_{k,k-1} & {\rm if}\ {\rm otherwise}
\end{cases}
\end{equation}
and $P(0)=P(1)=P(n+1)=0$. Since each $p_{i,j}$ is in $\{0,1\}$,
it holds for $k$, $l\in I$ such that $k\geq l$.
\begin{equation}\label{pineq}
P(k)\geq P(l),
\end{equation}
\begin{equation}\label{pineq2}
(k-l)+P(l)\geq P(k),
\end{equation}
except for the case $\mathfrak{g}$ is of type ${\rm D}_n$, $k=n$ and $l=n-1$. 

For $k\in\mathbb{Z}_{\geq1}$, we rewrite $x_k$, $\beta_k$ and $S_k$ in \ref{poly-uqm} as
\begin{equation}\label{rewrite}
x_k= x_{s,j},\quad S_k=S_{s,j},\quad \beta_k= \beta_{s,j} 
\end{equation}
if $i_k=j$ and $j$ is appearing $s$ times in $i_k$, $i_{k-1}$, $\cdots,i_1$.
For example, if $\iota=(\cdots,2,1,3,2,1,3,2,1,3)$ then
we rewrite $(\cdots,x_6,x_5,x_4,x_3,x_2,x_1)=(\cdots,x_{2,2},x_{2,1},x_{2,3},x_{1,2},x_{1,1},x_{1,3})$.
We will use the both notation $x_k$ and $x_{s,j}$.

\begin{defn}
Let us define the (partially) ordered sets $J_{\rm A}$, $J_{\rm B}$, $J_{\rm C}$ and $J_{\rm D}$ as follows:

\begin{itemize}
\item $J_{\rm A}:=\{1,2,\cdots,n,n+1\}$ with the order $1<2<\cdots<n<n+1$.
\item
$J_{\rm B}=J_{\rm C}:=\{1,2,\cdots,n,\overline{n},\cdots,\overline{2},\overline{1}\}$ with the order
\[
1<2<\cdots<n<\overline{n}<\cdots<\overline{2}<\overline{1}.\]
\item
$J_{\rm D}:=\{1,2,\cdots,n,\overline{n},\cdots,\overline{2},\overline{1}\}$ with the partial order 
\[
 1< 2<\cdots< n-1<\ ^{n}_{\overline{n}}\ < \overline{n-1}< \cdots< \overline
 2< \overline 1.\]
\end{itemize}
For $j\in\{1,2,\cdots,n\}$, we set $|j|=|\overline{j}|=j$.

\end{defn}

\begin{defn}\label{box-def}\cite{KaN}
\begin{enumerate}
\item For $j\in[1,n+1]$ and $s\in\mathbb{Z}$, we set
\[
\fbox{$j$}^{\rm A}_{s}:=x_{s+P(j),j}-x_{s+P(j-1)+1,j-1}\in (\QQ^{\ify})^*,\]
where $x_{m,0}=x_{m,n+1}=0$ for $m\in\mathbb{Z}$, and $x_{m,i}=0$ for $m\in\mathbb{Z}_{\leq0}$ and $i\in I$.
\item
For $j\in[1,n]$ and $s\in\mathbb{Z}$, we set
\[
\fbox{$j$}^{\rm B}_{s}
:=x_{s+P(j),j}-x_{s+P(j-1)+1,j-1}\in (\QQ^{\ify})^*,\]
\[
\fbox{$\overline{j}$}^{\rm B}_{s}
:=x_{s+P(j-1)+n-j+1,j-1}-x_{s+P(j)+n-j+1,j}\in (\QQ^{\ify})^*,
\]
where $x_{m,0}=0$ for $m\in\mathbb{Z}$, and $x_{m,i}=0$ for $m\in\mathbb{Z}_{\leq0}$ and $i\in I$.
\item
For $j\in[1,n-1]$ and $s\in\mathbb{Z}$, we set
\[
\fbox{$j$}^{\rm C}_{s}
:=x_{s+P(j),j}-x_{s+P(j-1)+1,j-1},
\quad
\fbox{$n$}^{\rm C}_{s}:=2x_{s+P(n),n}-x_{s+P(n-1)+1,n-1}\in (\QQ^{\ify})^*,\]
\[
\fbox{$\overline{n}$}^{\rm C}_{s}
:=x_{s+P(n-1)+1,n-1}-2x_{s+P(n)+1,n},
\quad
\fbox{$\overline{j}$}^{\rm C}_{s}
:=x_{s+P(j-1)+n-j+1,j-1}-x_{s+P(j)+n-j+1,j}\in (\QQ^{\ify})^*,
\]
\[
\fbox{$\overline{n+1}$}^{\rm C}_{s}
:=x_{s+P(n),n}
\in (\QQ^{\ify})^*,
\]
where $x_{m,0}=0$ for $m\in\mathbb{Z}$, and $x_{m,i}=0$ for $m\in\mathbb{Z}_{\leq0}$ and $i\in I$.
\item 
For $s\in\mathbb{Z}$, we set
\[
\fbox{$j$}^{\rm D}_{s}
:=x_{s+P(j),j}-x_{s+P(j-1)+1,j-1}\in (\QQ^{\ify})^*,\ \ (1\leq j\leq n-2,\ j=n),
\]
\[
\fbox{$n-1$}^{\rm D}_{s}
:=x_{s+P(n-1),n-1}+x_{s+P(n),n}-x_{s+P(n-2)+1,n-2}\in (\QQ^{\ify})^*,
\]
\[
 \fbox{$\overline{n}$}^{\rm D}_{s}
:=x_{s+P(n-1),n-1}-x_{s+P(n)+1,n}\in (\QQ^{\ify})^*,
\]
\[
\fbox{$\overline{n-1}$}^{\rm D}_{s}
:=x_{s+P(n-2)+1,n-2}-x_{s+P(n-1)+1,n-1}-x_{s+P(n)+1,n}\in (\QQ^{\ify})^*,
\]
\[
\fbox{$\overline{j}$}^{\rm D}_{s}
:=x_{s+P(j-1)+n-j,j-1}-x_{s+P(j)+n-j,j}\in (\QQ^{\ify})^*,\ \ (1\leq j\leq n-2),\]
\[
\fbox{$\overline{n+1}$}^{\rm D}_{s}
:=x_{s+P(n),n}\in (\QQ^{\ify})^*,
\]
where $x_{m,0}=0$ for $m\in\mathbb{Z}$, and $x_{m,i}=0$ for $m\in\mathbb{Z}_{\leq0}$ and $i\in I$.
\end{enumerate}
\end{defn}

\begin{lem}\label{box-lem}\cite{KaN}
\begin{enumerate}
\item In the case $\mathfrak{g}$ is of type A, the boxes $\fbox{$j$}^{\rm A}_{s}$ satisfy the following:
\begin{equation}\label{A-box}
\fbox{$j+1$}^{\rm A}_{s}=\fbox{$j$}^{\rm A}_{s} -\beta_{s+P(j),j} \qquad (1\leq j\leq n,\ s\geq 1-P(j)).
\end{equation}
\item In the case $\mathfrak{g}$ is of type B, the boxes $\fbox{$j$}^{\rm B}_{s}$
satisfy the following:
\begin{eqnarray}
\fbox{$j+1$}^{\rm B}_{s}&=&\fbox{$j$}^{\rm B}_{s} -\beta_{s+P(j),j}
 \qquad (1\leq j\leq n-1,\ s\geq 1-P(j)), \label{B-box1}\\
\fbox{$\overline{n}$}^{\rm B}_{s}&=&\fbox{$n$}^{\rm B}_{s} -\beta_{s+P(n),n}
 \qquad (s\geq 1-P(n)), \label{B-box2}\\
\fbox{$\overline{j-1}$}^{\rm B}_{s}&=&\fbox{$\overline{j}$}^{\rm B}_{s} -\beta_{s+P(j-1)+n-j+1,j-1}
\qquad (2\leq j\leq n,\ s\geq j-P(j-1)-n). \label{B-box3}
\end{eqnarray}
\item In the case $\mathfrak{g}$ is of type C, the boxes $\fbox{$j$}^{\rm C}_{s}$
satisfy the following:
\begin{eqnarray}
\fbox{$j+1$}^{\rm C}_{s}&=&\fbox{$j$}^{\rm C}_{s} -\beta_{s+P(j),j}\ \  (1\leq j\leq n-1,\ s\geq 1-P(j)), \label{C-box1}\\
\fbox{$\overline{n}$}^{\rm C}_{s}&=&\fbox{$n$}^{\rm C}_{s} -2\beta_{s+P(n),n}
\qquad (s\geq 1-P(n)), \label{C-box2}\\
\fbox{$\overline{j-1}$}^{\rm C}_{s}&=&\fbox{$\overline{j}$}^{\rm C}_{s} -\beta_{s+P(j-1)+n-j+1,j-1} 
\quad
 (2\leq j\leq n,\ s\geq j-P(j-1)-n),\ \ \ \ \ \label{C-box3}\\
\fbox{$\ovl{n+1}$}^{\rm C}_{l+1}+ \fbox{$\ovl{n}$}^{\rm C}_{l}
&=&\fbox{$\ovl{n+1}$}^{\rm C}_{l}-\beta_{l+P(n),n}
\qquad (l\geq 1-P(n)). \label{BC-pr3}
\end{eqnarray}

\item In the case $\mathfrak{g}$ is of type D, the boxes $\fbox{$j$}^{\rm D}_{s}$
satisfy the following:
\begin{eqnarray}
\fbox{$j+1$}^{\rm D}_{s}&=&\fbox{$j$}^{\rm D}_{s} -\beta_{s+P(j),j}
 \qquad (1\leq j\leq n-1,\ s\geq 1-P(j)), \label{D-box1}\\
\fbox{$\overline{n}$}^{\rm D}_{s}&=&\fbox{$n-1$}^{\rm D}_{s} -\beta_{s+P(n),n}
\qquad (s\geq 1-P(n)), \label{D-box2}\\
\fbox{$\overline{n-1}$}^{\rm D}_{s}&=&\fbox{$n$}^{\rm D}_{s} -\beta_{s+P(n),n} 
\qquad (s\geq 1-P(n)), \label{D-box3}\\
\fbox{$\overline{j-1}$}^{\rm D}_{s}&=&\fbox{$\overline{j}$}^{\rm D}_{s} -\beta_{s+P(j-1)+n-j,j-1} 
\ (2\leq j\leq n,\ s\geq 1+j-P(j-1)-n),
\qquad \quad \label{D-box4} \\
\fbox{$\overline{n+1}$}^{\rm D}_{l+2}
+\fbox{$\overline{n}$}^{\rm D}_{l+1}
+\fbox{$\overline{n-1}$}^{\rm D}_{l}
&=&\fbox{$\overline{n+1}$}^{\rm D}_{l} -\beta_{l+P(n),n}
\qquad (l\geq 1-P(n)). \label{D-box5}
\end{eqnarray}
\end{enumerate}
\end{lem}

\begin{defn}\label{tab-def}\cite{KaN}
\begin{enumerate}
\item
For ${\rm X}={\rm A}$, ${\rm B}$, ${\rm C}$ or ${\rm D}$ and $s\in\mathbb{Z}$, we set

\[
\begin{ytableau}
j_1 \\
j_2 \\
\vdots \\
\scriptstyle j_{k-1}\\
j_k
\end{ytableau}^{\rm X}_{s}
:=
\fbox{$
j_k$}_s^{X}
+\fbox{$j_{k-1}$}_{s+1}^{\rm X}+\cdots
+
\fbox{$
j_2 
$}^{\rm X}_{s+k-2}
+
\fbox{$
j_1
$}^{\rm X}_{s+k-1} \in (\QQ^{\ify})^*.
\]
\item
For ${\rm X}={\rm A}$, ${\rm B}$,
\[
{\rm Tab}_{{\rm X},\iota} := \{ \begin{ytableau}
j_1 \\
j_2 \\
\vdots \\
j_k
\end{ytableau}^{\rm X}_{s} | k\in I,\ j_i\in J_{\rm X},\ s\geq 1-P(k),\ (*)^{\rm X}_k \},
\]

$(*)^{\rm A}_k : 1\leq j_1<j_2<\cdots<j_k\leq n+1$, 

$(*)^{\rm B}_k : 
\begin{cases}
1\leq j_1<j_2<\cdots<j_k\leq \overline{1} & {\rm for}\ k<n,\\
1\leq j_1<j_2<\cdots<j_n\leq \overline{1},\ \ |j_l|\neq |j_m|\ (l\neq m) & {\rm for}\ k=n.
\end{cases}$



\[
{\rm Tab}_{{\rm C},\iota}:=
\{ 
\begin{ytableau}
j_1 \\
j_2 \\
\vdots \\
j_k
\end{ytableau}^{\rm C}_{s} 
| 
\begin{array}{l}
j_1\in J_{\rm C}\cup\{\ovl{n+1}\}, j_2,\cdots,j_k\in J_{\rm C}, \\
{\rm if}\ j_1\neq \ovl{n+1}\ {\rm then}\ k\in[1,n-1]\ {\rm and}\ 1\leq j_1<j_2<\cdots<j_k\leq \overline{1},\ s\geq 1-P(k), \\
{\rm if}\ j_1=\ovl{n+1}\ {\rm then}\ k\in[1,n+1],\ \overline{n}\leq j_2<\cdots<j_k\leq\overline{1},\ s\geq 1-P(n).\\
\end{array}
\}
\]



\[
{\rm Tab}_{{\rm D},\iota}:=
\{ 
\begin{ytableau}
j_1 \\
j_2 \\
\vdots \\
j_k
\end{ytableau}^{\rm D}_{s} 
| 
\begin{array}{l}
j_1\in J_{\rm D}\cup\{\ovl{n+1}\}, j_2,\cdots,j_k\in J_{\rm D},\\
{\rm if}\ j_1\neq\ovl{n+1}\ {\rm then}\ k\in[1,n-2]\ {\rm and}\ j_1\ngeq j_2\ngeq\cdots\ngeq j_k,\ s\geq 1-P(k),\\
{\rm if}\ j_1=\ovl{n+1}\ {\rm and}\ k\ {\rm is\ even} \ {\rm then}\ k\in[1,n+1],\ \overline{n}\leq j_2<\cdots<j_k\leq\overline{1},
\ s\geq 1-P(n-1), \\
{\rm if}\ j_1=\ovl{n+1}\ {\rm and}\ k\ {\rm is\ odd} \ {\rm then}\ k\in[1,n+1],\ \overline{n}\leq j_2<\cdots<j_k\leq\overline{1},
\ s\geq 1-P(n).
\end{array}
\}
\]

\end{enumerate}
\end{defn}

\begin{rem}

Similar notations to Definition \ref{box-def} and \ref{tab-def} (i) can be found
in \cite{Nj, NN}.

\end{rem}

\begin{thm}\label{thm2}\cite{KaN}
For ${\rm X}={\rm A}$, ${\rm B}$, ${\rm C}$ or ${\rm D}$,
we suppose that $\iota$ is adapted to the Cartan matrix of type ${\rm X}$.
Then
\[
\Xi_{\iota}={\rm Tab}_{{\rm X},\iota}.
\]
\end{thm}

\begin{thm}\label{thm1}\cite{KaN}
In the setting of Theorem \ref{thm2}, $\iota$ satisfies the positivity condition.
\end{thm}

\begin{cor}\label{cor1}\cite{KaN}
In the setting of Theorem \ref{thm2},
we have
\[
{\rm Im}(\Psi_{\iota}) =\{\textbf{x}\in\mathbb{Z}^{\infty}_{\iota} | \varphi(\textbf{x})\geq0,
\ {\rm for\ all}\ \varphi\in {\rm Tab}^{n}_{{\rm X},\iota},\ x_{m,i}=0\ {\rm for}\ m>n,\ i\in I \},
\]
where ${\rm Tab}^{n}_{{\rm X},\iota}:=
\{\begin{ytableau}
j_1 \\
\vdots \\
j_k
\end{ytableau}^{\rm X}_{s} \in{\rm Tab}_{{\rm X},\iota} | s\leq n \}$.
\end{cor}

\section{Tableaux descriptions of polyhedral realizations of $B(\lambda)$}

Just as in the previous section,
we take $\mathfrak{g}$ as a finite dimensional simple Lie algebra of type ${\rm A}_n$,
${\rm B}_n$, ${\rm C}_n$ or ${\rm D}_n$ and 
suppose $\iota=(\cdots,i_3,i_2,i_1)$ satisfies
$(\ref{seq-con})$ and is adapted to the Cartan matrix $A$ of $\mathfrak{g}$.
In this section we shall give a column tableaux description of polyhedral realizations of $B(\lambda)$.
We denote each tableau $\begin{ytableau}
j_1 \\
\vdots \\
j_l
\end{ytableau}^{\rm X}_{s}$ by $[j_1,\cdots,j_l]^{\rm X}_s$.

We consider the following two conditions on $k\in I$:
\begin{equation}\label{cond12}
(1)\ k<n\ {\rm and}\ \iota^{(k)}>\iota^{(k+1)},\qquad
(2)\ k>1\ {\rm and}\ \iota^{(k)}>\iota^{(k-1)}.
\end{equation}

\begin{defn}\label{tab-def-lam}
We fix a dominant integral weight $\lambda\in P_+$. For $k\in I$, we set 

\vspace{3mm}

${\rm Tab}_{{\rm A},\iota,k}[\lambda] :=$
\[
\begin{cases}
\{ -x_{1,k} + \lan h_k,\lambda\ran \} & {\rm if}\ {\rm (1),(2)\ do\ not\ hold}, \\
\{
\fbox{$t$}^{\rm A}_{1-P(k+1)} +\lan h_k,\lambda\ran  | t\in J_{{\rm A}},\ k+1\leq t \leq n+1 
\}
& {\rm if}\ {\rm only\ (1)\ holds}, \\
\{
[j_1,\cdots,j_{k-1},k+1,\cdots,n,n+1]^{\rm A}_{-P(k-1)-n+k}+\lan h_k,\lambda\ran | 
\begin{array}{l}
1\leq j_1<\cdots\\
<j_{k-1} \leq k,\ j_i\in J_{{\rm A}}
\end{array}
\}
& {\rm if}\ {\rm only\ (2)\ holds}, \\
\{
[j_1,\cdots,j_k]^{\rm A}_{-P(k-1)}+\lan h_k,\lambda\ran | j_i\in J_{{\rm A}},\ 1\leq j_1<\cdots<j_k \leq n+1,\ j_k>k
\}
& {\rm if}\ {\rm both\ (1)\ and\ (2)\ hold}. \\
\end{cases}
\]

For $k\in \{1,\cdots,n-1\}$, \vspace{-3mm}
\[
{\rm Tab}_{{\rm B},\iota,k}[\lambda] :=
\begin{cases}
\{ -x_{1,k} + \lan h_k,\lambda\ran \} &
 {\rm if}\ {\rm (1),(2)\ do\ not\ hold}, 
\\
\{
\fbox{$t$}^{\rm B}_{1-P(k+1)} +\langle h_k,\lambda\rangle  | t\in J_{{\rm B}},\ k+1\leq t \leq \overline{1} 
\}
& {\rm if}\ 
 {\rm only\ (1)\ holds},
\\
\{
\fbox{$t$}^{\rm B}_{-P(k-1)-n+k}+\langle h_k,\lambda\rangle | t\in J_{{\rm B}},\ \overline{k}\leq t\leq \overline{1} 
\}
& {\rm if}\ 
{\rm only\ (2)\ holds},
\\
\{
\begin{array}{l}
[j_1,\cdots,j_k]^{\rm B}_{-P(k-1)}\\+\langle h_k,\lambda\rangle
\end{array}
 |
\begin{array}{l}
j_1<\cdots<j_k ,\\ j_k>k,\ j_i\in J_{\rm B}
\end{array}
\}
& {\rm if}\ 
{\rm both\ (1)\ and\ (2)\ hold}, 
\end{cases}
\]
\vspace{-9mm}

\[
{\rm Tab}_{{\rm B},\iota,n}[\lambda] :=
\begin{cases}
\{ -x_{1,n} + \langle h_n,\lambda\rangle \} &
{\rm if}\ \ \iota^{(n)}<\iota^{(n-1)},
\\
\{
[j_1,\cdots,j_n]^{\rm B}_{-P(n-1)}+\langle h_n,\lambda\rangle |
\begin{array}{l}
j_1<\cdots<j_n,\\
  j_n>n,\ j_i\in J_{\rm B}\\
   |j_l|\neq |j_m|\ {\rm if\ }l\neq m 
\end{array}
\}
& {\rm if}
\ \iota^{(n)}>\iota^{(n-1)}.
\\
\end{cases}
\]

For $k\in \{1,\cdots,n-1\}$,
\[
{\rm Tab}_{{\rm C},\iota,k}[\lambda] :=
\begin{cases}
\{ -x_{1,k} + \lan h_k,\lambda\ran \} &
 {\rm if}\ {\rm (1),(2)\ do\ not\ hold},
\\
\{
\fbox{$t$}^{\rm C}_{1-P(k+1)} +\langle h_k,\lambda\rangle  | t\in J_{{\rm C}},\ k+1\leq t \leq \overline{1} 
\}
& {\rm if}\ 
 {\rm only\ (1)\ holds},
\\
\{
\fbox{$t$}^{\rm C}_{-P(k-1)-n+k}+\langle h_k,\lambda\rangle | t\in J_{{\rm C}},\ \overline{k}\leq t\leq \overline{1} 
\}
& {\rm if}\ 
 {\rm only\ (2)\ holds},
\\
\{
\begin{array}{l}
[j_1,\cdots,j_k]^{\rm C}_{-P(k-1)}\\+\langle h_k,\lambda\rangle 
\end{array}
|
\begin{array}{l}
j_1<\cdots<j_k ,\\ j_k>k,\ j_i\in J_{\rm C}
\end{array}
\}
& {\rm if}\ 
{\rm both\ (1)\ and\ (2)\ hold}, 
\\
\end{cases}
\]

\[
{\rm Tab}_{{\rm C},\iota,n}[\lambda] :=
\begin{cases}
\{ -x_{1,n} + \langle h_n,\lambda\rangle \} &
{\rm if}\ \ \iota^{(n)}<\iota^{(n-1)},
\\
\{
[\overline{n+1},j_2,\cdots,j_s]^{\rm C}_{-P(n-1)}+\langle h_n,\lambda\rangle |
\begin{array}{l}
2\leq s\leq n+1,\\ 
 \overline{n}\leq j_2<\cdots \\
 \cdots<j_s\leq \overline{1},\ j_i\in J_{\rm C}
 \end{array}
\}
& {\rm if}
\ \iota^{(n)}>\iota^{(n-1)}.
\\
\end{cases}
\]

For $k\in\{1,2,\cdots,n-3\}$, we set
\[
{\rm Tab}_{{\rm D},\iota,k}[\lambda] :=
\begin{cases}
\{ -x_{1,k} + \lan h_k,\lambda\ran\} & {\rm if}\ {\rm (1),(2)\ do\ not\ hold}, \\
\{
\fbox{$t$}^{\rm D}_{1-P(k+1)} +\lan h_k,\lambda\ran  | t\in J_{\rm D},\ k+1\leq t \leq \ovl{1} 
\}
& {\rm if}\ {\rm only\ (1)\ holds}, \\
\{
\fbox{$t$}^{\rm D}_{1-P(k-1)-n+k}+\lan h_k,\lambda\ran |t\in J_{\rm D},\ \ovl{k}\leq t\leq \ovl{1} 
\}
& {\rm if}\ {\rm only\ (2)\ holds} \\
\{
[j_1,\cdots,j_k]^{\rm D}_{-P(k-1)}+\lan h_k,\lambda\ran |
\begin{array}{l}
j_i\in J_{\rm D},\ k<j_k,\\ 
 j_1\ngeq \cdots\ngeq j_k
 \end{array}\}
& {\rm if}\ {\rm both\ (1)\ and\ (2)\ hold}. \\
\end{cases}
\]
For $t\in\{n-3,n-1,n\}$, we consider the following conditions:

\vspace{4mm}

\nd
$C_{t}$ : $\iota^{(t)}<\iota^{(n-2)}$ holds,\qquad $\ovl{C}_{t}$ : $\iota^{(t)}>\iota^{(n-2)}$ holds.

\begin{multline*}
{\rm Tab}_{{\rm D},\iota,n-2}[\lambda] :=\\
\begin{cases}
\{ -x_{1,n-2} + \lan h_{n-2},\lambda\ran \} & {\rm if}\ 
\ovl{C}_{n-3}, \ovl{C}_{n-1}, \ovl{C}_{n},
\\
\{
 -x_{1,n-2}+x_{1,n-1}+\lan h_{n-2},\lambda\ran,\ -x_{2,n-1}+\lan h_{n-2},\lambda\ran
\}
& {\rm if}\ 
\ovl{C}_{n-3}, C_{n-1}, \ovl{C}_{n},
\\
\{
 -x_{1,n-2}+x_{1,n}+\lan h_{n-2},\lambda\ran,\ -x_{2,n}+\lan h_{n-2},\lambda\ran
\}
& {\rm if}\ 
\ovl{C}_{n-3}, \ovl{C}_{n-1}, C_{n},
\\
\{
\fbox{$t$}^{\rm D}_{-1-P(n-3)}+\lan h_{n-2},\lambda\ran | t\in J_{\rm D},\ \ovl{n-2}\leq t\leq \ovl{1} 
\} & {\rm if}\ 
C_{n-3}, \ovl{C}_{n-1}, \ovl{C}_{n},
\\
\{
\fbox{$t$}^{\rm D}_{-P(n-2)}+\lan h_{n-2},\lambda\ran
|t\in J_{\rm D},\ n-1\leq t\leq \ovl{1},
\}
& {\rm if}\ 
\ovl{C}_{n-3}, C_{n-1}, C_{n},
\\
\left\{
[\ovl{n+1},j_2,\cdots,j_s]^{\rm D}_{-1-P(n-2)}
+\lan h_{n-2},\lambda\ran
\left| 
\begin{array}{l}
3\leq s\leq n+1,\ s\ {\rm is\ odd},\\
{\rm if\ }s=3\ {\rm then}\ j_3\geq\ovl{n-2},\\
\ovl{n}\leq j_2<\cdots<j_s\leq\ovl{1}\ {\rm in}\ J_{\rm D}
\end{array} \right.
\right\}
& {\rm if}\ 
C_{n-3}, C_{n-1}, \ovl{C}_{n},
\\
\left\{
[\ovl{n+1},j_2,\cdots,j_s]^{\rm D}_{-1-P(n-2)}
+\lan h_{n-2},\lambda\ran
\left| 
\begin{array}{l}
2\leq s\leq n+1,\ s\ {\rm is\ even},\\
{\rm if\ }s=2\ {\rm then}\ j_2\geq\ovl{n-2},\\
\ovl{n}\leq j_2<\cdots<j_s\leq\ovl{1}\ {\rm in}\ J_{\rm D}
\end{array} \right.
\right\}
& {\rm if}\ 
C_{n-3}, \ovl{C}_{n-1}, C_{n},
\\
\left\{
[j_1,\cdots,j_{n-2}]^{\rm D}_{-P(n-3)}+\lan h_{n-2},\lambda\ran
\left|
\begin{array}{l}
j_1,\cdots,j_{n-2} \in J_{\rm D},\\
j_1\ngeq \cdots \ngeq j_{n-2},\\
j_{n-2}\geq n-1,
\end{array} \right.
\right\}
& {\rm if}\ 
C_{n-3}, C_{n-1}, C_{n},
\end{cases}
\end{multline*}
\[
{\rm Tab}_{{\rm D},\iota,n-1}[\lambda] :=
\begin{cases}
\{ -x_{1,n-1} + \lan h_{n-1},\lambda\ran \} & {\rm if}\ 
C_{n-1}, 
\\
\left\{
[\ovl{n+1},j_2,\cdots,j_{s}]^{\rm D}_{-P(n-2)}+\lan h_{n-1},\lambda\ran
\left|
\begin{array}{l}
2\leq s\leq n+1,\ s\ {\rm is\ even},\\
j_2,\cdots,j_{s} \in J_{\rm D},\\
\ovl{n}\leq j_2< \cdots < j_{s}\leq \ovl{1},\\
{\rm if}\ s=2\ {\rm then}\ j_{2}\geq \ovl{n-1}
\end{array}\right.
\right\}
& {\rm if}\ 
\ \ovl{C}_{n-1},
\end{cases}
\]
\[
{\rm Tab}_{{\rm D},\iota,n}[\lambda] :=
\begin{cases}
\{ -x_{1,n} + \lan h_n,\lambda\ran \} & {\rm if}\ 
C_{n}, 
\\
\left\{
[\ovl{n+1},j_2,\cdots,j_{s}]^{\rm D}_{-P(n-2)}+\lan h_n,\lambda\ran
\left|
\begin{array}{l}
3\leq s\leq n+1,\ s\ {\rm is\ odd},\\
j_2,\cdots,j_{s} \in J_{\rm D},\\
\ovl{n}\leq j_2< \cdots < j_{s}\leq \ovl{1}
\end{array} \right.
\right\}
& {\rm if}\ 
\ovl{C}_{n}.
\end{cases}
\]

\end{defn}

For ${\rm X}={\rm A,B,C\ or\ D}$, we set
\[
{\rm Tab}_{{\rm X,\iota}}[\lambda] :=
\left(
\bigcup_{k\in I}
{\rm Tab}_{{\rm X},\iota,k}[\lambda]\right)\cup\{0\}.
\]

\begin{thm}\label{thm3}
Let $\mathfrak{g}$ be of type ${\rm X}={\rm A,B,C\ or\ D}$ and $\lambda\in P_+$.
If $\iota$ is adapted to the Cartan matrix of $\mathfrak{g}$ then the pair $(\iota,\lambda)$ satisfies the ample condition
and we have
\[
\Xi_{\iota}[\lambda]={\rm Tab}_{{\rm X,\iota}}[\lambda]\cup {\rm Tab}_{{\rm X,\iota}}.
\]
\end{thm}
\nd
The following corollary follows from Theorem \ref{Nthm1}, Corollary \ref{cor1} and Theorem \ref{thm3}.

\begin{cor}
In the setting of Theorem \ref{thm3}, we get
\[
{\rm Im}(\Psi^{(\lambda)}_{\iota})=
\{\textbf{x}\in\mathbb{Z}^{\infty}_{\iota}[\lambda] | \varphi(\textbf{x})\geq0,
\forall\varphi\in {\rm Tab}_{{\rm X,\iota}}[\lambda]\cup {\rm Tab}^{n}_{{\rm X},\iota},\ x_{m,i}=0\ (\forall i\in I,m>n) \}.
\]
\end{cor}

\begin{ex}\label{ex-3}

Let $\mathfrak{g}$ be the Lie algebra of type ${\rm A}_3$ and $\iota=(\cdots,3,1,2,3,1,2)$.
The sequence $\iota$ is adapted to the Cartan matrix of type ${\rm A}_3$. We obtain
$p_{2,1}=1$, $p_{3,2}=0$, $P(2)=P(3)=1$ and
\begin{eqnarray}
{\rm Tab}^3_{{\rm A},\iota}&=&\{\fbox{$j$}^{\rm A}_s| 1\leq s\leq 3,\ j\in[1,4] \}
\cup \left\{
[i,j]^{\rm A}_{s} 
\left| 
\begin{array}{l}
 0\leq s\leq 3,\\
 1\leq i<j\leq 4. 
\end{array}\right.
\right\}
\cup \left\{
[i,j,k]^{\rm A}_{s} 
\left| 
\begin{array}{l}
 0\leq s\leq 3,\\
 1\leq i<j<k\leq 4. 
\end{array}\right.
\right\}\nonumber \\
&=& 
\{x_{s,1},\ x_{s+1,2}-x_{s+1,1},\ x_{s+1,3}-x_{s+2,2},\ -x_{s+2,3} | 1\leq s\leq 3 \}\nonumber\\
&\cup&
\left\{x_{s+1,2},\ x_{s+1,3}-x_{s+2,2}+x_{s+1,1},\ x_{s+1,1}-x_{s+2,3},\ x_{s+1,3}-x_{s+2,1},\right.\nonumber\\
& &\left. x_{s+2,2}-x_{s+2,1}-x_{s+2,3},\ -x_{s+3,2} | 0\leq s\leq 3 \right\} \nonumber \\
&\cup&
\{x_{s+1,3},\ x_{s+2,2}-x_{s+2,3},\ x_{s+2,1}-x_{s+3,2},\ -x_{s+3,1} | 0\leq s\leq 3 \}.\nonumber
\end{eqnarray}
Since $\iota^{(1)}=2$, $\iota^{(2)}=1$, $\iota^{(3)}=3$, we get
\[
{\rm Tab}_{{\rm A},\iota,2}[\lambda]=\{-x_{1,2}+\lambda_2\},
\]
and
\begin{eqnarray*}
{\rm Tab}_{{\rm A},\iota,1}[\lambda]
&=&\{\fbox{2}^{\rm A}_0+ \lambda_1,\ \fbox{3}^{\rm A}_0+ \lambda_1,\ 
\fbox{4}^{\rm A}_0+ \lambda_1\}\\
&=&\{x_{1,2}-x_{1,1}+\lambda_1,\ x_{1,3}-x_{2,2}+\lambda_1,\ -x_{2,3}+\lambda_1\},
\end{eqnarray*}
\begin{eqnarray*}
{\rm Tab}_{{\rm A},\iota,3}[\lambda]
&=&\{[1,2,4]^{\rm A}_{-1}+ \lambda_3,\ [1,3,4]^{\rm A}_{-1}+ \lambda_3,\ 
[2,3,4]^{\rm A}_{-1}+ \lambda_3\}\\
&=&\{x_{1,2}-x_{1,3}+\lambda_3,\ x_{1,1}-x_{2,2}+\lambda_3,\ -x_{2,1}+\lambda_3\},
\end{eqnarray*}
where we put $\lambda_k:=\lan h_k,\lambda\ran$ ($k=1,2,3$). We get
\[
{\rm Im}(\Psi^{(\lambda)}_{\iota})
=\{\textbf{x}\in\mathbb{Z}^{\infty}_{\iota} | x_{m,i}=0\ (m\in\mathbb{Z}_{\geq4},i\in I),\ \varphi(\textbf{x})\geq0, \ \ \forall\varphi\in {\rm Tab}^3_{{\rm A},\iota}\cup
{\rm Tab}_{{\rm A},\iota,1}[\lambda] \cup
{\rm Tab}_{{\rm A},\iota,2}[\lambda]\cup
{\rm Tab}_{{\rm A},\iota,3}[\lambda] \}.
\]
For $\textbf{x}=(\cdots,x_{2,3},x_{2,1},x_{2,2},x_{1,3},x_{1,1},x_{1,2})\in {\rm Im}(\Psi_{\iota})$,
combining inequalities $x_{s,1}\geq0$ $(1\leq s\leq3)$, $-x_{s+3,1}\geq0$ $(0\leq s\leq3)$ in ${\rm Tab}_{{\rm A},\iota}$, 
we obtain
$x_{3,1}=0$.
Similarly, by $x_{s+1,2}\geq0$, $-x_{s+3,2}\geq0$ $(0\leq s\leq3)$,
we get $x_{3,2}=0$. We also get $x_{3,3}=0$.
Hence, simplifying the inequalities, we obtain
\[
{\rm Im}(\Psi^{(\lambda)}_{\iota})=
\left\{\textbf{x}\in \mathbb{Z}^{\infty}_{\iota} \left|
\begin{array}{l}
 x_{s,1}=x_{s,2}=x_{s,3}=0\ {\rm for}\ s\in\mathbb{Z}_{\geq3},\ 
\ x_{2,2}- x_{2,1}\geq x_{2,3}\geq0,\\
x_{1,3}- x_{2,2}+x_{1,1}\geq0,\ x_{1,1}\geq x_{2,3}\geq0,\ x_{1,3}\geq x_{2,1}\geq0,\ x_{1,2}\geq0,\\
\lambda_2\geq x_{1,2},\ \lambda_1\geq x_{1,1}-x_{1,2},\ \lambda_1\geq x_{2,2}-x_{1,3},\ \lambda_1\geq x_{2,3},\\
\lambda_3\geq x_{1,3}-x_{1,2},\ \lambda_3\geq x_{2,2}-x_{1,1},\ \lambda_3\geq x_{2,1}
\end{array}\right.
\right\}.
\]

\end{ex}

\begin{ex}\label{ex-4}

Let $\mathfrak{g}$ be the Lie algebra of type ${\rm C}_3$ and $\iota=(\cdots,3,1,2,3,1,2)$.
The sequence $\iota$ is adapted to the Cartan matrix of type ${\rm C}_3$. We get
$p_{2,1}=1$, $p_{3,2}=0$, $P(2)=P(3)=1$ and
\begin{eqnarray}
{\rm Tab}^{3}_{{\rm C},\iota}&=&\{\fbox{$j$}^{\rm C}_s| 1\leq s\leq 3,\ 1\leq j\leq \ovl{1}\}
\cup\{
\begin{ytableau}
i \\
j 
\end{ytableau}^{\rm C}_{s} 
| 
\begin{array}{l}
 0\leq s\leq 3,\\
 1\leq i<j\leq \ovl{1}. 
\end{array}
\}
\cup\{
\begin{ytableau}
\ovl{4} \\
j_2 \\
\vdots\\
j_k 
\end{ytableau}^{\rm C}_{s} 
| 
\begin{array}{l}
 0\leq s\leq 3,\\
 k\in[1,4],\ 
 \ovl{3}\leq j_2<\cdots<j_k\leq \ovl{1}. 
\end{array}
\}\nonumber \\
&=& 
\{x_{s,1},\ x_{s+1,2}-x_{s+1,1},\ 2x_{s+1,3}-x_{s+2,2},\ x_{s+2,2}-2x_{s+2,3},\ x_{s+2,1}-x_{s+3,2},\ -x_{s+3,1}| 1\leq s\leq 3 \}\nonumber\\
&\cup&
\left\{x_{s+1,2},\ 2x_{s+1,3}-x_{s+2,2}+x_{s+1,1},\ x_{s+1,1}+x_{s+2,2}-2x_{s+2,3},\ x_{s+1,1}+x_{s+2,1}-x_{s+3,2},
\right.\nonumber\\
& &\left. x_{s+1,1}-x_{s+3,1},\ 2x_{s+1,3}-x_{s+2,1},\ 2x_{s+2,2}-x_{s+2,1}-2x_{s+2,3},\ x_{s+2,2}-x_{s+3,2},\ 
x_{s+2,2}-x_{s+2,1} -x_{s+3,1}, \right. \nonumber \\
& & \left. 2x_{s+2,3}-2x_{s+3,2} + x_{s+2,1},\ 2x_{s+2,3}-x_{s+3,2} - x_{s+3,1},\ x_{s+2,1}-2x_{s+3,3},\right. \label{ex-2} \\
& & \left. x_{s+3,2}-2x_{s+3,3}-x_{s+3,1},\ -x_{s+4,2} | 0\leq s\leq 3 \right\} \nonumber \\
&\cup&
\left\{x_{s+1,3},\ x_{s+2,2}-x_{s+2,3},\ x_{s+2,3}+x_{s+2,1}-x_{s+3,2},\ x_{s+2,3}-x_{s+3,1},\ x_{s+2,1}-x_{s+3,3}, \right. \nonumber\\
& & \left. x_{s+3,2}-x_{s+3,1}-x_{s+3,3},\ x_{s+3,3}-x_{s+4,2},\ -x_{s+4,3} | 0\leq s\leq 3 \right\}.\nonumber
\end{eqnarray}
Because $\iota^{(1)}=2$, $\iota^{(2)}=1$, $\iota^{(3)}=3$, we get
\[
{\rm Tab}_{{\rm C},\iota,2}[\lambda]=\{-x_{1,2}+\lambda_2\},
\]
and
\begin{eqnarray*}
{\rm Tab}_{{\rm C},\iota,1}[\lambda]
&=&\{\fbox{2}^{\rm C}_0+ \lambda_1,\ \fbox{3}^{\rm C}_0+ \lambda_1,\ 
\fbox{$\ovl{3}$}^{\rm C}_0+ \lambda_1,\ \fbox{$\ovl{2}$}^{\rm C}_0+ \lambda_1,\ 
\fbox{$\ovl{1}$}^{\rm C}_0+ \lambda_1 \}\\
&=&\{x_{1,2}-x_{1,1}+\lambda_1,\ 2x_{1,3}-x_{2,2}+\lambda_1,\ x_{2,2}-2x_{2,3}+\lambda_1,\ 
 x_{2,1}-x_{3,2}+\lambda_1,\ -x_{3,1}+\lambda_1\},
\end{eqnarray*}
\begin{eqnarray*}
{\rm Tab}_{{\rm C},\iota,3}[\lambda]
&=&\{[\ovl{4},\ovl{3}]^{\rm C}_{-1}+ \lambda_3,\ [\ovl{4},\ovl{2}]^{\rm C}_{-1}+ \lambda_3,\ 
[\ovl{4},\ovl{1}]^{\rm C}_{-1}+ \lambda_3,\ 
[\ovl{4},\ovl{3},\ovl{2}]^{\rm C}_{-1}+ \lambda_3,\\ 
& &[\ovl{4},\ovl{3},\ovl{1}]^{\rm C}_{-1}+ \lambda_3,
[\ovl{4},\ovl{2},\ovl{1}]^{\rm C}_{-1}+ \lambda_3,
[\ovl{4},\ovl{3},\ovl{2},\ovl{1}]^{\rm C}_{-1}+ \lambda_3
\}\\
&=&\{x_{1,2}-x_{1,3}+\lambda_3,\ x_{1,3}+x_{1,1}-x_{2,2}+\lambda_3,\ x_{1,3}-x_{2,1}+\lambda_3
,\ x_{1,1}-x_{2,3}+\lambda_3, \\
& & x_{2,2}-x_{2,1}-x_{2,3}+\lambda_3,\ x_{2,3}-x_{3,2}+\lambda_3,
\ -x_{3,3}+\lambda_3\},
\end{eqnarray*}
where we put $\lambda_k:=\lan h_k,\lambda\ran$ ($k=1,2,3$).
Thus. we get
\[
{\rm Im}(\Psi^{(\lambda)}_{\iota})
=\{\textbf{x}\in\mathbb{Z}^{\infty}_{\iota} |x_{m,i}=0\ (m\in\mathbb{Z}_{\geq4},i\in I),\  \varphi(\textbf{x})\geq0, \ \ \forall\varphi\in {\rm Tab}^3_{{\rm C},\iota}\cup
{\rm Tab}_{{\rm C},\iota,1}[\lambda] \cup
{\rm Tab}_{{\rm C},\iota,2}[\lambda]\cup
{\rm Tab}_{{\rm C},\iota,3}[\lambda] \}.
\]
Simplifying the inequalities, we obtain
\[
{\rm Im}(\Psi^{(\lambda)}_{\iota})\\
=
\left\{\textbf{x}\in \mathbb{Z}^{\infty}_{\iota} \left|
\begin{array}{l} x_{s,1}=x_{s,2}=x_{s,3}=0\ {\rm for}\ s\in\mathbb{Z}_{\geq4},\ 
\ x_{2,2}\geq x_{2,1}\geq0, \ 
2x_{2,3}\geq x_{3,2}\geq x_{3,1}\geq 0, \\
 x_{3,2}-2x_{3,3}\geq0,\ 2x_{1,3}-x_{2,2}+x_{1,1}\geq0, \\ 
 x_{1,1}+x_{2,2}-2x_{2,3}\geq0,\ 
 x_{1,1}+x_{2,1}-x_{3,2}\geq0,\ x_{1,1}\geq x_{3,1}\geq0,\ 
2x_{1,3}-x_{2,1}\geq0, \\
 2x_{2,2}-x_{2,1}-2x_{2,3}\geq0,\ x_{2,2}-x_{3,2}\geq0,\ x_{2,2}-x_{2,1}-x_{3,1}\geq0,\ 2x_{2,3}-2x_{3,2}+x_{2,1}\geq0, 
 \\
 2x_{2,3}- x_{3,2}-x_{3,1}\geq0,\ x_{2,1}- 2x_{3,3}\geq0,\ x_{3,2}- 2x_{3,3}-x_{3,1}\geq0, \\
 x_{2,2}-x_{2,3}\geq0,\ x_{2,3}+x_{2,1}-x_{3,2}\geq0,\ x_{2,3}\geq x_{3,1}\geq0,\ x_{2,1}\geq x_{3,3}\geq0,\\
x_{3,2}-x_{3,1}-x_{3,3}\geq0,\ 
x_{1,2}\geq0,\ x_{1,3}\geq0,\\
\lambda_2\geq x_{1,2},\ \lambda_1\geq x_{1,1}-x_{1,2},\ \lambda_1\geq x_{2,2}-2x_{1,3},\ 
\lambda_1\geq 2x_{2,3}-x_{2,2},\ \lambda_1\geq x_{3,2}-x_{2,1},\\ \lambda_1\geq x_{3,1},\ 
\lambda_3\geq x_{1,3}-x_{1,2},\ 
\lambda_3\geq x_{2,2}-x_{1,1}-x_{1,3},\ 
\lambda_3\geq x_{2,1}-x_{1,3},\ 
\lambda_3\geq x_{2,3}-x_{1,1},\\
\lambda_3\geq x_{2,1}+x_{2,3}-x_{2,2},\ 
\lambda_3\geq x_{3,2}-x_{2,3},\ 
\lambda_3\geq x_{3,3}\ 
\end{array}\right. \right\}.
\]


\end{ex}

\begin{thm}\label{thm4}
Let $\mathfrak{g}$ be of type ${\rm X}={\rm A,B,C\ or\ D}$.
If $\iota$ is adapted to the Cartan matrix of $\mathfrak{g}$ then $\iota$
satisfies the strict positivity condition. In particular, for $i\in I$ and $x\in \Sigma_{\iota}$, we get
\[
\varepsilon^*_i(x) = {\rm max}\{-\varphi(x) | \varphi\in{\rm Tab}_{{\rm X},\iota,i}[0]\}.
\]
\end{thm}

\begin{ex}

Let $\mathfrak{g}$ be the Lie algebra of type ${\rm A}_3$, $\iota=(\cdots,3,1,2,3,1,2)$ and
\[
b:=(\cdots,b_{3,3},b_{3,1},b_{3,2},b_{2,3},b_{2,1},b_{2,2},b_{1,3},b_{1,1},b_{1,2})=(\cdots,0,0,0,2,1,3,1,2,1)\in \Sigma_{\iota}={\rm Im}(\Psi_{\iota}).
\]
Following 2.4 of \cite{NZ},
we can calculate the action of $*$ on $b$ as
\[
b^*=\tilde{f}_2\tilde{f}^2_1\tilde{f}_3\tilde{f}^3_2\tilde{f}_1\tilde{f}^2_3 u_{\infty}
=(\cdots,0,0,1,2,3,2,1,1)
\]
and
\[
\varepsilon_1^*(b)=\varepsilon_1(b^*)=2,\ \varepsilon_2^*(b)=\varepsilon_2(b^*)=1,\ \ 
\varepsilon_3^*(b)=\varepsilon_3(b^*)=1.
\]
On the other hand, we have seen in Example \ref{ex-3} that
\[
{\rm max}\{-\varphi(b) | \varphi\in{\rm Tab}_{{\rm A},\iota,1}[0]\}=
{\rm max}\{0,b_{1,1}-b_{1,2}, b_{2,2}-b_{1,3}, b_{2,3} \}={\rm max}\{0,1,2\}=2,
\]
\[
{\rm max}\{-\varphi(b) | \varphi\in{\rm Tab}_{{\rm A},\iota,2}[0]\}={\rm max}\{0,1\}=1,
\]
\[
{\rm max}\{-\varphi(b) | \varphi\in{\rm Tab}_{{\rm A},\iota,3}[0]\}=
{\rm max}\{0,b_{1,3}-b_{1,2}, b_{2,2}-b_{1,1}, b_{2,1} \}={\rm max}\{0,1\}=1.
\]
Thus, it holds $\varepsilon^*_i(b) = {\rm max}\{-\varphi(b) | \varphi\in{\rm Tab}_{{\rm A},\iota,i}[0]\}$.

\end{ex}

\begin{ex}

Let $\mathfrak{g}$ be the Lie algebra of type ${\rm C}_3$, $\iota=(\cdots,3,1,2,3,1,2)$ and
\[
b:=(\cdots,b_{3,3},b_{3,1},b_{3,2},b_{2,3},b_{2,1},b_{2,2},b_{1,3},b_{1,1},b_{1,2})=(\cdots,0,0,2,4,2,7,2,3,1)\in \Sigma_{\iota}={\rm Im}(\Psi_{\iota}).
\]
Calculating $b^*$ as
\[
b^*=\tilde{f}_2\tilde{f}^3_1\tilde{f}^2_3\tilde{f}^7_2\tilde{f}^2_1\tilde{f}^4_3 \tilde{f}^2_2 u_{\infty}
=(\cdots,0,0,0,2,3,8,4,2,2),
\]
we see that
\[
\varepsilon^*_1(b)=3,\ \ 
\varepsilon^*_2(b)=1,\ \ 
\varepsilon^*_3(b)=2. 
\]
By Example \ref{ex-4},
\[
{\rm max}\{-\varphi(b) | \varphi\in{\rm Tab}_{{\rm C},\iota,1}[0]\}=
{\rm max}\{0,b_{1,1}-b_{1,2}, 
b_{2,2}-2b_{1,3},
2b_{2,3}-b_{2,2}, b_{3,2}-b_{2,1},b_{3,1}
 \}={\rm max}\{0,2,3,1\}=3,
\]
\[
{\rm max}\{-\varphi(b) | \varphi\in{\rm Tab}_{{\rm C},\iota,2}[0]\}={\rm max}\{0,1\}=1,
\]
\begin{eqnarray*}
& &{\rm max}\{-\varphi(b) | \varphi\in{\rm Tab}_{{\rm C},\iota,3}[0]\}\\
&=&{\rm max}\{0,b_{1,3}-b_{1,2}, b_{2,2}-b_{1,1}-b_{1,3}, b_{2,1}-b_{1,3}, b_{2,3}-b_{1,1},
b_{2,1}+b_{2,3}-b_{2,2}, b_{3,2}-b_{2,3}, b_{3,3} \}\\
&=&{\rm max}\{0,1,2,-1,-2\}=2.
\end{eqnarray*}
Hence, we have $\varepsilon^*_i(b) = {\rm max}\{-\varphi(b) | \varphi\in{\rm Tab}_{{\rm C},\iota,i}[0]\}$.

\end{ex}

\section{Actions of operators $\what{S}_{m,j}$}\label{Stabsect}

As in the previous section, we denote each tableau $\begin{ytableau}
j_1 \\
\vdots \\
j_k
\end{ytableau}^{\rm X}_{s}$ by $[j_1,\cdots,j_k]_s^{\rm X}$.
When we see the condition $j_l\neq t$ with $t\in\{1,2,\cdots,n,n+1,\ovl{n},\cdots,\ovl{1}\}$
for $[j_1,\cdots,j_k]_s^{\rm X}$ it means $j_l\neq t$ with $l\in[1,k]$ or $l>k$ or $l<1$.

\subsection{Actions of operators $\what{S}_{m,j}$ for type A}

\begin{prop}\label{closednessA}

We suppose that $j_1<\cdots<j_k$ ($j_1,\ \cdots,\ j_k\in J_{\rm A}$) and
put $T:=[j_1,\cdots,j_k]_s^{\rm A}$ with $s\in\mathbb{Z}$. For $m\in\mathbb{Z}_{\geq1}$ and $j\in I$,
\begin{eqnarray*}
\hspace{-15mm}& &\what{S}_{m,j}T=\\
\hspace{-15mm}& &\begin{cases}
[j_1,\cdots,j_{i-1},j+1,j_{i+1},\cdots,j_k]_s^{\rm A} & {\rm if}\ j_i=j,\ j_{i+1}\neq j+1,\ m=s+k-i+P(j)
\ {\rm for\ some}\ i\in[1,k],\\
[j_1,\cdots,j_{i-1},j,j_{i+1},\cdots,j_k]_s^{\rm A} & {\rm if}\ j_i=j+1,\ j_{i-1}\neq j,\ m=s+k-i+1+P(j)>1\ {\rm for\ some}\ i\in[1,k], \\
T+\beta^{(-)}_{1,j} & {\rm if}\ j_i=j+1,\ j_{i-1}\neq j,\ m=s+k-i+1+P(j)=1\ {\rm for\ some}\ i\in[1,k], \\
T & {\rm otherwise}.
\end{cases}
\end{eqnarray*}
\end{prop}

\nd
{\it Proof.}

\nd
We see that
\begin{equation}\label{A-pr01}
 T=
[j_1,\cdots,j_k]^{\rm A}_s 
=\sum^{k}_{i=1} \fbox{$j_i$}^{\rm A}_{s+k-i} 
=\sum^{k}_{i=1} (x_{s+k-i+P(j_i),j_i}-x_{s+k-i+1+P(j_i-1),j_{i}-1}),
\end{equation}
where we set $x_{t,l}:=0$ for $t\in\mathbb{Z}_{\leq0}$ and $l\in I$.
Note that since
\begin{eqnarray*}
& &
\fbox{$j_i$}^{\rm A}_{s+k-i}
+
\fbox{$j_{i+1}$}^{\rm A}_{s+k-i-1}\\ 
&=&
x_{s+k-i+P(j_i),j_i}-x_{s+k-i+1+P(j_i-1),j_i-1}
+x_{s+k-i-1+P(j_{i+1}),j_{i+1}}-x_{s+k-i+P(j_{i+1}-1),j_{i+1}-1}, 
\end{eqnarray*}
if $j_{i+1}=j_{i}+1$ then we get
\begin{equation}\label{A-pr02}
\fbox{$j_i$}^{\rm A}_{s+k-i}
+
\fbox{$j_{i}+1$}^{\rm A}_{s+k-i-1}=
x_{s+k-i-1+P(j_{i}+1),j_{i}+1}-x_{s+k-i+1+P(j_i-1),j_i-1}.
\end{equation}
It follows from (\ref{A-pr01}) and (\ref{A-pr02}) that for $m\in\mathbb{Z}_{\geq1}$ and $j\in I$,
$x_{m,j}$ has non-zero coefficient in $T$ if and only if
the pair $(m,j)$ belongs to
\begin{multline*}
\{ (s+k-i+P(j_i),j_i) | i=1,2,\cdots,k,\ \ j_{i+1}> j_i+1\} \\
\cup \{(s+k-i+1+P(j_i-1),j_{i}-1) | i=1,2,\cdots,k,\ \ j_i-1> j_{i-1}\},
\end{multline*}
where we set $j_{k+1}=n+2$ and $j_0=0$. 

If $(m,j)=(s+k-i+P(j_i),j_i)$ with $j_{i+1}> j_i+1$ then $x_{m,j}$ has coefficient $1$ in $T$
and by Lemma \ref{box-lem} (\ref{A-box}) and the definition of $\what{S}_{m,j}$,
\begin{eqnarray*}
\what{S}_{m,j} T
&=&T-\beta_{m,j}
= [j_1,\cdots,j_{i-1},j,j_{i+1},\cdots,j_k]^{\rm A}_s -\beta_{m,j}\\
&=&\fbox{$j_1$}^{\rm A}_{s+k-1}+\cdots+
\fbox{$j_{i-1}$}^{\rm A}_{s+k-i+1}+
\fbox{$j$}^{\rm A}_{s+k-i}+
\fbox{$j_{i+1}$}^{\rm A}_{s+k-i-1}+
\cdots +\fbox{$j_k$}^{\rm A}_{s} -\beta_{m,j}\\
&=&
\fbox{$j_1$}^{\rm A}_{s+k-1}+\cdots+
\fbox{$j_{i-1}$}^{\rm A}_{s+k-i+1}+
\fbox{$j+1$}^{\rm A}_{s+k-i}+
\fbox{$j_{i+1}$}^{\rm A}_{s+k-i-1}+
\cdots +\fbox{$j_k$}^{\rm A}_{s} \\
&=&
[j_1,\cdots,j_{i-1},j+1,j_{i+1},\cdots,j_k]^{\rm A}_s.
\end{eqnarray*}
If $(m,j)=(s+k-i+1+P(j_i-1),j_{i}-1)$ with $j_i-1> j_{i-1}$ then $x_{m,j}$ has coefficient $-1$ in $T$.
If $m>1$ then
\begin{eqnarray*}
\what{S}_{m,j} T&=&\what{S}_{m,j}
[j_1,\cdots,j_i,\cdots,j_k]^{\rm A}_s 
=
[j_1,\cdots,j_i,\cdots,j_k]^{\rm A}_s +\beta_{m-1,j}\\
&=&  
\fbox{$j_1$}^{\rm A}_{s+k-1}+\cdots+\fbox{$j_i$}^{\rm A}_{s+k-i}+\cdots +\fbox{$j_k$}^{\rm A}_{s}+\beta_{m-1,j}\\
&=&  
\fbox{$j_1$}^{\rm A}_{s+k-1}+\cdots+\fbox{$j+1$}^{\rm A}_{s+k-i}+\cdots +\fbox{$j_k$}^{\rm A}_{s}+\beta_{m-1,j}\\
&=&\fbox{$j_1$}^{\rm A}_{s+k-1}+\cdots+\fbox{$j$}^{\rm A}_{s+k-i}+\cdots +\fbox{$j_k$}^{\rm A}_{s}
=[j_1,\cdots,j,\cdots,j_k]^{\rm A}_s.
\end{eqnarray*}
If $m=1$ then the definition of $\what{S}_{m,j}$ means $\what{S}_{m,j}T=T+\beta^{(-)}_{1,j}$.
The definition of $\what{S}_{m,j}$ also means that if $x_{m,j}$ is not a summand of $T$ then $\what{S}_{m,j}T=T$.
Consequently, we get our claim. \qed

\subsection{Actions of operators $\what{S}_{m,j}$ for type B, C}

In this subsection, we consider type B, C cases.

\begin{prop}\label{closednessBC}
For each $T=[j_1,\cdots,j_k]_{-P(k-1)}^{\rm X}+\lan h_k,\lambda\ran \in {\rm Tab}_{{\rm X},\iota,k}\setminus \{\lambda^{(k)}\}$ 
$({\rm X}={\rm B}$ or ${\rm C},\ k\in I)$,
$j\in I$ and $m\in\mathbb{Z}_{\geq1}$, we consider the following conditions for the triple $(T,j,m)$:
\begin{enumerate}
\item[(1)] $j<n$ and there exists $i\in [1,k]$ such that $j_i=j$, $j_{i+1}\neq j+1$ and $m=-P(k-1)+k-i+P(j)$,
\item[(2)] $j<n$ and there exists $i'\in [1,k]$ such that $j_{i'}=\overline{j+1}$, $j_{i'+1}\neq \overline{j}$
and $m=-P(k-1)+k-i'+n-j+P(j)$,
\item[(3)] $j=n$ and there exists $i\in [1,k]$ such that $j_i=n$, $j_{i+1}\neq \overline{n}$ and $m=-P(k-1)+k-i+P(n)$,
\item[(4)] $j<n$ and there exists $i\in [1,k]$ such that $j_{i-1}\neq j$, $j_{i}= j+1$ and $m=1-P(k-1)+k-i+P(j)$,
\item[(5)] $j<n$ and there exists $i'\in [1,k]$ such that $j_{i'-1}\neq\overline{j+1}$, $j_{i'}= \overline{j}$
and $m=1-P(k-1)+k-i'+n-j+P(j)$,
\item[(6)] $j=n$ and there exists $i\in [1,k]$ such that $j_{i-1}\neq n$, $j_{i}= \overline{n}$ and $m=1-P(k-1)+k-i+P(n)$.
\end{enumerate}
\begin{enumerate}
\item
We suppose $j_1\neq \ovl{n+1}$ and $k\in[2,n-1]$ satisfies the both conditions (1) and (2).
Then we have
\begin{eqnarray*}
\hspace{-15mm}& &\what{S}_{m,j}T= \\
\hspace{-15mm}& &\begin{cases}
[j_1,\cdots,j_{i-1},j+1,j_{i+1},\cdots,j_k]_{-P(k-1)}^{\rm X}+\lan h_k,\lambda\ran & {\rm if\ (1)\ holds\ and\ (2),\ (5)\ do\ not\ hold}, \\
[j_1,\cdots,j_{i'-1},\ovl{j},j_{i'+1},\cdots,j_k]_{-P(k-1)}^{\rm X}+\lan h_k,\lambda\ran & {\rm if\ (2)\ holds\ and\ (1),\ (4)\ do\ not\ hold}, \\
[j_1,\cdots,j_{i-1},j+1,j_{i+1},\cdots,j_{i'-1},\ovl{j},j_{i'+1},\cdots,j_k]_{-P(k-1)}^{\rm X}+\lan h_k,\lambda\ran & {\rm if\ (1)\ and\ (2)\ hold}, \\
[j_1,\cdots,j_{i-1},\ovl{n},j_{i+1}\cdots,,j_k]_{-P(k-1)}^{\rm X}+\lan h_k,\lambda\ran & {\rm if\ (3)\ holds}, \\
[j_1,\cdots,j_{i-1},j,j_{i+1},\cdots,j_k]_{-P(k-1)}^{\rm X}+\lan h_k,\lambda\ran & {\rm if\ (4)\ holds\ and\ (2),\ (5)\ do\ not\ hold}, \ \ \\
[j_1,\cdots,j_{i'-1},\ovl{j+1},j_{i'+1},\cdots,j_k]_{-P(k-1)}^{\rm X}+\lan h_k,\lambda\ran & {\rm if\ (5)\ holds\ and\ (1),\ (4)\ do\ not\ hold}, \\
[j_1,\cdots,j_{i-1},j,j_{i+1},\cdots,j_{i'-1},\ovl{j+1},j_{i'+1},\cdots,j_k]_{-P(k-1)}^{\rm X}+\lan h_k,\lambda\ran & {\rm if\ (4)\ and\ (5)\ hold}, \\
[j_1,\cdots,j_{i-1},n,j_{i+1}\cdots,,j_k]_{-P(k-1)}^{\rm X}+\lan h_k,\lambda\ran & {\rm if\ (6)\ holds}, \\
T & {\rm otherwise}.
\end{cases}
\end{eqnarray*}
\item We set $k=n$ and suppose $\iota^{(n)}>\iota^{(n-1)}$.
For each $T=[j_1,\cdots,j_n]_{-P(n-1)}^{\rm B}+\lan h_n,\lambda\ran\in {\rm Tab}_{{\rm B},\iota,n}\setminus \{\lambda^{(n)}\}$,
 $m\in\mathbb{Z}_{\geq1}$ and $j\in I$, we have
\[ \what{S}_{m,j}T=
\begin{cases}
[j_1,\cdots,j_{i-1},j+1,j_{i+1},\cdots,j_{i'-1},\ovl{j},j_{i'+1},\cdots,j_n]_{-P(n-1)}^{\rm B}+\lan h_n,\lambda\ran & {\rm if\ (1)\ and\ (2)\ hold}, \\
[j_1,\cdots,j_{i-1},\ovl{n},j_{i+1}\cdots,,j_n]_{-P(n-1)}^{\rm B}+\lan h_n,\lambda\ran & {\rm if\ (3)\ holds}, \\
[j_1,\cdots,j_{i-1},j,j_{i+1},\cdots,j_{i'-1},\ovl{j+1},j_{i'+1},\cdots,j_n]_{-P(n-1)}^{\rm B}+\lan h_n,\lambda\ran & {\rm if\ (4)\ and\ (5)\ hold}, \\
[j_1,\cdots,j_{i-1},n,j_{i+1}\cdots,,j_n]_{-P(n-1)}^{\rm B}+\lan h_n,\lambda\ran & {\rm if\ (6)\ holds}, \\
T & {\rm otherwise}.
\end{cases}
\]
\item We suppose $\iota^{(n)}>\iota^{(n-1)}$.
For each $T=[\ovl{n+1},j_2,j_3,\cdots,j_k]_{-P(n-1)}^{\rm C}+\lan h_n,\lambda\ran\in {\rm Tab}_{{\rm C},\iota,n}
\setminus \{\lambda^{(n)}\}$ with $k\in [2,n+1]$,
$m\in\mathbb{Z}_{\geq1}$ and $j\in I$, we have
\[ \what{S}_{m,j}T=
\begin{cases}
[\ovl{n+1},j_2,\cdots,j_{i'-1},\ovl{j},j_{i'+1},\cdots,j_k]_{-P(n-1)}^{\rm C}+\lan h_n,\lambda\ran & {\rm if\ (2)'\ holds}, \\
[\ovl{n+1},j_2,\cdots,j_{i'-1},\ovl{j+1},j_{i'+1},\cdots,j_k]_{-P(n-1)}^{\rm C}+\lan h_n,\lambda\ran & {\rm if\ (5)'\ holds}, \\
[\ovl{n+1},j_3,\cdots,j_k]_{-P(n-1)}^{\rm C}+\lan h_n,\lambda\ran & {\rm if\ (6)'\ holds}, \\
[\ovl{n+1},\ovl{n},j_2,\cdots,j_k]_{-P(n-1)}^{\rm C}+\lan h_n,\lambda\ran & {\rm if\ (7)\ holds}, \\
T & {\rm otherwise},
\end{cases}
\]
where the conditions (2)', (5)', (6)', (7) are as follows:
\begin{enumerate}
\item[(2)'] $j<n$ and there exists $i'\in [1,k]$ such that $j_{i'}=\overline{j+1}$, $j_{i'+1}\neq \overline{j}$
and $m=-P(n-1)+k-i'+n-j+P(j)$,
\item[(5)'] $j<n$ and there exists $i'\in [1,k]$ such that $j_{i'-1}\neq\overline{j+1}$, $j_{i'}= \overline{j}$
and $m=1-P(n-1)+k-i'+n-j+P(j)$,
\item[(6)'] $j=n$, $j_{2}= \overline{n}$ and $m=-P(n-1)-1+k+P(n)$,
\item[(7)] $j=n$, $j_2\neq\ovl{n}$ and $m=-P(n-1)-1+k+P(n)$.
\end{enumerate}
\end{enumerate}

\end{prop}

\nd
{\it Proof.}

\nd
(i) In this setting, we get
\[
{\rm Tab}_{{\rm X},\iota,k}[\lambda]=
\{
[j_1,\cdots,j_k]^{\rm X}_{-P(k-1)}\\+\langle h_k,\lambda\rangle |
\begin{array}{l}
j_1<\cdots<j_k ,\\ j_k>k,\ j_i\in J_{\rm X}
\end{array}
\}.
\]
For $T=[j_1,\cdots,j_k]_{-P(k-1)}^{\rm X}+\lan h_k,\lambda\ran \in {\rm Tab}_{{\rm X},\iota,k}\setminus \{\lambda^{(k)}\}$,
let us recall that
$[j_1,\cdots,j_k]^{\rm X}_{-P(k-1)}=\sum_{i=1}^{k}\fbox{$j_i$}^{\rm X}_{-P(k-1)+k-i}$, and
by Definition \ref{box-def}, we obtain
\begin{equation}
\fbox{$j_i$}^{\rm X}_{-P(k-1)+k-i}=
\begin{cases}\label{BC-box}
c(j_i)x_{-P(k-1)+k-i+P(j_i),j_i}-x_{-P(k-1)+k-i+P(j_i-1)+1,j_i-1} & {\rm if}\ j_i\leq n,\\
x_{-P(k-1)+k-i+P(|j_i|-1)+n-|j_i|+1,|j_i|-1}-c(j_i)x_{-P(k-1)+k-i+P(|j_i|)+n-|j_i|+1,|j_i|} & {\rm if}\ j_i\geq \ovl{n},
\end{cases}
\end{equation}
where if $\mathfrak{g}$ is of type C and $j_i\in\{n,\ovl{n}\}$ then $c(j_i)=2$,
otherwise $c(j_i)=1$.


Since we supposed (1), (2) in (\ref{cond12}) hold, one obtain $p_{k+1,k}=1$, $p_{k,k-1}=0$ so that
$P(k+1)=p_{k+1,k}+P(k)=1+P(k)=1+P(k-1)$. 
Thus, for $i\in[1,k-1]$ such that $j_i\leq n$,
\begin{equation}\label{BC-pr01}
-P(k-1)+k-i+P(j_i)\geq -P(k-1)+1+(k-1)-i+P(i)\geq -P(k-1)+1+P(k-1)=1,
\end{equation}
where (\ref{pineq}), (\ref{pineq2}) in the second inequality.
If $j_k\leq n$ then since $j_k\geq k+1$,
\begin{equation}\label{BC-pr02}
-P(k-1)+k-k+P(j_k)=-P(k-1)+P(j_k)\geq -P(k-1)+P(k+1)=1.
\end{equation}
For $i\in[1,k-1]$ such that $j_i\leq n$, $j_i>j_{i-1}+1$ (we set $j_0:=0$),
\begin{multline}
-P(k-1)+k-i+P(j_i-1)+1 \label{BC-pr05}\\
\geq -P(k-1)+2+(k-1)-i+P(j_{i-1}+1)\geq -P(k-1)+2+(k-1)-i+P(i)\geq -P(k-1)+2+P(k-1)=2.
\end{multline}
If $j_k\leq n$ such that $j_k>j_{k-1}+1$ then we have $j_k>k+1$ since $j_k=k+1$ yields $j_{k-1}=k-1$, $j_{k-2}=k-2$, $\cdots$, $j_1=1$ and
$T=[1,2,\cdots,k-1,k+1]+\lan h_k,\lambda\ran=\lambda^{(k)}$, which contradicts
$T\in {\rm Tab}_{{\rm X},\iota,k}\setminus \{\lambda^{(k)}\}$.
Hence,
\begin{equation}\label{BC-pr03}
-P(k-1)+k-k+P(j_k-1)+1 =-P(k+1)+P(j_k-1)+2 \geq -P(k+1)+P(k+1)+2=2.
\end{equation}
For $i\in[1,k]$ such that $j_i\geq \ovl{n}$,
using $i\leq k$ and $P(k+1)\leq P(n)$, we also see that
\begin{multline}\label{BCm1}
-P(k-1)+k-i+P(|j_i|-1)+n-|j_i|+1\\
=-P(k-1)+k-i+P(|j_i|-1)+n-(|j_i|-1)
\geq 1-P(k+1)+P(n)\geq 1
\end{multline}
and
\begin{equation}\label{BC-pr04}
-P(k-1)+k-i+P(|j_i|)+n-|j_i|+1 = -P(k+1)+k-i+P(|j_i|)+n-|j_i|+2
\geq -P(k+1)+k-i+P(n)+2\geq2.
\end{equation}
Hence, it follows from (\ref{BC-pr01}), (\ref{BC-pr02}) and (\ref{BCm1}) that
the left indices of 
\[
x_{-P(k-1)+k-i+P(j_i),j_i},\quad x_{-P(k-1)+k-i+P(|j_i|-1)+n-|j_i|+1,|j_i|-1}\]
 in (\ref{BC-box})
are positive. Because of (\ref{BC-pr05}) and (\ref{BC-pr03}), we also see that if $j_i>j_{i-1}+1$ then
the left indices of $x_{-P(k-1)+k-i+P(j_i-1)+1,j_i-1}$ in (\ref{BC-box}) are greater than or equal to $2$.
Furthermore, the inequality (\ref{BC-pr04}) means the
left indices of $x_{-P(k-1)+k-i+P(|j_i|)+n-|j_i|+1,|j_i|}$ in (\ref{BC-box}) are greater than or equal to $2$.
By a similar argument to the proof of Proposition 4.2 (i) in \cite{KaN},
we can prove our claim.

\vspace{3mm}

\nd
(ii) 
We can also prove (ii) by a similar argument
to the proof of Proposition 4.2 (ii) in \cite{KaN}.

\vspace{3mm}

\nd
(iii) 
We take an element $T=[\ovl{n+1},j_2,\cdots,j_k]^{\rm C}_{-P(n-1)}+\lan h_n,\lambda\ran
\in {\rm Tab}_{{\rm C},\iota,n}[\lambda]\setminus \{\lambda^{(n)}\}$.
One can describe $\lambda^{(n)}$ as $\lambda^{(n)}=x_{1,n}-x_{1,n-1}+\lan h_n,\lambda\ran
=[\ovl{n+1},\ovl{n}]^{\rm C}_{1-P(n-1)}+\lan h_n,\lambda\ran$. 
Taking into account that $\ovl{n}\leq j_2<\cdots<j_k$ and $T\neq \lambda^{(n)}$, we see $j_k\neq \ovl{n}$.

We can explicitly write $[\ovl{n+1},j_2,\cdots,j_k]^{\rm C}_{-P(n-1)}$ as
\begin{eqnarray}
& &[\ovl{n+1},j_2,\cdots,j_k]^{\rm C}_{-P(n-1)}
=\fbox{$\ovl{n+1}$}^{\rm C}_{-P(n-1)+k-1} +
\sum^{k}_{i=2} \fbox{$j_i$}^{\rm C}_{-P(n-1)+k-i}\nonumber\\
&=& x_{-P(n-1)+k-1+P(n),n}  \nonumber\\
& & + \sum^{k}_{i=2} (x_{-P(n-1)+k-i+P(|j_i|-1)+n-|j_i|+1,|j_i|-1}-c(|j_i|)x_{-P(n-1)+k-i+P(|j_i|)+n-|j_i|+1,|j_i|}),
\qquad \ \ \label{T-ex}
\end{eqnarray}
where $c(n)=2$ and $c(t)=1$ for $t\in[1,n-1]$. By $k\geq2$, we see that 
\begin{equation}\label{BC-pr06}
-P(n-1)+k-1+P(n)\geq1,
\end{equation}
and for $i\in[2,k]$,
\begin{multline}
-P(n-1)+k-i+P(|j_i|-1)+n-|j_i|+1= \\ \label{BC-pr07}
-P(n-1)+k-i+P(|j_i|-1)+(n-1)-(|j_i|-1)+1\geq -P(n-1)+P(n-1)+1=1.
\end{multline}
\nd
In the case $i<k$, it follows
\begin{equation}\label{BC-pr08}
-P(n-1)+k-i+P(|j_i|)+n-|j_i|+1\geq -P(n-1)+k-i+P(n)+1\geq2,
\end{equation}
and in the case $i=k$, using $j_k\neq \ovl{n}$,
\begin{equation}\label{BC-pr09}
-P(n-1)+k-k+P(|j_k|)+n-|j_k|+1= -P(n-1)+P(|j_k|)+(n-1)-|j_k|+2\geq
-P(n-1)+P(n-1)+2=2.
\end{equation}
The inequalities (\ref{BC-pr06}), (\ref{BC-pr07}) imply
the left indices of
\[
x_{-P(n-1)+k-1+P(n),n},\ \ x_{-P(n-1)+k-i+P(|j_i|-1)+n-|j_i|+1,|j_i|-1}
\]
in (\ref{T-ex}) are positive, and (\ref{BC-pr08}), (\ref{BC-pr09}) imply the left indices of
$x_{-P(n-1)+k-i+P(|j_i|)+n-|j_i|+1,|j_i|}$ 
in (\ref{T-ex}) are greater than or equal to $2$.
By a similar argument to the proof of Proposition 4.2 (iii) in\cite{KaN}, we obtain our claim (iii). 
\qed

\subsection{Actions of operators $\what{S}_{m,j}$ for type D}

\begin{prop}\label{closednessD}
\begin{enumerate}
\item 
For each $T=[j_1,\cdots,j_k]_{-P(k-1)}^{\rm D}+\lan h_k,\lambda\ran\in {\rm Tab}_{{\rm D},\iota,k}[\lambda]\setminus\{\lambda^{(k)}\}$
($k\in [1,n-2]$),
 $j\in I$ and
$m\in\mathbb{Z}_{\geq1}$, we consider the following conditions for the triple $(T,j,m)$:
\begin{enumerate}
\item[(1)] $j<n$ and there exists $i\in [1,k]$ such that $j_i=j$, $j_{i+1}\neq j+1$ and $m=-P(k-1)+k-i+P(j)$,
\item[(2)] $j<n$ and there exists $i'\in [1,k]$ such that $j_{i'}=\overline{j+1}$, $j_{i'+1}\neq \overline{j}$, $n$
and $m=-1-P(k-1)+k-i'+n-j+P(j)$,
\item[(3)] $j<n$ and there exists $i\in [1,k]$ such that $j_{i}= j+1$, $j_{i-1}\neq j$, $\ovl{n}$ and $m=-P(k-1)+k-i+1+P(j)$,
\item[(4)] $j<n$ and there exists $i'\in [1,k]$ such that $j_{i'}= \overline{j}$, $j_{i'-1}\neq\ovl{j+1}$
and $m=-P(k-1)+k-i'+n-j+P(j)$.
\end{enumerate}
We suppose
\begin{equation}\label{lam-cond-D}
\begin{cases}
\iota^{(1)}>\iota^{(2)} & {\rm if}\ k=1,\\
\iota^{(k)}>\iota^{(k-1)},\ \iota^{(k)}>\iota^{(k+1)} & {\rm if}\ 1<k<n-2,\\
\iota^{(n-2)}>\iota^{(n-3)},\ \iota^{(n-2)}>\iota^{(n-1)},\ \iota^{(n-2)}>\iota^{(n)} & {\rm if}\ k=n-2.  
\end{cases}
\end{equation}
If $j<n$ then
\begin{eqnarray*}
\hspace{-15mm}& &\what{S}_{m,j}T= \\
\hspace{-15mm}& &\begin{cases}
[j_1,\cdots,j_{i-1},j+1,j_{i+1},\cdots,j_k]_{-P(k-1)}^{\rm D}+\lan h_k,\lambda\ran & {\rm if\ (1)\ holds\ and\ (2),\ (4)\ do\ not\ hold}, \\
[j_1,\cdots,j_{i'-1},\ovl{j},j_{i'+1},\cdots,j_k]_{-P(k-1)}^{\rm D}+\lan h_k,\lambda\ran & {\rm if\ (2)\ holds\ and\ (1),\ (3)\ do\ not\ hold}, \\
[j_1,\cdots,j_{i-1},j+1,j_{i+1},\cdots,j_{i'-1},\ovl{j},j_{i'+1},\cdots,j_k]_{-P(k-1)}^{\rm D}+\lan h_k,\lambda\ran & {\rm if\ (1)\ and\ (2)\ hold}, \\
[j_1,\cdots,j_{i-1},j,j_{i+1},\cdots,j_k]_{-P(k-1)}^{\rm D}+\lan h_k,\lambda\ran & {\rm if\ (3)\ holds\ and\ (2),\ (4)\ do\ not\ hold}, \ \ \\
[j_1,\cdots,j_{i'-1},\ovl{j+1},j_{i'+1},\cdots,j_k]_{-P(k-1)}^{\rm D}+\lan h_k,\lambda\ran & {\rm if\ (4)\ holds\ and\ (1),\ (3)\ do\ not\ hold}, \\
[j_1,\cdots,j_{i-1},j,j_{i+1},\cdots,j_{i'-1},\ovl{j+1},j_{i'+1},\cdots,j_k]_{-P(k-1)}^{\rm D}+\lan h_k,\lambda\ran & {\rm if\ (3)\ and\ (4)\ hold}, \\
T & {\rm otherwise}.
\end{cases}
\end{eqnarray*}
\item 
For each $T=[j_1,\cdots,j_k]_{-P(k-1)}^{\rm D}+\lan h_k,\lambda\ran\in {\rm Tab}_{{\rm D},\iota,k}[\lambda]\setminus\{\lambda^{(k)}\}$
$(k\in [1,n-2])$ and
$m\in\mathbb{Z}_{\geq1}$,
we consider the following conditions for the pair $(T,m)$:
\begin{enumerate}
\item[(5)] there exists $i\in [1,k]$ such that $j_i=n-1$, $j_{i+1}\neq \overline{n}$, $\ovl{n-1}$ and $m=-P(k-1)+k-i+P(n)$,
\item[(6)] there exists $i\in [1,k]$ such that $j_i=n$, $j_{i+1}\neq \ovl{n},\ \overline{n-1}$ and $m=-P(k-1)+k-i+P(n)$,
\item[(7)] there exists $i\in [1,k]$ such that $j_{i}= \overline{n}$, $j_{i-1}\neq n-1$, $n$ and $m=-P(k-1)+k-i+1+P(n)$,
\item[(8)] there exists $i\in [1,k]$ such that $j_{i}= \overline{n-1}$, $j_{i-1}\neq n-1$, $n$ and $m=-P(k-1)+k-i+1+P(n)$.
\end{enumerate}
We suppose (\ref{lam-cond-D}). Then
\begin{eqnarray*}
\what{S}_{m,n}T&=&\begin{cases}
[j_1,\cdots,j_{i-1},\ovl{n},j_{i+1}\cdots,,j_k]_{-P(k-1)}^{\rm D}+\lan h_k,\lambda\ran & {\rm if\ (5)\ holds}, \\
[j_1,\cdots,j_{i-1},\ovl{n-1},j_{i+1}\cdots,,j_k]_{-P(k-1)}^{\rm D}+\lan h_k,\lambda\ran & {\rm if\ (6)\ holds}, \\
[j_1,\cdots,j_{i-1},n-1,j_{i+1}\cdots,,j_k]_{-P(k-1)}^{\rm D}+\lan h_k,\lambda\ran & {\rm if\ (7)\ holds}, \\
[j_1,\cdots,j_{i-1},n,j_{i+1}\cdots,,j_k]_{-P(k-1)}^{\rm D}+\lan h_k,\lambda\ran & {\rm if\ (8)\ holds}, \\
T & {\rm otherwise}.
\end{cases}
\end{eqnarray*}
\item 
We suppose $t=n-1$ or $t=n$ and $\iota^{(t)}>\iota^{(n-2)}$.
For each $T=[\ovl{n+1},j_2,j_3,\cdots,j_k]_{-P(n-2)}^{\rm D}+\lan h_t,\lambda\ran\in {\rm Tab}_{{\rm D},\iota,t}[\lambda]\setminus\{\lambda^{(t)}\}$,
$j\in I$ and
$m\in\mathbb{Z}_{\geq1}$,
one consider the following conditions (2)', (4)', (9), (10) for the triple $(T,j,m)$:
\begin{enumerate}
\item[(2)'] $j<n$ and there exists $i'\in [1,k]$ such that $j_{i'}=\overline{j+1}$, $j_{i'+1}\neq \overline{j}$
and $m=-1-P(n-2)+k-i'+n-j+P(j)$,
\item[(4)'] $j<n$ and there exists $i'\in [1,k]$ such that $j_{i'}= \overline{j}$, $j_{i'-1}\neq\ovl{j+1}$
and $m=-P(n-2)+k-i'+n-j+P(j)$,
\item[(9)] $j=n$, $j_1=\ovl{n+1}$, $j_2\neq\ovl{n}$, $\ovl{n-1}$ and $m=-P(n-2)+k-1+P(n)$,
\item[(10)] $j=n$, $j_1=\ovl{n+1}$, $j_2=\ovl{n}$, $j_3=\ovl{n-1}$ and $m=-P(n-2)+k-2+P(n)$.
\end{enumerate}
Then we have
\[\what{S}_{m,j}T=
\begin{cases}
[\ovl{n+1},j_2,\cdots,j_{i'-1},\ovl{j},j_{i'+1},\cdots,j_k]_{-P(n-2)}^{\rm D}+\lan h_t,\lambda\ran & {\rm if\ (2)'\ holds}, \\
[\ovl{n+1},j_2,\cdots,j_{i'-1},\ovl{j+1},j_{i'+1},\cdots,j_k]_{-P(n-2)}^{\rm D}+\lan h_t,\lambda\ran & {\rm if\ (4)'\ holds}, \\
[\ovl{n+1},\ovl{n},\ovl{n-1},j_2,\cdots,j_k]_{-P(n-2)}^{\rm D}+\lan h_t,\lambda\ran & {\rm if\ (9)\ holds}, \\
[\ovl{n+1},j_4,\cdots,j_k]_{-P(n-2)}^{\rm D}+\lan h_t,\lambda\ran & {\rm if\ (10)\ holds}, \\
T & {\rm otherwise}.
\end{cases}
\]
\end{enumerate}

\end{prop}

\nd
{\it Proof.}

\nd
(i), (ii)
For $T=[j_1,\cdots,j_k]^{\rm D}_{-P(k-1)}+\lan h_k,\lambda\ran\in {\rm Tab}_{{\rm D},\iota,k}[\lambda]\setminus \{\lambda^{(k)}\}$ with $k\in[1,n-2]$ and $j_1\neq \ovl{n+1}$,
let us consider the action of $\what{S}_{m,j}$ $(m\in\mathbb{Z}_{\geq1},\ j\in I)$.
Recall that 
$[j_1,\cdots,j_k]^{\rm D}_{-P(k-1)}=\sum_{i=1}^{k}\fbox{$j_i$}^{\rm D}_{-P(k-1)+k-i}$, and
by Definition \ref{box-def}, we obtain
\begin{eqnarray}\label{D-box}
& &\fbox{$j_i$}^{\rm D}_{-P(k-1)+k-i}= \nonumber\\
& &\begin{cases}
x_{-P(k-1)+k-i+P(j_i),j_i}-x_{-P(k-1)+k-i+P(j_i-1)+1,j_i-1} & {\rm if}\ j_i\in[1,n-2]\cup\{n\},\\
x_{-P(k-1)+k-i+P(n-1),n-1}+x_{-P(k-1)+k-i+P(n),n}-x_{-P(k-1)+k-i+P(n-2)+1,n-2} & {\rm if}\ j_i=n-1,\\
x_{-P(k-1)+k-i+P(n-1),n-1}-x_{-P(k-1)+k-i+P(n)+1,n} & {\rm if}\ j_i= \ovl{n}, \\ 
x_{-P(k-1)+k-i+P(n-2)+1,n-2}-x_{-P(k-1)+k-i+P(n-1)+1,n-1}
-x_{-P(k-1)+k-i+P(n)+1,n} & {\rm if}\ j_i= \ovl{n-1},\\
x_{-P(k-1)+k-i+P(|j_i|-1)+n-|j_i|,|j_i|-1}-x_{-P(k-1)+k-i+P(|j_i|)+n-|j_i|,|j_i|} & {\rm if}\ j_i\geq \ovl{n-2}.
\end{cases}
\end{eqnarray}
Just as in the proof of Proposition \ref{closednessBC} (i), we can prove for $i\in[1,k]$ such that
$j_i\in[1,n-2]$,
\begin{equation}\label{pr-01d}
-P(k-1)+k-i+P(j_i)\geq1.
\end{equation}
In the case $k<n-2$, it is easy to check
$-P(k-1)+k-i+P(n)=1-P(k+1)+k-i+P(n)\geq1$. In the case $k=n-2$, since we supposed
$\iota^{(k)}=\iota^{(n-2)}>\iota^{(n-1)},\iota^{(n)}$, one obtain $P(n)=P(n-1)$ and 
$-P(n-3)+(n-2)-i+P(n)=1-P(n-1)+(n-2)-i+P(n)\geq1$. Thus, for any $k\in[1,n-2]$ and $i\in[1,k]$, we get
\begin{equation}\label{pr-03d}
-P(k-1)+k-i+P(n)\geq1.
\end{equation}
Since we assume $k\leq n-2$, it follows for any $i\in[1,k]$
\begin{equation}\label{pr-02d}
-P(k-1)+k-i+P(n-1)\geq1,\ \ -P(k-1)+k-i+P(n-2)+1\geq1.
\end{equation}
We also get if $j_i\geq \ovl{n-2}$ then
\begin{eqnarray}
-P(k-1)+k-i+P(|j_i|-1)+n-|j_i| &=& -P(k-1)+k-i+P(|j_i|-1)+(n-1)-(|j_i|-1)\nonumber\\
&\geq & -P(k-1)+k-i+P(n-1)\geq1. \label{pr-04d}
\end{eqnarray}
The inequalities (\ref{pr-01d})-(\ref{pr-04d}) mean that
the left indices of
\begin{multline*}
x_{-P(k-1)+k-i+P(j_i),j_i},\ \ \ x_{-P(k-1)+k-i+P(n-1),n-1},\ \ \ x_{-P(k-1)+k-i+P(n),n},\\ 
x_{-P(k-1)+k-i+P(n-2)+1,n-2},\qquad x_{-P(k-1)+k-i+P(|j_i|-1)+n-|j_i|,|j_i|-1}
\end{multline*}
in (\ref{D-box}) are positive.

By a similar way to the proof of Proposition \ref{closednessBC} (i),
we see that for $i\in[1,k]$ such that $j_i\in[1,n]$ and $j_i>j_{i-1}+1$ ($j_0=0$), it holds
\begin{equation}\label{pr-02e}
-P(k-1)+k-i+P(j_i-1)+1 \geq2.
\end{equation}
In the case $k< n-2$, for any $i\in[1,k]$, we get
$-P(k-1)+k-i+P(n)+1\geq2$. In the case $k=n-2$, by $\iota^{(k)}=\iota^{(n-2)}>\iota^{(n-1)}$, $\iota^{(n)}$,
it holds $P(n-1)=P(n)$, which yields $-P(n-3)+(n-2)-i+P(n)+1\geq2$ for any $i\in[1,n-2]$. Therefore,
for any $i\in[1,k]$, 
\begin{equation}\label{pr-02f}
-P(k-1)+k-i+P(n)+1\geq2.
\end{equation}
It is easy to check for any $i\in[1,k]$,
\begin{equation}\label{pr-02g}
-P(k-1)+k-i+P(n-1)+1 \geq2.
\end{equation}
For $i\in[1,k]$ such that $j_i\geq\ovl{n-2}$, 
\begin{eqnarray}
-P(k-1)+k-i+P(|j_i|)+n-|j_i|
&=& 1-P(k-1)+k-i+P(|j_i|)+(n-1)-|j_i|\nonumber\\
&\geq& 1-P(k-1)+k-i+P(n-1)\geq 2. \label{pr-02h}
\end{eqnarray}

Hence, by (\ref{pr-02e})-(\ref{pr-02h}), 
the left indices of
\begin{multline*}
x_{-P(k-1)+k-i+P(j_i-1)+1,j_i-1},\ \ \ x_{-P(k-1)+k-i+P(n-2)+1,n-2},\ \ \ x_{-P(k-1)+k-i+P(n)+1,n},\\
x_{-P(k-1)+k-i+P(n-1)+1,n-1},\qquad x_{-P(k-1)+k-i+P(|j_i|)+n-|j_i|,|j_i|}
\end{multline*}
in (\ref{D-box}) are greater than or equal to $2$.

By a similar argument to the proof of Proposition 4.3 (i),(ii) in
\cite{KaN}, we can prove our claims (i),(ii).

\nd
(iii) By Definition \ref{box-def} (iv) and Definition \ref{tab-def}, we get
\begin{equation}\label{D-pr-04-a}
\lambda^{(n-1)}=-x_{1,n-1}+x_{1,n-2}+\lan h_{n-1},\lambda\ran =[\ovl{n+1},\ovl{n-1}]_{-P(n-2)}+\lan h_{n-1},\lambda\ran,
\end{equation}
\begin{equation}\label{D-pr-04-b}
\lambda^{(n)}=-x_{1,n}+x_{1,n-2}+\lan h_n,\lambda\ran =[\ovl{n+1},\ovl{n},\ovl{n-1}]_{-P(n-2)}+\lan h_n,\lambda\ran.
\end{equation}

\nd
We can explicitly write $T$ as
\[
T-\lan h_t,\lambda\ran=\fbox{$\ovl{n+1}$}^{\rm D}_{-P(n-2)+k-1} + \sum^{k}_{i=2} \fbox{$j_i$}^{\rm D}_{-P(n-2)+k-i}.
\]
Recall that
\begin{equation}\label{D-box-ag}
\fbox{$j_i$}^{\rm D}_{-P(n-2)+k-i}=
\begin{cases}
x_{-P(n-2)+k-i+P(n),n} & {\rm if}\ j_i=\ovl{n+1}\\
x_{-P(n-2)+k-i+P(n-1),n-1}-x_{-P(n-2)+k-i+P(n)+1,n} & {\rm if}\ j_i= \ovl{n}, \\ 
x_{k-i+1,n-2}-x_{-P(n-2)+k-i+P(n-1)+1,n-1}
-x_{-P(n-2)+k-i+P(n)+1,n} & {\rm if}\ j_i= \ovl{n-1},\\
x_{-P(n-2)+k-i+P(|j_i|-1)+n-|j_i|,|j_i|-1}-x_{-P(n-2)+k-i+P(|j_i|)+n-|j_i|,|j_i|} & {\rm if}\ j_i\geq \ovl{n-2}.
\end{cases}
\end{equation}

If $j_i=\ovl{n+1}$ then $i=1$ by the definition of ${\rm Tab}_{{\rm D},\iota,t}$ ($t=n-1$, $n$).
Combining with $k>1$, we have
\begin{equation}\label{D-pr-03-a}
-P(n-2)+k-i+P(n)\geq1.
\end{equation}

In the case $j_i=\ovl{n}$, it follows from the conditions of ${\rm Tab}_{{\rm D},\iota,t}$
(Definition \ref{tab-def-lam}) that $i<k$. Hence, 
\begin{equation}\label{D-pr-03-b}
-P(n-2)+k-i+P(n-1)\geq1,\ \ \ -P(n-2)+k-i+P(n)+1\geq2.
\end{equation}

For $i\in[1,k]$ such that $j_i= \ovl{n-1}$, if $i=k$ then by the conditions of ${\rm Tab}_{{\rm D},\iota,t}$,
it holds $T=\lambda^{(t)}$, which contradicts our assumption. Thus, it holds $i<k$ and
\begin{equation}\label{D-pr-03-c}
k-i+1\geq1,\ \ \ -P(n-2)+k-i+P(n-1)+1\geq2,\ \ \ -P(n-2)+k-i+P(n)+1\geq2.
\end{equation}
We also get if $j_i\geq \ovl{n-2}$ then
\[
-P(n-2)+k-i+P(|j_i|-1)+n-|j_i|= 1-P(n-2)+k-i+P(|j_i|-1)+(n-2)-(|j_i|-1)\geq 1-P(n-2)+k-i+P(n-2)\geq1,
\]
\[
-P(n-2)+k-i+P(|j_i|)+n-|j_i|=2-P(n-2)+k-i+P(|j_i|)+(n-2)-|j_i|\geq 2-P(n-2)+k-i+P(n-2)\geq2,\]
therefore,
\begin{equation}\label{D-pr-03-d}
-P(n-2)+k-i+P(|j_i|-1)+n-|j_i|\geq1,\ \ \ -P(n-2)+k-i+P(|j_i|)+n-|j_i|\geq2.
\end{equation}
The inequalities (\ref{D-pr-03-a})-(\ref{D-pr-03-d}) mean
the left indices of
\[
x_{-P(n-2)+k-i+P(n),n},\ \ x_{-P(n-2)+k-i+P(n-1),n-1},\ \ x_{k-i+1,n-2},\ \ x_{-P(n-2)+k-i+P(|j_i|-1)+n-|j_i|,|j_i|-1}
\]
in (\ref{D-box-ag}) are positive, and the left indices of
\[
x_{-P(n-2)+k-i+P(n)+1,n},\ \ x_{-P(n-2)+k-i+P(n-1)+1,n-1},\ \ 
x_{-P(n-2)+k-i+P(n)+1,n},\ \ x_{-P(n-2)+k-i+P(|j_i|)+n-|j_i|,|j_i|}
\]
in (\ref{D-box-ag}) are greater than or equal to $2$.
By a similar argument to the proof of Proposition 4.3 (iii) in \cite{KaN},
we can prove our claim (iii). \qed

\section{Proof of Theorem \ref{thm3}, \ref{thm4}}

In this section, we prove our main result Theorem \ref{thm3}.
For $k\in I$,
we set
\[
\Xi_{\iota,k}[\lambda]:=
\{\what{S}_{j_t}\cdots \what{S}_{j_1} \lambda^{(k)} | t\in\mathbb{Z}_{\geq0},\ j_1,\cdots,j_t\in \mathbb{Z}_{\geq1} \}.
\]
Note that the definition (\ref{xilamdef}) means 
\[
\Xi_{\iota}[\lambda]=\{\what{S}_{j_l}\cdots \what{S}_{j_1}x_{j_0} | l\in\mathbb{Z}_{\geq0}, j_0,\cdots,j_l\in \mathbb{Z}_{\geq1} \}
\cup \bigcup_{k\in I} \Xi_{\iota,k}[\lambda].
\]
In \cite{KaN}, we shown that 
$\{S_{j_l}\cdots S_{j_1}x_{j_0} | l\in\mathbb{Z}_{\geq0}, j_0,\cdots,j_l\in \mathbb{Z}_{\geq1} \}={\rm Tab}_{{\rm X},\iota}$ and $\iota$
satisfies the positivity condition (Theorem \ref{thm2}, \ref{thm1}). By the definitions of $S$, $\what{S}$ ((\ref{Sk}),(\ref{Shat})),
it holds $\{\what{S}_{j_l}\cdots \what{S}_{j_1}x_{j_0} | l\in\mathbb{Z}_{\geq0}, j_0,\cdots,j_l\in \mathbb{Z}_{\geq1} \}={\rm Tab}_{{\rm X},\iota}$.
Thus, we need to prove $\Xi_{\iota,k}[\lambda]={\rm Tab}_{{\rm X},\iota,k}[\lambda]\cup\{0\}$.
In what follows, we consider the conditions (1), (2) in (\ref{cond12}).

\subsection{Proof of Theorem \ref{thm3} for type A-case}\label{A-proof}

\nd
In the case both (1) and (2) in (\ref{cond12}) do not hold,
by (\ref{lmi-def}), we have $\lambda^{(k)}=-x_{1,k}+\lan h_k,\lambda \ran$.
It follows from (\ref{betadef}) and (\ref{Shat}) that
\[\what{S}_{l,j}\lambda^{(k)}=
\begin{cases}
0 & {\rm if}\ (l,j)=(1,k), \\
\lambda^{(k)} & {\rm otherwise},
\end{cases}
\]
which yields $\Xi_{\iota,k}[\lambda]=\{0,\lambda^{(k)}\}$.

Next, let us consider the case 
only (1) holds, which means $p_{k+1,k}=1$. Note that
$P(k+1)=p_{k+1,k}+P(k)=1+P(k)$.
Taking (\ref{lmi-def}) and Definition \ref{box-def} (i) into account, we obtain $\lambda^{(k)}=-x_{1,k}+x_{1,k+1}+\lan h_k,\lambda \ran=
\fbox{$k+1$}^{\rm A}_{1-P(k+1)}+\lan h_k,\lambda \ran$.
By Proposition \ref{closednessA}, it holds
\[\what{S}_{l,j}\lambda^{(k)}=
\begin{cases}
0 & {\rm if}\ (l,j)=(1,k), \\
\fbox{$k+2$}^{\rm A}_{1-P(k+1)}+\lan h_k,\lambda \ran& {\rm if}\ (l,j)=(1,k+1), \\
\lambda^{(k)}=\fbox{$k+1$}^{\rm A}_{1-P(k+1)}+\lan h_k,\lambda\ran & {\rm otherwise}.
\end{cases}
\]
Note that,
for $t\in[k+2,n+1]$, it holds $\fbox{$t$}^{\rm A}_{1-P(k+1)}=x_{1-P(k+1)+P(t),t}-x_{2-P(k+1)+P(t-1),t-1}$ and by (\ref{pineq}),
\begin{equation}\label{A-ineq1}
1-P(k+1)+P(t)\geq 1-P(k+1)+P(k+2)\geq1,\ \ 2-P(k+1)+P(t-1)\geq 2-P(k+1)+P(k+1)=2.
\end{equation}
By Proposition \ref{closednessA}, we get
\[\what{S}_{l,j} \fbox{$t$}^{\rm A}_{1-P(k+1)}=
\begin{cases}
\fbox{$t+1$}^{\rm A}_{1-P(k+1)} & {\rm if}\ (l,j)=(1-P(k+1)+P(t),t), \\
\fbox{$t-1$}^{\rm A}_{1-P(k+1)} & {\rm if}\ (l,j)=(2-P(k+1)+P(t-1),t-1), \\
\fbox{$t$}^{\rm A}_{1-P(k+1)} & {\rm otherwise},
\end{cases}
\]
which yields $\Xi_{\iota,k}[\lambda]=\{0\}\cup\{\fbox{$t$}^{\rm A}_{1-P(k+1)}+\lan h_k,\lambda \ran| k+1\leq t\leq n+1\}$.

Next, we consider the case 
only (2) holds, which means $p_{k+1,k}=0$, $p_{k,k-1}=0$. In this setting, by Definition \ref{tab-def-lam},
\[
{\rm Tab}_{{\rm A},\iota,k}[\lambda]=\{[j_1,\cdots,j_{k-1},k+1,\cdots,n+1]^{\rm A}_{-P(k-1)-n+k}+
\lan h_k,\lambda \ran| 1\leq j_1<\cdots<j_{k-1}\leq k\}.
\]
By (\ref{lmi-def}), Definition \ref{box-def} (i) and Definition \ref{tab-def} (i), it holds $\lambda^{(k)}=-x_{1,k}+x_{1,k-1}+\lan h_k,\lambda \ran=
[1,2,\cdots,k-1]^{\rm A}_{1-P(k-1)}-x_{1,k}+\lan h_k,\lambda \ran=[1,2,\cdots,k-1,k+1,\cdots,n+1]^{\rm A}_{-P(k-1)-n+k}+\lan h_k,\lambda \ran$.
Thus, it holds $\lambda^{(k)}\in {\rm Tab}_{{\rm A},\iota,k}[\lambda]$.
Considering Proposition \ref{closednessA}, we can verify that
\[\what{S}_{r,j}\lambda^{(k)}=
\begin{cases}
0 & {\rm if}\ (r,j)=(1,k), \\
[1,\cdots,k-2,k,k+1,\cdots,n+1]^{\rm A}_{-P(k-1)-n+k}+\lan h_k,\lambda \ran& {\rm if}\ (r,j)=(1,k-1), \\
\lambda^{(k)} & {\rm otherwise},
\end{cases}
\]
which means $\what{S}_{r,j}\lambda^{(k)}\in {\rm Tab}_{{\rm A},\iota,k}[\lambda]\cup\{0\}$ for any $(r,j)$.
Note that each element $T=[j_1,\cdots,j_{k-1},k+1,\cdots,n+1]^{\rm A}_{-P(k-1)-n+k}+\lan h_k,\lambda \ran\in {\rm Tab}_{{\rm A},\iota,k}[\lambda]$
other than $\lambda^{(k)}$
 can be written
as $T=[1,\cdots,l-1,l+1,\cdots,n+1]^{\rm A}_{-P(k-1)-n+k}+\lan h_k,\lambda \ran$ with some $l\in[1,k-1]$,
which implies $j_{l-1}=l-1$, $j_{l}=l+1$. Putting
$m:=(-P(k-1)-n+k)+n-(l-1)+P(l-1)$, $m':=(-P(k-1)-n+k)+n-l+1+P(l)$
we obtain
\[
m=-P(k-1)+1+(k-1)-(l-1)+P(l-1)\geq -P(k-1)+1+P(k-1)=1,
\]
\[
m'=-P(k-1)+(k-1)-l+2+P(l)\geq -P(k-1)+2+P(k-1)=2,
\]
where we use (\ref{pineq2}) in the above inequalities.
Thus, using Proposition \ref{closednessA},
\[\what{S}_{r,j}[1,\cdots,l-1,l+1,\cdots,n+1]^{\rm A}_{-P(k-1)-n+k}=
\begin{cases}
[1,\cdots,l-2,l,\cdots,n+1]^{\rm A}_{-P(k-1)-n+k} & {\rm if}\ (r,j)=(m,l-1), \\
[1,\cdots,l,l+2,\cdots,n+1]^{\rm A}_{-P(k-1)-n+k} & {\rm if}\ (r,j)=(m',l), \\
[1,\cdots,l-1,l+1,\cdots,n+1]^{\rm A}_{-P(k-1)-n+k} & {\rm otherwise},
\end{cases}
\]
which yields $\what{S}_{r,j}[1,\cdots,l-1,l+1,\cdots,n+1]^{\rm A}_{-P(k-1)-n+k}+\lan h_k,\lambda \ran
\in {\rm Tab}_{{\rm A},\iota,k}[\lambda]\cup\{0\}$ for any $(r,j)$.
Thus, we get $\Xi_{\iota,k}[\lambda]={\rm Tab}_{{\rm A},\iota,k}[\lambda]\cup\{0\}$.

Finally, let us turn to the case both (1) and (2) hold, which means $p_{k+1,k}=1$, $p_{k,k-1}=0$, $P(k+1)=1+P(k)$ and $P(k)=P(k-1)$.
First, we prove $\Xi_{\iota,k}[\lambda]\subset {\rm Tab}_{{\rm A},\iota,k}[\lambda] \cup\{0\}$.
We see that $\lambda^{(k)}=x_{1,k-1}+x_{1,k+1}-x_{1,k}+\lan h_k,\lambda \ran
=[1,2,\cdots,k-1]^{\rm A}_{1-P(k-1)} + \fbox{$k+1$}^{\rm A}_{-P(k-1)}+\lan h_k,\lambda \ran
=[1,2,\cdots,k-1,k+1]^{\rm A}_{-P(k-1)}+\lan h_k,\lambda \ran\in {\rm Tab}_{{\rm A},\iota,k}[\lambda]
=\{
[j_1,\cdots,j_k]_{-P(k-1)}+\lan h_k,\lambda\ran | 1\leq j_1<\cdots<j_k \leq n+1,\ j_k>k
\}$. Considering Proposition \ref{closednessA}, it holds
\[
\what{S}_{m,j}\lambda^{(k)}=
\begin{cases}
[1,2,\cdots,k-1,k+2]^{\rm A}_{-P(k-1)}+\lan h_k,\lambda \ran & {\rm if}\ (m,j)=(1,k+1), \\
[1,2,\cdots,k-2,k,k+1]^{\rm A}_{-P(k-1)}+\lan h_k,\lambda \ran & {\rm if}\ (m,j)=(1,k-1), \\
0 & {\rm if}\ (m,j)=(1,k),\\
\lambda^{(k)} & {\rm otherwise}.
\end{cases}
\]
According to Proposition \ref{closednessA}, for each $T=[j_1,\cdots,j_k]^{\rm A}_{-P(k-1)}+\lan h_k,\lambda\ran\in {\rm Tab}_{{\rm A},\iota,k}[\lambda] $, except for the following case, we have
$\what{S}_{m,j}T\in{\rm Tab}_{{\rm A},\iota,k}[\lambda]\cup\{0\}$:
\begin{equation}\label{pr-A-1}
j_i=j+1,\ j_{i-1}\neq j,\ m=-P(k-1)+k-i+1+P(j)=1\ {\rm for\ some}\ i\in[1,k].
\end{equation}
In the case $i\in[1,k-1]$, the condition (\ref{pr-A-1}) does not hold. To prove it, 
we assume (\ref{pr-A-1}) holds and deduce a contradiction from this assumption.
If $j_i=i$ then the condition $j_1<\cdots<j_k$ of ${\rm Tab}_{{\rm A},\iota,k}[\lambda] $ means
$j_1=1,\ j_2=2,\ \cdots,\ j_{i-1}=i-1=j$, which contradicts $j_{i-1}\neq j$ in (\ref{pr-A-1}) and if
$j_i>i$ then
\begin{eqnarray}
m=-P(k-1)+k-i+1+P(j)
&=&-P(k-1)+k-i+1+P(j_i-1) \nonumber\\
&\geq& -P(k-1)+k-i+1+P(i)\nonumber\\
&=& -P(k-1)+(k-1)-i+2+P(i)\nonumber\\
&\geq& -P(k-1)+2+P(k-1)=2, \label{pr-A-01}
\end{eqnarray}
which contradicts $m=1$ in (\ref{pr-A-1}), where the above two inequalities follow from
(\ref{pineq}), (\ref{pineq2}). Thus, we proved in the case $i\in[1,k-1]$, the condition (\ref{pr-A-1}) does not hold.
In the case $i=k$, if $(j_i=)j_k>k+1$ then (\ref{pr-A-1}) does not hold. To prove it, we assume (\ref{pr-A-1}) holds.
In conjunction with $j_i>k+1$, the condition $j_i=j+1$ means $j>k$, which yields that
\begin{equation}\label{pr-A-2}
m=-P(k-1)+k-k+1+P(j)=-P(k-1)+1+P(j)\geq -P(k-1)+1+P(k+1)=2,
\end{equation}
which contradicts $m=1$ in (\ref{pr-A-1}).
Hence, if $(j_i=)j_k>k+1$ then (\ref{pr-A-1}) does not hold.
Only in the case $j_i=j_k=k+1$ and $j_{k-1}=k-1$, $j_{k-2}=k-2$, $\cdots$, $j_1=1$, the condition (\ref{pr-A-1}) holds,
which implies $T=[1,2,\cdots,k-1,k+1]^{\rm A}_{-P(k-1)}+\lan h_k,\lambda\ran=\lambda^{(k)}$.
By the above argument, we see that
${\rm Tab}_{{\rm A},\iota,k}[\lambda] \cup\{0\}$ is closed under the action of $\what{S}_{m,j}$ ($m\in\mathbb{Z}_{\geq1}$, $j\in I$).
Since we know $\lambda^{(k)}\in{\rm Tab}_{{\rm A},\iota,k}[\lambda] \cup\{0\}$,
one obtain $\Xi_{\iota,k}[\lambda]\subset {\rm Tab}_{{\rm A},\iota,k}[\lambda] \cup\{0\}$.

Next, we show ${\rm Tab}_{{\rm A},\iota,k}[\lambda] \cup\{0\}\subset\Xi_{\iota,k}[\lambda]$.
For each $T=[j_1,\cdots,j_k]^{\rm A}_{-P(k-1)}+\lan h_k,\lambda\ran\in {\rm Tab}_{{\rm A},\iota,k}[\lambda] $, we show
$T\in \Xi_{\iota,k}[\lambda]$
using induction on the value $j_1+\cdots+j_k$.
By $1\leq j_1<\cdots<j_k\leq n+1$ and $j_k>k$, the minimal value of $j_1+\cdots+j_k$
is $1+2+\cdots+(k-1)+(k+1)$. In this case, we can easily check that
$j_1=1$, $j_2=2$, $\cdots$, $j_{k-1}=k-1$, $j_k=k+1$
and
$T=[1,2,\cdots,k-1,k+1]^{\rm A}_{-P(k-1)}+\lan h_k,\lambda\ran
=\lambda^{(k)}\in \Xi_{\iota,k}[\lambda]$.

Next, we assume $j_1+\cdots+j_k>1+2+\cdots+(k-1)+(k+1)$.
Either the following (i) or (ii) holds:

\nd
(i) $j_i>j_{i-1}+1$ for some $i\in[1,k-1]$ (we set $j_0=0$),

\nd
(ii) $j_l=l$ for $l\in[1,k-1]$ and $j_k>k+1$.

If (i) holds then $j_i-1>j_{i-1}\geq i-1$ so that $j_i-1\geq i$.
Putting $m:=-P(k-1)+k-i+P(j_i-1)$, we see that
\begin{multline*}
m=-P(k-1)+k-i+P(j_i-1)\geq -P(k-1)+k-i+P(i)= \\
-P(k-1)+1+(k-1)-i+P(i) \geq -P(k-1)+1 +P(k-1)=1. 
\end{multline*}
It follows by Proposition \ref{closednessA} that
\begin{equation}\label{A-pr1}
\what{S}_{m,j_{i}-1}[j_1,\cdots,j_{i-1},j_i-1,j_{i+1},\cdots,j_k]^{\rm A}_{-P(k-1)}=
[j_1,\cdots,j_{i-1},j_i,j_{i+1},\cdots,j_k]^{\rm A}_{-P(k-1)}.
\end{equation}
Since $j_i>j_{i-1}+1$, we obtain $[j_1,\cdots,j_{i-1},j_i-1,j_{i+1},\cdots,j_k]^{\rm A}_{-P(k-1)}+\lan h_k,\lambda\ran\in{\rm Tab}_{{\rm A},\iota,k}[\lambda]$.
By the induction assumption,
we get $[j_1,\cdots,j_i-1,\cdots,j_k]^{\rm A}_{-P(k-1)}+\lan h_k,\lambda\ran\in \Xi_{\iota,k}[\lambda]$, which implies
$[j_1,\cdots,j_i-1,\cdots,j_k]^{\rm A}_{-P(k-1)}+\lan h_k,\lambda\ran=\what{S}_{l_p}\cdots \what{S}_{l_2}\what{S}_{l_1}\lambda^{(k)}$
with some $l_1,\ \cdots,\ l_p\in\mathbb{Z}_{\geq1}$.
In conjunction with (\ref{A-pr1}), we obtain $[j_1,\cdots,j_i,\cdots,j_k]^{\rm A}_{-P(k-1)}+\lan h_k,\lambda\ran\in \Xi_{\iota,k}[\lambda]$. 

If (ii) holds then putting $m:=-P(k-1)+P(j_k-1)$, we obtain $m=-P(k-1)+P(j_k-1)\geq -P(k-1)+P(k+1)=1$.
Using Proposition \ref{closednessA}, one obtain
\begin{equation}\label{A-pr2}
\what{S}_{m,j_{k}-1}[j_1,\cdots,j_{k-1},j_k-1]^{\rm A}_{-P(k-1)}=
[j_1,\cdots,j_{k-1},j_k]^{\rm A}_{-P(k-1)}.
\end{equation}
By the induction assumption,
we see $[j_1,\cdots,j_{k-1},j_k-1]^{\rm A}_{-P(k-1)}+\lan h_k,\lambda\ran\in \Xi_{\iota,k}[\lambda]$ and 
(\ref{A-pr2}) yields $[j_1,\cdots,j_k]^{\rm A}_{-P(k-1)}+\lan h_k,\lambda\ran\in \Xi_{\iota,k}[\lambda]$.
Therefore, the inclusion ${\rm Tab}_{{\rm A},\iota,k}[\lambda] \cup\{0\}\subset\Xi_{\iota,k}[\lambda]$ follows.

\qed


\subsection{Type B-case}

\vspace{2mm}

\nd
\underline{Case 1 : $k<n$.}

\vspace{2mm}

\nd
First, we suppose $k<n$.

\vspace{2mm}

\underline{Case 1-1 : the case both (1) and (2) do not hold}

\vspace{2mm}

In this case,
using (\ref{lmi-def}), we have $\lambda^{(k)}=-x_{1,k}+\lan h_k,\lambda \ran$
and
$\Xi_{\iota,k}[\lambda]=\{0,\lambda^{(k)}\}$ by a similar way to the type A-case.

\vspace{2mm}

\underline{Case 1-2 : the case only (1) holds} 

\vspace{2mm}

In this case, we have $p_{k+1,k}=1$. Thus, it holds
$P(k+1)=p_{k+1,k}+P(k)=1+P(k)$.
The Definition \ref{box-def} (ii) and (\ref{lmi-def}) say $\lambda^{(k)}=-x_{1,k}+x_{1,k+1}+\lan h_k,\lambda \ran=
\fbox{$k+1$}^{\rm B}_{1-P(k+1)}+\lan h_k,\lambda \ran$.
We obtain
\[
\fbox{$t$}^{\rm B}_{1-P(k+1)}=
\begin{cases}
x_{1-P(k+1)+P(t),t}-x_{-P(k+1)+P(t-1)+2,t-1} & {\rm if}\ t\leq n,\\
x_{-P(k+1)+P(|t|-1)+n-|t|+2,|t|-1}-x_{-P(k+1)+P(|t|)+n-|t|+2,|t|} & {\rm if}\ t\geq \ovl{n}.
\end{cases}
\]
If $t\in[k+1,n]$ then it is easy to see
\begin{equation}\label{B-6-2-1}
1-P(k+1)+P(t)\geq1,
\end{equation}
and if $t\in[k+2,n]$ then
\begin{equation}\label{B-6-2-2}
-P(k+1)+P(t-1)+2\geq2.
\end{equation}
If $t\geq \ovl{n}$ then
\begin{equation}\label{B-6-2-3}
-P(k+1)+P(|t|-1)+n-|t|+2=1-P(k+1)+P(|t|-1)+n-(|t|-1)\geq 1-P(k+1)+P(n) \geq1,
\end{equation}
\begin{equation}\label{B-6-2-4}
-P(k+1)+P(|t|)+n-|t|+2\geq 2-P(k+1)+P(n) \geq2.
\end{equation}
Putting $s:=1-P(k+1)$, $\what{S}(j):=\what{S}_{s+P(j),j}$, $\what{S}'(j):=\what{S}_{s+P(j)+1,j}$ for $j\in[k+1,n]$ and
$\what{S}(\ovl{j}):=\what{S}_{s+P(j)+n-j,j}$, $\what{S}'(\ovl{j}):=\what{S}_{s+P(j)+n-j+1,j}$ for $j\in[1,n-1]$,
by Lemma \ref{box-lem}, we obtain the following diagram of actions of $\what{S}$:
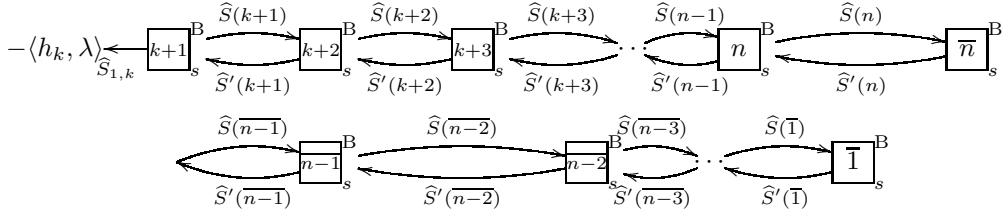
\[
\begin{xy}
(-16,0) *{-\lan h_k,\lambda\ran}="00",
(0,0) *{\begin{ytableau}
\scriptstyle k+1
\end{ytableau}^{\rm B}_{s}}="1",
(20,0)*{\begin{ytableau}
\scriptstyle k+2
\end{ytableau}^{\rm B}_{s}}="2",
(40,0)*{
\begin{ytableau}
\scriptstyle k+3
\end{ytableau}^{\rm B}_{s}}="3",
(60,0)*{\ \cdots\ }="dots",
(75,0)*{\begin{ytableau}
n
\end{ytableau}^{\rm B}_{s}}="r",
(105,0)*{\begin{ytableau}
\overline{n}
\end{ytableau}^{\rm B}_{s}}="r+1",
(0,-15) *{\ }="0",
(20,-15) *{\begin{ytableau}
\scriptstyle \overline{n-1}
\end{ytableau}^{\rm B}_{s}}="r-1r-1",
(55,-15)*{\begin{ytableau}
\scriptstyle \overline{n-2}
\end{ytableau}^{\rm B}_{s}}="r-2r-2",
(70,-15)*{
\ 
\cdots
\ }="dotdot",
(90,-15)*{\begin{ytableau}
\overline{1}
\end{ytableau}^{\rm B}_{s}}="11",
\ar@{->} "1";"00"^{\what{S}_{1,k}}
\ar@/^/ @{->} "1";"2"^{\what{S}(k+1)}
\ar@/^/ @{->} "2";"3"^{\what{S}(k+2)}
\ar@/^/ @{->} "3";"dots"^{\ \what{S}(k+3)}
\ar@/^/ @{->} "dots";"r"^{\what{S}(n-1)}
\ar@/^/ @{->} "r";"r+1"^{\what{S}(n)}
\ar@/^/ @{->} "r-1r-1";"r-2r-2"^{\what{S}(\ovl{n-2})}
\ar@/^/ @{->} "r-2r-2";"dotdot"^{\what{S}(\ovl{n-3})}
\ar@/^/ @{->} "dotdot";"11"^{\what{S}(\ovl{1})}
\ar@/^/ @{->} "0";"r-1r-1"^{\what{S}(\ovl{n-1})}
\ar@/_/ @{<-} "1";"2"_{\what{S}'(k+1)}
\ar@/_/ @{<-} "2";"3"_{\what{S}'(k+2)}
\ar@/_/ @{<-} "3";"dots"_{\ \what{S}'(k+3)}
\ar@/_/ @{<-} "dots";"r"_{\what{S}'(n-1)}
\ar@/_/ @{<-} "r";"r+1"_{\what{S}'(n)}
\ar@/_/ @{<-} "r-1r-1";"r-2r-2"_{\what{S}'(\ovl{n-2})}
\ar@/_/ @{<-} "r-2r-2";"dotdot"_{\what{S}'(\ovl{n-3})}
\ar@/_/ @{<-} "dotdot";"11"_{\what{S}'(\ovl{1})}
\ar@/_/ @{<-} "0";"r-1r-1"_{\what{S}'(\ovl{n-1})}
\end{xy}
\]
Other actions of $\what{S}$ are trivial.
Therefore, it holds
$\Xi_{\iota,k}[\lambda]=\{0\}\cup\{\fbox{$t$}^{\rm B}_{1-P(k+1)}+\lan h_k,\lambda\ran| k+1\leq t\leq \overline{1}\}$.

\vspace{2mm}

\underline{Case 1-3 : the case only (2) holds}

\vspace{2mm}

In this case, we obtain $p_{k,k-1}=0$ so that $P(k)=p_{k,k-1}+P(k-1)=P(k-1)$.
We see that $\lambda^{(k)}=-x_{1,k}+x_{1,k-1}+\lan h_k,\lambda \ran=
\fbox{$\ovl{k}$}^{\rm B}_{-P(k-1)-n+k}+\lan h_k,\lambda \ran$ by Definition \ref{box-def} (ii).
For $t\in[1,k]$,
\[
\fbox{$\overline{t}$}^{\rm B}_{-P(k-1)-n+k}
=x_{-P(k-1)+k+P(t-1)-t+1,t-1}-x_{-P(k-1)+k+P(t)-t+1,t}
\]
and
\begin{equation}\label{B-6-2-5}
-P(k-1)+k+P(t-1)-t+1\geq 1.
\end{equation}
If $t\leq k-1$ then
\begin{equation}\label{B-6-2-6}
-P(k-1)+k+P(t)-t+1=-P(k-1)+P(t)+(k-1)-t+2\geq -P(k-1)+P(k-1)+2= 2.
\end{equation}
By Lemma \ref{box-lem}, putting $s:=-P(k-1)-n+k$ and
$\what{S}(\ovl{j}):=\what{S}_{s+P(j)+n-j,j}$, $\what{S}'(\ovl{j}):=\what{S}_{s+P(j)+n-j+1,j}$ for $j\in[1,k-1]$,
we get the following diagram of actions of $\what{S}$:
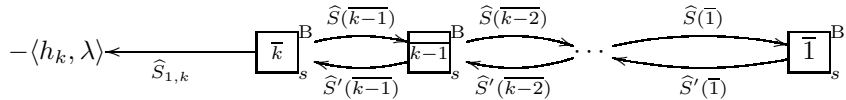
\[
\begin{xy}
(-10,0) *{-\lan h_k,\lambda\ran}="00",
(0,0) *{\ }="0",
(20,0) *{\begin{ytableau}
\scriptstyle \overline{k}
\end{ytableau}^{\rm B}_{s}}="r-1r-1",
(40,0)*{\begin{ytableau}
\scriptstyle \overline{k-1}
\end{ytableau}^{\rm B}_{s}}="r-2r-2",
(60,0)*{
\ 
\cdots
\ }="dotdot",
(90,0)*{\begin{ytableau}
\overline{1}
\end{ytableau}^{\rm B}_{s}}="11",
\ar@{->} "r-1r-1";"00"^{\what{S}_{1,k}}
\ar@/^/ @{->} "r-1r-1";"r-2r-2"^{\what{S}(\ovl{k-1})}
\ar@/^/ @{->} "r-2r-2";"dotdot"^{\what{S}(\ovl{k-2})}
\ar@/^/ @{->} "dotdot";"11"^{\what{S}(\ovl{1})}
\ar@/_/ @{<-} "r-1r-1";"r-2r-2"_{\what{S}'(\ovl{k-1})}
\ar@/_/ @{<-} "r-2r-2";"dotdot"_{\what{S}'(\ovl{k-2})}
\ar@/_/ @{<-} "dotdot";"11"_{\what{S}'(\ovl{1})}
\end{xy}
\]
Thus, we get $\Xi_{\iota,k}[\lambda]=\{0\}\cup
\{\fbox{$t$}^{\rm B}_{-P(k-1)-n+k}+\lan h_k,\lambda \ran| \ovl{k}\leq t\leq \ovl{1}\}$.

\vspace{2mm}

\underline{Case 1-4: the case both (1) and (2) hold}

\vspace{2mm}

In this case, we obtain $p_{k,k-1}=0$, $p_{k+1,k}=1$ so that $P(k+1)=1+P(k)=1+P(k-1)$
and $\lambda^{(k)}=-x_{1,k}+x_{1,k-1}+x_{1,k+1}+\lan h_k,\lambda\ran
=[1,\cdots,k-1,k+1]^{\rm B}_{-P(k-1)}+\lan h_k,\lambda\ran$ by Definition \ref{box-def} (ii), \ref{tab-def} (i).
It follows from (\ref{betadef}) and (\ref{Shat}) that
\[\what{S}_{l,j}\lambda^{(k)}=
\begin{cases}
0 & {\rm if}\ (l,j)=(1,k), \\
[1,\cdots,k-2,k,k+1]^{\rm B}_{-P(k-1)}+\lan h_k,\lambda \ran& {\rm if}\ (l,j)=(1,k-1), \\
[1,\cdots,k-1,k+2]^{\rm B}_{-P(k-1)}+\lan h_k,\lambda \ran& {\rm if}\ (l,j)=(1,k+1), \\
\lambda^{(k)} & {\rm otherwise}.
\end{cases}
\]
Combining with Proposition \ref{closednessBC} (i), we see that
${\rm Tab}_{{\rm B},\iota,k}[\lambda]\cup\{0\}$ is closed under the action of $\what{S}_{l,j}$ for all $(l,j)\in\mathbb{Z}_{\geq1}\times I$.
By $\lambda^{(k)}=[1,\cdots,k-1,k+1]^{\rm B}_{-P(k-1)}+\lan h_k,\lambda\ran\in {\rm Tab}_{{\rm B},\iota,k}[\lambda]\cup\{0\}$,
we obtain $\Xi_{\iota,k}[\lambda]\subset{\rm Tab}_{{\rm B},\iota,k}[\lambda]\cup\{0\}$.

We can prove the inclusion ${\rm Tab}_{{\rm B},\iota,k}[\lambda]\cup\{0\}\subset \Xi_{\iota,k}[\lambda]$ by a similar way
to the proof of Lemma 5.6 in \cite{KaN}.

\vspace{2mm}

\nd
\underline{Case 2 : $k=n$.}

\nd
Next, let us turn to the case $k=n$. The condition (1) does not hold.

\vspace{2mm}

\underline{Case 2-1 : the case the condition (2) does not hold}

\vspace{2mm}

We have
$\lambda^{(n)}=-x_{1,n}+\lan h_n,\lambda \ran$ and
$\Xi_{\iota,n}[\lambda]=\{0,\lambda^{(n)}\}$ by a similar argument to Case 1-1.

\vspace{2mm}

\underline{Case 2-2 : the case the condition (2) holds}

\vspace{2mm}

If the condition (2) holds then $p_{n,n-1}=0$ so that $P(n)=P(n-1)$ and
$\lambda^{(n)}=-x_{1,n}+2x_{1,n-1}+\lan h_n,\lambda \ran=[1,2,\cdots,n-1,\ovl{n}]_{-P(n-1)}
+\lan h_n,\lambda \ran$ by Definition \ref{box-def} (ii), \ref{tab-def} (i).
By a direct calculation, we can verify
\[\what{S}_{l,j}\lambda^{(n)}=
\begin{cases}
0 & {\rm if}\ (l,j)=(1,n), \\
[1,\cdots,n-2,n,\ovl{n-1}]^{\rm B}_{-P(n-1)}+\lan h_n,\lambda \ran& {\rm if}\ (l,j)=(1,n-1), \\
\lambda^{(n)} & {\rm otherwise}.
\end{cases}
\]
In conjunction with Proposition \ref{closednessBC} (ii), we see that
${\rm Tab}_{{\rm B},\iota,n}[\lambda]\cup\{0\}$ is closed under the
action of $\what{S}_{l,j}$ for all $(l,j)\in\mathbb{Z}_{\geq1}\times I$ and
$\Xi_{\iota,n}[\lambda]\subset{\rm Tab}_{{\rm B},\iota,n}[\lambda]\cup\{0\}$.

We can also show the inclusion ${\rm Tab}_{{\rm B},\iota,n}[\lambda]\cup\{0\}\subset \Xi_{\iota,n}[\lambda]$ by a similar way
to the proof of Lemma 5.7 (i) in \cite{KaN}. \qed

\subsection{Type C-case}

\vspace{2mm}

\nd
\underline{Case 1 : $k<n$.}

\vspace{2mm}

\nd
First, we suppose $k<n$.

\vspace{2mm}

\underline{Case 1-1 : the case both (1) and (2) do not hold}

\vspace{2mm}

By (\ref{lmi-def}), we have $\lambda^{(k)}=-x_{1,k}+\lan h_k,\lambda \ran$
and
$\Xi_{\iota,k}[\lambda]=\{0,\lambda^{(k)}\}$ by a similar way to the type A, B-cases.

\vspace{2mm}

\underline{Case 1-2 : the case only (1) holds}

\vspace{2mm}

In this case, we have $p_{k+1,k}=1$ so that
$P(k+1)=p_{k+1,k}+P(k)=1+P(k)$.
The Definition \ref{box-def} (iii) and (\ref{lmi-def}) mean $\lambda^{(k)}=-x_{1,k}+(1+\delta_{k+1,n})x_{1,k+1}
+\lan h_k,\lambda \ran=
\fbox{$k+1$}^{\rm C}_{1-P(k+1)}+\lan h_k,\lambda \ran$.
Hence, by Lemma \ref{box-lem} and a similar argument to type A, B, we obtain
$\Xi_{\iota,k}[\lambda]=\{0\}\cup\{\fbox{$t$}^{\rm C}_{1-P(k+1)}+\lan h_k,\lambda\ran| k+1\leq t\leq \overline{1}\}$.

\vspace{2mm}

\underline{Case 1-3 : the case only (2) holds}

\vspace{2mm}

In this case, we have $p_{k,k-1}=0$ so that
$P(k)=P(k-1)$. Considering the Definition \ref{box-def} (iii) and (\ref{lmi-def}),
we obtain
$\lambda^{(k)}=-x_{1,k}+x_{1,k-1}
+\lan h_k,\lambda \ran=
\fbox{$\ovl{k}$}^{\rm C}_{-P(k-1)-n+k}+\lan h_k,\lambda \ran$. By a similar argument to Case 1-3 of the type B-case,
we can verify
$\Xi_{\iota,k}[\lambda]=\{0\}\cup\{\fbox{$t$}^{\rm C}_{-P(k-1)-n+k}+\lan h_k,\lambda\ran| \ovl{k}\leq t\leq \overline{1}\}$.

\vspace{2mm}

\underline{Case 1-4 : the case both (1), (2) hold}

\vspace{2mm}

In this case, we have $p_{k+1,k}=1$, $p_{k,k-1}=0$ so that
$P(k+1)=p_{k+1,k}+P(k)=1+P(k)=1+P(k-1)$.
It holds $\lambda^{(k)}=-x_{1,k}+(1+\delta_{k+1,n})x_{1,k+1}+x_{1,k-1}
+\lan h_k,\lambda \ran=
[1,2,\cdots,k-1,k+1]^{\rm C}_{-P(k-1)}+\lan h_k,\lambda \ran$.
A similar argument to Case 1-4 of the type B-case shows
${\rm Tab}_{{\rm B},\iota,k}[\lambda]\cup\{0\}= \Xi_{\iota,k}[\lambda]$.

\vspace{2mm}

\nd
\underline{Case 2 : $k=n$.}

\vspace{2mm}

\underline{Case 2-1 : the case the condition (2) does not hold}

\vspace{2mm}

We have
$\lambda^{(n)}=-x_{1,n}+\lan h_n,\lambda \ran$ and
$\Xi_{\iota,n}[\lambda]=\{0,\lambda^{(n)}\}$.

\vspace{2mm}

\underline{Case 2-2 : the case the condition (2) holds}

\vspace{2mm}

In this case, we obtain $p_{n,n-1}=0$ so that $P(n)=P(n-1)$ and
$\lambda^{(n)}=-x_{1,n}+x_{1,n-1}+\lan h_n,\lambda \ran=[\ovl{n+1},\ovl{n}]^{\rm C}_{-P(n-1)}+\lan h_n,\lambda \ran$.
By Lemma \ref{box-lem}, we obtain
\[\what{S}_{l,j}\lambda^{(n)}=
\begin{cases}
0 & {\rm if}\ (l,j)=(1,n), \\
[\ovl{n+1},\ovl{n-1}]^{\rm C}_{-P(n-1)}+\lan h_n,\lambda \ran & {\rm if}\ (l,j)=(1,n-1), \\
\lambda^{(n)} & {\rm otherwise}.
\end{cases}
\]
By Proposition \ref{closednessBC} (iii), we see that
${\rm Tab}_{{\rm C},\iota,n}[\lambda]\cup\{0\}$ is closed under the
action of $\what{S}_{l,j}$ for all $(l,j)\in\mathbb{Z}_{\geq1}\times I$,
which yields
$\Xi_{\iota,n}[\lambda]\subset{\rm Tab}_{{\rm C},\iota,n}[\lambda]\cup\{0\}$.

We can also get the inclusion ${\rm Tab}_{{\rm C},\iota,n}[\lambda]\cup\{0\}\subset \Xi_{\iota,n}[\lambda]$ by a similar way
to the proof of Lemma 5.7 (ii) in \cite{KaN}. \qed

\subsection{Type D-case}

\vspace{2mm}

\nd
\underline{Case 1 : $k<n-2$.}

\vspace{2mm}

\nd
We can show $\Xi_{\iota,k}[\lambda]={\rm Tab}_{{\rm D},\iota,k}[\lambda]\cup\{0\}$ just as in Case 1 of the proofs of type B,C-cases.

\nd
\underline{Case 2 : $k=n-2$}

\vspace{2mm}

\underline{Case 2-1: $\iota^{(n-2)}<\iota^{(n-3)},\ \iota^{(n-1)},\ \iota^{(n)}$}

\vspace{2mm}

In this case, it is easy to check $\lambda^{(n-2)}=-x_{1,n-2}+\lan h_{n-2},\lambda\ran$ and
$\Xi_{\iota,n-2}[\lambda]=\{0,\lambda^{(n-2)}\}$.

\vspace{2mm}

\underline{Case 2-2 : $\iota^{(n-2)}<\iota^{(n-3)},\ \iota^{(n-2)}<\iota^{(n)}$ and $\iota^{(n-2)}>\iota^{(n-1)}$}

\vspace{2mm}

By $\lambda^{(n-2)}=-x_{1,n-2}+x_{1,n-1}+\lan h_{n-2},\lambda\ran$ and a direct calculation,
it holds
\[
\what{S}_{l,j}\lambda^{(n-2)}=
\begin{cases}
-x_{2,n-1} +\lan h_{n-2},\lambda\ran & {\rm if}\ (l,j)=(1,n-1), \\
0 & {\rm if}\ (l,j)=(1,n-2), \\
\lambda^{(n-2)} & {\rm otherwise},
\end{cases}
\]
and
\[
\what{S}_{l,j}(-x_{2,n-1}+\lan h_{n-2},\lambda\ran)=
\begin{cases}
\lambda^{(n-2)} & {\rm if}\ (l,j)=(2,n-1), \\
-x_{2,n-1}+\lan h_{n-2},\lambda\ran & {\rm otherwise}.
\end{cases}
\]
Therefore, it follows $\Xi_{\iota,n-2}[\lambda]=\{0,\lambda^{(n-2)}, -x_{2,n-1}+\lan h_{n-2},\lambda\ran\}$.

\vspace{2mm}

\underline{Case 2-3 : $\iota^{(n-2)}<\iota^{(n-3)},\ \iota^{(n-2)}<\iota^{(n-1)}$ and $\iota^{(n-2)}>\iota^{(n)}$}

\vspace{2mm}

Just as in Case 2-2, we can show
$\Xi_{\iota,n-2}[\lambda]=\{0,\lambda^{(n-2)}, -x_{2,n}+\lan h_{n-2},\lambda\ran\}$.

\vspace{2mm}

\underline{Case 2-4 : $\iota^{(n-2)}<\iota^{(n-1)},\ \iota^{(n-2)}<\iota^{(n)}$ and $\iota^{(n-2)}>\iota^{(n-3)}$}

\vspace{2mm}

It holds 
$\lambda^{(n-2)}=-x_{1,n-2}+x_{1,n-3}+\lan h_{n-2},\lambda\ran= \fbox{$\ovl{n-2}$}^{\rm D}_{-1-P(n-3)}
+\lan h_{n-2},\lambda\ran$.
By Lemma \ref{box-lem}, putting $s:=-1-P(n-3)$ and
$\what{S}(\ovl{j}):=\what{S}_{s+P(j)+n-j-1,j}$, $\what{S}'(\ovl{j}):=\what{S}_{s+P(j)+n-j,j}$ for $j\in[1,n-2]$,
we get the following diagram of actions of $\what{S}$:
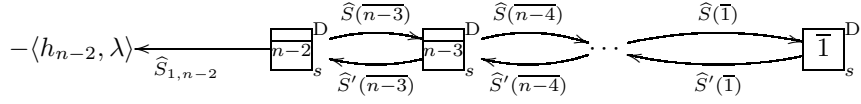
\[
\begin{xy}
(-10,0) *{-\lan h_{n-2},\lambda\ran}="00",
(0,0) *{\ }="0",
(20,0) *{\begin{ytableau}
\scriptstyle \overline{n-2}
\end{ytableau}^{\rm D}_{s}}="r-1r-1",
(40,0)*{\begin{ytableau}
\scriptstyle \overline{n-3}
\end{ytableau}^{\rm D}_{s}}="r-2r-2",
(60,0)*{
\ 
\cdots
\ }="dotdot",
(90,0)*{\begin{ytableau}
\overline{1}
\end{ytableau}^{\rm D}_{s}}="11",
\ar@{->} "r-1r-1";"00"^{\what{S}_{1,n-2}}
\ar@/^/ @{->} "r-1r-1";"r-2r-2"^{\what{S}(\ovl{n-3})}
\ar@/^/ @{->} "r-2r-2";"dotdot"^{\what{S}(\ovl{n-4})}
\ar@/^/ @{->} "dotdot";"11"^{\what{S}(\ovl{1})}
\ar@/_/ @{<-} "r-1r-1";"r-2r-2"_{\what{S}'(\ovl{n-3})}
\ar@/_/ @{<-} "r-2r-2";"dotdot"_{\what{S}'(\ovl{n-4})}
\ar@/_/ @{<-} "dotdot";"11"_{\what{S}'(\ovl{1})}
\end{xy}
\]
Other actions are trivial.
Hence $\Xi_{\iota,n-2}[\lambda]=
\{0\}\cup\{\fbox{$t$}^{\rm D}_{-1-P(n-3)}+\lan h_{n-2},\lambda\ran|\ \ovl{n-2}\leq t\leq \ovl{1} \}$.

\vspace{2mm}

\underline{Case 2-5 : $\iota^{(n-2)}<\iota^{(n-3)}$ and $\iota^{(n-2)}>\iota^{(n-1)}$, $\iota^{(n)}$}

\vspace{2mm}

In this case, we obtain $p_{n-1,n-2}=p_{n,n-2}=1$ and $P(n)=P(n-1)=1+P(n-2)$.
We also get
$\lambda^{(n-2)}=-x_{1,n-2}+x_{1,n-1}+x_{1,n}+\lan h_{n-2},\lambda\ran= \fbox{$n-1$}^{\rm D}_{-P(n-2)}
+\lan h_{n-2},\lambda\ran$. By Lemma \ref{box-lem}, putting $s:=-P(n-2)$ and
$\what{S}(\ovl{j}):=\what{S}_{s+P(j)+n-j-1,j}$, $\what{S}'(\ovl{j}):=\what{S}_{s+P(j)+n-j,j}$ for $j\in[1,n-1]$,
one obtain the following diagram of actions of $\what{S}$:
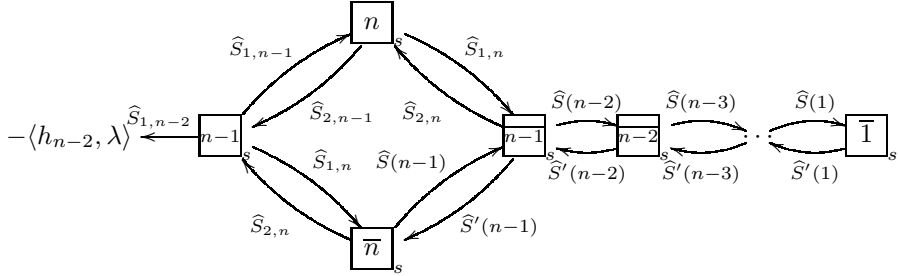
\[
\begin{xy}
(15,0)*{-\lan h_{n-2},\lambda \ran\ }="dots",
(35,0)*{\begin{ytableau}
\scriptstyle n-1
\end{ytableau}_{s}}="r-1",
(55,15)*{\begin{ytableau}
n
\end{ytableau}_{s}}="r",
(55,-15)*{\begin{ytableau}
\overline{n}
\end{ytableau}_{s}}="rr",
(75,0) *{\begin{ytableau}
\scriptstyle \overline{n-1}
\end{ytableau}_{s}}="r-1r-1",
(90,0) *{\begin{ytableau}
\scriptstyle \overline{n-2}
\end{ytableau}_{s}}="r-2r-2",
(105,0)*{\ \cdots\ }="dotdot",
(120,0)*{\begin{ytableau}
\overline{1}
\end{ytableau}_{s}}="11",
\ar@/^/@{->} "r-1";"r"^{\what{S}_{1,n-1}}
\ar@/^/@{->} "r-1";"rr"^{\what{S}_{1,n}}
\ar@/^/@{->} "r";"r-1r-1"^{\what{S}_{1,n}}
\ar@/^/@{->} "rr";"r-1r-1"^{\what{S}(n-1)}
\ar@/^/@{->} "r-1r-1";"r-2r-2"^{\what{S}(n-2)}
\ar@/^/@{->} "r-2r-2";"dotdot"^{\what{S}(n-3)}
\ar@/^/@{->} "dotdot";"11"^{\what{S}(1)}
\ar@{<-} "dots";"r-1"^{\ \ \what{S}_{1,n-2}}
\ar@/_/@{<-} "r-1";"r"_{\what{S}_{2,n-1}}
\ar@/_/@{<-} "r-1";"rr"_{\what{S}_{2,n}}
\ar@/_/@{<-} "r";"r-1r-1"_{\what{S}_{2,n}}
\ar@/_/@{<-} "rr";"r-1r-1"_{\what{S}'(n-1)}
\ar@/_/@{<-} "r-1r-1";"r-2r-2"_{\what{S}'(n-2)}
\ar@/_/@{<-} "r-2r-2";"dotdot"_{\what{S}'(n-3)}
\ar@/_/@{<-} "dotdot";"11"_{\what{S}'(1)}
\end{xy}
\]

Hence, it holds $\Xi_{\iota,n-2}[\lambda]=
\{0\}\cup\{\fbox{$t$}^{\rm D}_{-P(n-2)}+\lan h_{n-2},\lambda\ran\ | n-1\leq t\leq \ovl{1} \}$.

\vspace{2mm}

\underline{Case 2-6 : $\iota^{(n-2)}<\iota^{(n)}$ and $\iota^{(n-2)}>\iota^{(n-3)}$, $\iota^{(n-1)}$}

\vspace{2mm}

In this setting, we get $p_{n-2,n-3}=0$, $p_{n-1,n-2}=1$, $p_{n,n-2}=0$ so that $P(n-2)=P(n-3)$, $P(n-1)=P(n-2)+1$ and
$P(n)=P(n-2)$ by (\ref{LaP}). Thus, one obtain
$\lambda^{(n-2)}=-x_{1,n-2}+x_{1,n-3}+x_{1,n-1}+\lan h_{n-2},\lambda\ran=
[\ovl{n+1},\ovl{n},\ovl{n-2}]^{\rm D}_{-1-P(n-2)}
+\lan h_{n-2},\lambda\ran\in {\rm Tab}_{{\rm D},\iota,n-2}[\lambda]\cup \{0\} $.
Using Lemma \ref{box-lem}, we get
\[
\what{S}_{l,j} \lambda^{(n-2)}=
\begin{cases}
[\ovl{n+1},\ovl{n},\ovl{n-3}]^{\rm D}_{-1-P(n-2)}+\lan h_{n-2},\lambda\ran & {\rm if}\ (l,j)=(1,n-3),\\
[\ovl{n+1},\ovl{n-1},\ovl{n-2}]^{\rm D}_{-1-P(n-2)}+\lan h_{n-2},\lambda\ran & {\rm if}\ (l,j)=(1,n-1),\\
0 & {\rm if}\ (l,j)=(1,n-2),\\
\lambda^{(n-1)} & {\rm otherwise}.
\end{cases}
\]
By Definition \ref{tab-def-lam},
each $T=[\ovl{n+1},j_2,\cdots,j_k]_{-1-P(n-2)}
+\lan h_{n-2},\lambda\ran \in {\rm Tab}_{{\rm D},\iota,n-2}[\lambda]$
satisfies $3\leq k\leq n+1$, $k$ is odd, $\ovl{n}\leq j_2<\cdots<j_k\leq\ovl{1}$
and
if $k=3$ then $j_3\geq \ovl{n-2}$.
It holds
\[
T= \fbox{$\ovl{n+1}$}^{\rm D}_{-P(n-2)+k-2} + \sum^{k}_{i=2} \fbox{$j_i$}^{\rm D}_{-1-P(n-2)+k-i}
+\lan h_{n-2},\lambda\ran.
\]
Recall that
\begin{equation}\label{D-box-ag-2}
\fbox{$\ovl{n+1}$}^{\rm D}_{-P(n-2)+k-2}=
x_{-P(n-2)+k-2+P(n),n},
\end{equation}
\begin{equation}\label{D-box-ag-3}
\fbox{$j_i$}^{\rm D}_{-1-P(n-2)+k-i}=
\begin{cases}
x_{-1-P(n-2)+k-i+P(n-1),n-1}-x_{-P(n-2)+k-i+P(n),n} & {\rm if}\ j_i= \ovl{n}, \\ 
x_{k-i,n-2}-x_{-P(n-2)+k-i+P(n-1),n-1}
-x_{-P(n-2)+k-i+P(n),n} & {\rm if}\ j_i= \ovl{n-1},\\
x_{-1-P(n-2)+k-i+P(|j_i|-1)+n-|j_i|,|j_i|-1}-x_{-1-P(n-2)+k-i+P(|j_i|)+n-|j_i|,|j_i|} & {\rm if}\ j_i\geq \ovl{n-2}.
\end{cases}
\end{equation}
It follows by $k\geq3$ that 
\begin{equation}\label{D-box-ag-4}
-P(n-2)+k-2+P(n)=k-2\geq1.
\end{equation}
By the conditions $k\geq3$ and $\ovl{n}\leq j_2<\cdots<j_k\leq \ovl{1}$
in ${\rm Tab}_{{\rm D},\iota,n-2}[\lambda]$, if
$j_i=\ovl{n}$ then $i=2<k$, which yields
\begin{equation}\label{D-box-ag-5}
-1-P(n-2)+k-i+P(n-1)=k-i\geq1.
\end{equation}
Similarly, by the conditions in ${\rm Tab}_{{\rm D},\iota,n-2}[\lambda]$,
if $j_i= \ovl{n-1}$ then $i<k$ and
\begin{equation}\label{D-box-ag-6}
k-i\geq1.
\end{equation}
If $j_i\geq \ovl{n-2}$ then
\begin{eqnarray}
& &-1-P(n-2)+k-i+P(|j_i|-1)+n-|j_i| \label{D-box-ag-7}\\
&=&-P(n-2)+k-i+P(|j_i|-1)+(n-3)-(|j_i|-1)+1
\geq -P(n-2)+P(n-3)+1 =1.\nonumber
\end{eqnarray}
The inequalities (\ref{D-box-ag-4})-(\ref{D-box-ag-7}) mean the left indices in
(\ref{D-box-ag-2}), (\ref{D-box-ag-3}) are positive.

If $j_i=\ovl{n}$ then $i=2$ and $-x_{-P(n-2)+k-2+P(n),n}$ is cancelled in $T$ by (\ref{D-box-ag-2}).
If $j_i=\ovl{n-1}$ then $i=2$ or $3$. In the case $i=2$, we see that
$-x_{-P(n-2)+k-2+P(n),n}$ is cancelled in $T$ and $-P(n-2)+k-2+P(n-1)=k-1\geq 2$.
In the case $i=3$, it holds $k\geq5$ so that $-P(n-2)+k-3+P(n-1)=k-2\geq2$ and
$-P(n-2)+k-3+P(n)=k-3\geq2$. 
We similarly see that if $j_i\geq\ovl{n-2}$ and $T\neq \lambda^{(n-2)}$ then
$-1-P(n-2)+k-i+P(|j_i|)+n-|j_i|\geq2$.

Just as in the proof of Proposition 4.3 (iii) in \cite{KaN},
we see that for each 
$T=[\ovl{n+1},j_2,\cdots,j_k]_{-1-P(n-2)}+\lan h_{n-2},\lambda \ran\in {\rm Tab}_{{\rm D},\iota,n-2}[\lambda]
\setminus\{\lambda^{(n-2)}\}$,
$j\in I$ and
$m\in\mathbb{Z}_{\geq1}$,
\begin{multline*}
\what{S}_{m,j}(T-\lan h_{n-2},\lambda \ran)=\\
\begin{cases}
[\ovl{n+1},j_2,\cdots,j_{i-1},\ovl{j},j_{i+1},\cdots,j_k]_{-1-P(n-2)}^{\rm D} & {\rm if}\ 
\begin{array}{l}
j<n\ {\rm and\ for\ some}\ i\in[2,k],\ j_i=\ovl{j+1},\ j_{i+1}\neq \ovl{j},\\
 m=-P(n-2)+k-i+P(j)+n-j-2,
\end{array} \\
[\ovl{n+1},j_2,\cdots,j_{i'-1},\ovl{j+1},j_{i'+1},\cdots,j_k]_{-1-P(n-2)}^{\rm D} & {\rm if}\ 
\begin{array}{l}
j<n\ {\rm and\ for\ some}\ i\in[2,k],\ j_i=\ovl{j},\ j_{i+1}\neq \ovl{j+1},\\
 m=-P(n-2)+k-i+P(j)+n-j-1,
\end{array}\\
[\ovl{n+1},\ovl{n},\ovl{n-1},j_2,\cdots,j_k]_{-1-P(n-2)}^{\rm D} & {\rm if\ }
j=n,\ j_2\neq\ovl{n},\ \ovl{n-1}\ {\rm and}\ m=k-2,
 \\
[\ovl{n+1},j_4,\cdots,j_k]_{-1-P(n-2)}^{\rm D} & {\rm if\ }
j=n,\ j_2=\ovl{n},\ j_3=\ovl{n-1}\ {\rm and}\ m=k-3,
\\
T-\lan h_{n-2},\lambda \ran & {\rm otherwise}.
\end{cases}
\end{multline*}
Thus, we can verify $\Xi_{\iota,n-2}[\lambda]\subset{\rm Tab}_{{\rm D},\iota,n-2}[\lambda]\cup \{0\}$.
By a similar way to the proof of Lemma 5.13 in \cite{KaN},
we can also verify ${\rm Tab}_{{\rm D},\iota,n-2}[\lambda]\cup \{0\}\subset \Xi_{\iota,n-2}[\lambda]$.

\vspace{2mm}

\underline{Case 2-7 : $\iota^{(n-2)}<\iota^{(n-1)}$ and $\iota^{(n-2)}>\iota^{(n-3)}$, $\iota^{(n)}$}

\vspace{2mm}

In this setting, we get $p_{n-2,n-3}=0$, $p_{n-1,n-2}=0$, $p_{n,n-2}=1$ so that $P(n-2)=P(n-3)$, $P(n-1)=P(n-2)$ and
$P(n)=P(n-2)+1$ by (\ref{LaP}). Thus,
$\lambda^{(n-2)}=-x_{1,n-2}+x_{1,n-3}+x_{1,n}+\lan h_{n-2},\lambda\ran=
[\ovl{n+1},\ovl{n-2}]^{\rm D}_{-1-P(n-2)}
+\lan h_{n-2},\lambda\ran$.
By a similar way to Case 2-6, we can verify
$\Xi_{\iota,n-2}[\lambda]={\rm Tab}_{{\rm D},\iota,n-2}[\lambda]\cup \{0\}$.

\vspace{2mm}

\underline{Case 2-8 : $\iota^{(n-2)}>\iota^{(n-3)}$, $\iota^{(n-1)}$, $\iota^{(n)}$}

\vspace{2mm}

It holds $p_{n-2,n-3}=0$, $p_{n-1,n-2}=p_{n,n-2}=1$ so that $P(n-2)=P(n-3)$, $P(n)=P(n-1)=P(n-2)+1$ and
$\lambda^{(n-2)}=-x_{1,n-2}+x_{1,n-3}+x_{1,n-1}+x_{1,n}+\lan h_{n-2},\lambda\ran=
[1,2,\cdots,n-3,n-1]^{\rm D}_{-P(n-3)}
+\lan h_{n-2},\lambda\ran$. The following is a consequence of Lemma \ref{box-lem}:
\[
\what{S}_{l,j}(\lambda^{(n-2)})=
\begin{cases}
[1,2,\cdots,n-4,n-2,n-1]^{\rm D}_{-P(n-3)}+\lan h_{n-2},\lambda\ran & {\rm if}\ (l,j)=(1,n-3),\\
[1,2,\cdots,n-3,n]^{\rm D}_{-P(n-3)}+\lan h_{n-2},\lambda\ran & {\rm if}\ (l,j)=(1,n-1),\\
[1,2,\cdots,n-3,\ovl{n}]^{\rm D}_{-P(n-3)}+\lan h_{n-2},\lambda\ran & {\rm if}\ (l,j)=(1,n),\\
0 & {\rm if}\ (l,j)=(1,n-2),\\
\lambda^{(n-2)} & {\rm otherwise}.
\end{cases}
\]
Considering Proposition \ref{closednessD} (i), (ii), we see that
${\rm Tab}_{{\rm D},\iota,n-2}[\lambda]\cup \{0\}$ is closed under the action of $\what{S}_{l,j}$
for all $(l,j)\in\mathbb{Z}_{\geq1}\times I$, which yields 
$\Xi_{\iota,n-2}[\lambda]\subset{\rm Tab}_{{\rm D},\iota,n-2}[\lambda]\cup \{0\}$.
The inclusion ${\rm Tab}_{{\rm D},\iota,n-2}[\lambda]\cup \{0\}\subset \Xi_{\iota,n-2}[\lambda]$
follows from a similar argument to the proof of Lemma 5.12 in \cite{KaN}.

\vspace{2mm}

\nd
\underline{Case 3 : $k=n-1$}

\vspace{2mm}

\underline{Case 3-1 : $\iota^{(n-1)}< \iota^{(n-2)}$}

\vspace{2mm}

In this case, we can easily check $\lambda^{(n-1)}=-x_{1,n-1}+\lan h_{n-1},\lambda\ran$ and
$\Xi_{\iota,n-1}[\lambda]=\{0,\lambda^{(n-1)}\}$.

\vspace{2mm}

\underline{Case 3-2 : $\iota^{(n-1)}> \iota^{(n-2)}$}

\vspace{2mm}

We get
$\lambda^{(n-1)}=-x_{1,n-1}+x_{1,n-2}+\lan h_{n-1},\lambda\ran
=[\ovl{n+1},\ovl{n-1}]^{\rm D}_{-P(n-2)}+\lan h_{n-1},\lambda\ran$.
By Lemma \ref{box-lem}, we see that
\[
\what{S}_{l,j} \lambda^{(n-1)}=
\begin{cases}
[\ovl{n+1},\ovl{n-2}]^{\rm D}_{-P(n-2)}+\lan h_{n-1},\lambda\ran & {\rm if}\ (l,j)=(1,n-2),\\
0 & {\rm if}\ (l,j)=(1,n-1),\\
\lambda^{(n-1)} & {\rm otherwise}.
\end{cases}
\]
Combining with Proposition \ref{closednessD} (iii), we also see that
${\rm Tab}_{{\rm D},\iota,n-1}[\lambda]\cup\{0\}$ is closed under the action of $\what{S}_{l,j}$ for all $(l,j)\in\mathbb{Z}_{\geq1}\times I$,
which means $\Xi_{\iota,n-1}[\lambda]\subset {\rm Tab}_{{\rm D},\iota,n-1}[\lambda]\cup\{0\}$.
The inclusion ${\rm Tab}_{{\rm D},\iota,n-1}[\lambda]\cup\{0\}\subset\Xi_{\iota,n-1}[\lambda]$ can be proved
just as in the proof of Lemma 5.13 (i) in \cite{KaN}. 

\vspace{2mm}

\nd
\underline{Case 4 : $k=n$}

\vspace{2mm}

\underline{Case 4-1 : $\iota^{(n)}< \iota^{(n-2)}$}

\vspace{2mm}

In this case, we obtain $\lambda^{(n)}=-x_{1,n}+\lan h_n,\lambda\ran$ and
$\Xi_{\iota,n}[\lambda]=\{0,\lambda^{(n)}\}$.

\vspace{2mm}

\underline{Case 4-2 : $\iota^{(n)}> \iota^{(n-2)}$}

\vspace{2mm}

It holds
$\lambda^{(n)}=-x_{1,n}+x_{1,n-2}+\lan h_n,\lambda\ran
=[\ovl{n+1},\ovl{n},\ovl{n-1}]^{\rm D}_{-P(n-2)}+\lan h_n,\lambda\ran$ and
\[
\what{S}_{l,j} \lambda^{(n)}=
\begin{cases}
[\ovl{n+1},\ovl{n},\ovl{n-2}]^{\rm D}_{-P(n-2)}+\lan h_n,\lambda\ran & {\rm if}\ (l,j)=(1,n-2),\\
0 & {\rm if}\ (l,j)=(1,n),\\
\lambda^{(n)} & {\rm otherwise}.
\end{cases}
\]
Just as in Case 3-2, our claim
$\Xi_{\iota,n}[\lambda]={\rm Tab}_{{\rm D},\iota,n}[\lambda]\cup\{0\}$ follows. \qed

\vspace{3mm}

Hence, we obtain $\Xi_{\iota}[\lambda]={\rm Tab}_{{\rm X,\iota}}[\lambda]\cup {\rm Tab}_{{\rm X,\iota}}$.
It is easy to verify that for any $\lambda\in P_+$, the pair $(\iota,\lambda)$ satisfies the ample condition by the explicit forms
of ${\rm Tab}_{{\rm X,\iota}}$, ${\rm Tab}_{{\rm X,\iota}}[\lambda]$ in Definition \ref{tab-def}, \ref{tab-def-lam}.
Consequently,
the proof of Theorem \ref{thm3} is completed.

\subsection{Proof of Theorem \ref{thm4}}

We proved in \cite{KaN} that
for $\varphi=\sum_{i}\varphi_ix_i\in \Xi^{(\infty)}_{\iota}$,
if $i^{(-)} = 0$ then $\varphi_i\geq0$. 
 Fixing $k\in I$, we show for $\varphi=\sum_{i}\varphi_ix_i\in {\rm Tab}_{X,\iota,k}[\lambda]\setminus\{\lambda^{(k)}\}
 -\langle h_k,\lambda\rangle$,
\begin{equation}\label{cl-1}
{\rm if}\ i^{(-)} = 0\ {\rm then}\ \varphi_i\geq0,
\end{equation}
where for a set $S$ of linear functions, the set $S-\lan h_k,\lambda\ran$ is defined as
$S-\lan h_k,\lambda\ran=\{f-\lan h_k,\lambda\ran | f\in S \}$.
 In the previous subsections, we proved that
\begin{eqnarray*}
{\rm Tab}_{{\rm X},\iota,k}[\lambda]\setminus\{\lambda^{(k)}\} - \lan h_k,\lambda\ran
&=&\Xi_{\iota,k}[\lambda]\setminus\{\lambda^{(k)},0\}-\lan h_k,\lambda\ran\\
&=&\{\what{S}_{j_l}\cdots\what{S}_{j_1}\lambda^{(k)} -\lan h_k,\lambda\ran\
|l\in\mathbb{Z}_{\geq1},j_1,\cdots,j_l\in\mathbb{Z}_{\geq1}\} \setminus\{\xi^{(k)},-\lan h_k,\lambda\ran \}.
\end{eqnarray*}
 for any $\lambda\in P_+$. Hence, if the condition (\ref{cl-1}) holds for all $k\in I$ then
by the definitions of $S_j$, $\what{S}_j$ ((\ref{Sk}),(\ref{Shat})), it follows
\[
{\rm Tab}_{{\rm X},\iota,k}[\lambda]\setminus\{\lambda^{(k)}\} - \lan h_k,\lambda\ran
=
\{S_{j_l}\cdots S_{j_1}\xi^{(k)}
|l\in\mathbb{Z}_{\geq0},j_1,\cdots,j_l\in\mathbb{Z}_{\geq1}\} \setminus\{\xi^{(k)} \}=\Xi_{\iota}^{(k)}\setminus\{\xi^{(k)}\},
\]
and strict positivity condition holds. In particular, we obtain ${\rm Tab}_{{\rm X},\iota,k}[0]=\Xi_{\iota}^{(k)}$. Therefore, all we need to prove is the condition (\ref{cl-1}).
Let us recall the conditions (1), (2) in (\ref{cond12}).

\vspace{2mm}

\nd
\underline{Type A-case}

\vspace{2mm}

\underline{Case 1 : both (1) and (2) do not hold}

\vspace{2mm}

In this case, by ${\rm Tab}_{{\rm A},\iota,k}[\lambda]\setminus\{\lambda^{(k)}\}=\phi$, our claim (\ref{cl-1}) is clear.

\vspace{2mm}

\underline{Case 2 : (1) holds and (2) does not hold}

\vspace{2mm}

In this case, it follows from $\lambda^{(k)}=\fbox{$k+1$}^{\rm A}_{1-P(k+1)}+\lan h_k,\lambda\ran$ and Definition \ref{tab-def-lam} that
${\rm Tab}_{{\rm A},\iota,k}[\lambda]\setminus\{\lambda^{(k)}\}=
\{
\fbox{$t$}^{\rm A}_{1-P(k+1)} +\langle h_k,\lambda\rangle  | k+1< t \leq n+1
\}$.
It is easy to see that $\fbox{$t$}^{\rm A}_{1-P(k+1)}=x_{1-P(k+1)+P(t),t}-x_{2-P(k+1)+P(t-1),t-1}$ and
$2-P(k+1)+P(t-1)\geq2$, which means (\ref{cl-1}).

\vspace{2mm}

\underline{Case 3 : (2) holds and (1) does not hold}

\vspace{2mm}

In \ref{A-proof}, we have seen that 
each element $T=[j_1,\cdots,j_{k-1},k+1,\cdots,n+1]^{\rm A}_{-P(k-1)-n+k}+\lan h_k,\lambda
\ran\in {\rm Tab}_{{\rm A},\iota,k}[\lambda]\setminus\{\lambda^{(k)}\}$
can be written
as $T=[1,\cdots,l-1,l+1,\cdots,n+1]^{\rm A}_{-P(k-1)-n+k}+\lan h_k,\lambda \ran$ with some $l\in[1,k-1]$.
We see that $T=x_{-P(k-1)-l+k+1+P(l-1),l-1}-x_{-P(k-1)-l+k+P(l)+1,l}$ and
\[
-P(k-1)-l+k+P(l)+1
=-P(k-1)+(k-1)-l+2+P(l)\geq -P(k-1)+2+P(k-1)=2,
\]
which implies (\ref{cl-1}).

\vspace{2mm}

\underline{Case 4 : both (1) and (2) hold}

\vspace{2mm}

For each $T=[j_1,\cdots,j_k]^{\rm A}_{-P(k-1)}+\lan h_k,\lambda\ran\in {\rm Tab}_{{\rm A},\iota,k}[\lambda]\setminus\{\lambda^{(k)}\}$,
it follows
\[
T-\lan h_k,\lambda\ran=\sum^{k}_{i=1} \fbox{$j_i$}^{\rm A}_{-P(k-1)+k-i} 
=\sum^{k}_{i=1} (x_{-P(k-1)+k-i+P(j_i),j_i}-x_{-P(k-1)+k-i+1+P(j_i-1),j_{i}-1}).
\]
For $i\in[1,k-1]$ such that $j_i=i$, the condition $1\leq j_1<\cdots<j_i$ in ${\rm Tab}_{{\rm A},\iota,k}[\lambda]$
implies $j_1=1$, $j_2=2$, $\cdots$, $j_{i-1}=i-1$. In particular, it holds $j_{i}=j_{i-1}+1$ and
we have seen in (\ref{A-pr02}) that
$x_{-P(k-1)+k-i+1+P(j_i-1),j_{i}-1}$ is cancelled in $T$.
For $i\in[1,k-1]$ such that $j_i>i$, we can show $-P(k-1)+k-i+1+P(j_i-1)\geq2$ just as in (\ref{pr-A-01}).

For $i=k$, if $j_k=j_{k-1}+1$ then $x_{-P(k-1)+k-k+1+P(j_k-1),j_{k}-1}$ is cancelled in $T$.
If $j_k>j_{k-1}+1$ then we see that $j_k>k+1$ by $T\neq \lambda^{(k)}$, which yields
\[
-P(k-1)+k-k+1+P(j_k-1)=-P(k-1)+1+P(j_k-1)\geq -P(k-1)+1+P(k+1) =2.
\]
Therefore, the condition (\ref{cl-1}) holds.

\vspace{2mm}

\nd
\underline{Type B-case}

\vspace{2mm}

We fix $k\in[1,n]$.

\vspace{2mm}

\underline{Case 1 : both (1) and (2) do not hold}

\vspace{2mm}

By ${\rm Tab}_{{\rm B},\iota,k}[\lambda]\setminus\{\lambda^{(k)}\}=\phi$, our claim (\ref{cl-1}) is clear.

\vspace{2mm}

\underline{Case 2 : (1) holds and (2) does not hold}

\vspace{2mm}

In this setting, it holds $p_{k+1,k}=1$ so that $P(k+1)=P(k)+1$.
By Definition \ref{tab-def-lam}, we have
${\rm Tab}_{{\rm B},\iota,k}[\lambda] =\{\fbox{$t$}^{\rm B}_{1-P(k+1)} +\langle h_k,\lambda\rangle
| k+1\leq t \leq \overline{1} \}$. By Definition \ref{box-def} (ii), we see that
$\fbox{$k+1$}^{\rm B}_{1-P(k+1)}=x_{1,k+1}-x_{1,k}=\lambda^{(k)}-\langle h_k,\lambda\rangle$.
Hence,
\[
{\rm Tab}_{{\rm B},\iota,k}[\lambda]\setminus\{\lambda^{(k)}\}
=\{\fbox{$t$}^{\rm B}_{1-P(k+1)} +\langle h_k,\lambda\rangle
| k+1< t \leq \overline{1} \}.
\]
We also see that for $t$ ($k+1< t \leq \overline{1}$), 
\[
\fbox{$t$}^{\rm B}_{1-P(k+1)}=
\begin{cases}
x_{1-P(k+1)+P(t),t}-x_{-P(k+1)+P(t-1)+2,t-1} & {\rm if}\ t\leq n,\\
x_{-P(k+1)+P(|t|-1)+n-|t|+2,|t|-1}-x_{-P(k+1)+P(|t|)+n-|t|+2,|t|} & {\rm if}\ t\geq \ovl{n},
\end{cases}
\]
and if $t\leq n$ then
\[
-P(k+1)+P(t-1)+2 \geq2,\]
if $t\geq \ovl{n}$ then
\[ -P(k+1)+P(|t|)+n-|t|+2\geq -P(k+1)+P(n)+2\geq2.
\]
Hence, the condition (\ref{cl-1}) holds.

\vspace{2mm}

\underline{Case 3 : (2) holds and (1) does not hold}

\vspace{2mm}

In this case, we have $p_{k,k-1}=0$ so that $P(k)=P(k-1)$ and
\[
{\rm Tab}_{{\rm B},\iota,k}[\lambda] =
\{
\fbox{$t$}^{\rm B}_{-P(k-1)-n+k}+\langle h_k,\lambda\rangle | \overline{k}\leq t\leq \overline{1}
\}.
\]
By Definition \ref{box-def} (ii), we see that
$\fbox{$\ovl{k}$}^{\rm B}_{-P(k-1)-n+k}=x_{1,k-1}-x_{1,k}=\lambda^{(k)}-\langle h_k,\lambda\rangle$.
Thus,
\[
{\rm Tab}_{{\rm B},\iota,k}[\lambda]\setminus\{\lambda^{(k)}\} =
\{
\fbox{$t$}^{\rm B}_{-P(k-1)-n+k}+\langle h_k,\lambda\rangle | \overline{k}< t\leq \overline{1}
\}.
\]
For $t$ ($\overline{k}< t\leq \overline{1}$), we have
\[
\fbox{$t$}^{\rm B}_{-P(k-1)-n+k}=x_{-P(k-1)+P(|t|-1)+k-|t|+1,|t|-1}-x_{-P(k-1)+P(|t|)+k-|t|+1,|t|}
\]
and
$-P(k-1)+P(|t|)+k-|t|+1=-P(k-1)+P(|t|)+(k-1)-|t|+2\geq -P(k-1)+P(k-1)+2=2$.
So the condition (\ref{cl-1}) holds.

\vspace{2mm}

\underline{Case 4 : both (1) and (2) hold}

\vspace{2mm}

For $T=[j_1,\cdots,j_k]_{-P(k-1)}^{\rm B}+\lan h_k,\lambda\ran \in {\rm Tab}_{{\rm B},\iota,k}\setminus \{\lambda^{(k)}\}$,
we have
$[j_1,\cdots,j_k]^{\rm B}_{-P(k-1)}=\sum_{i=1}^{k}\fbox{$j_i$}^{\rm B}_{-P(k-1)+k-i}$ and
Definition \ref{box-def} says
\[
\fbox{$j_i$}^{\rm B}_{-P(k-1)+k-i}=
\begin{cases}
x_{-P(k-1)+k-i+P(j_i),j_i}-x_{-P(k-1)+k-i+P(j_i-1)+1,j_i-1} & {\rm if}\ j_i\leq n,\\
x_{-P(k-1)+k-i+P(|j_i|-1)+n-|j_i|+1,|j_i|-1}-x_{-P(k-1)+k-i+P(|j_i|)+n-|j_i|+1,|j_i|} & {\rm if}\ j_i\geq \ovl{n}.
\end{cases}
\]
For $i\in[1,k]$ such that $j_i\leq n$, if $j_i=j_{i-1}+1$ ($j_0:=0$) then
the summand $x_{-P(k-1)+k-i+P(j_i-1)+1,j_i-1}$ in $[j_1,\cdots,j_k]^{\rm B}_{-P(k-1)}$ is cancelled
by a similar argument to (\ref{A-pr02}).
The inequalities
(\ref{BC-pr05}), (\ref{BC-pr03}) and (\ref{BC-pr04}) in the proof of Proposition \ref{closednessBC}
mean that the condition (\ref{cl-1}) holds.

\vspace{2mm}

\nd
\underline{Type C-case}

\vspace{2mm}

We can check the condition (\ref{cl-1}) for $k\in[1,n-1]$ just as in Type B-case. Hence,
let us check the condition for $k=n$. 
If $\iota^{(n)}<\iota^{(n-1)}$ then ${\rm Tab}_{{\rm C},\iota,n}\setminus\{\lambda^{(n)}\}=\phi$ and
(\ref{cl-1}) is clear.
Thus, we suppose $\iota^{(n)}>\iota^{(n-1)}$. In this setting, by the argument (in particular (\ref{BC-pr08}), (\ref{BC-pr09}))
in the proof of 
Proposition \ref{closednessBC} (iii), 
we see that the condition (\ref{cl-1}) holds.

\vspace{2mm}

\nd
\underline{Type D-case}

\vspace{2mm}

For $k\in[1,n-3]$, the condition (\ref{cl-1}) holds by the inequalities
(\ref{pr-02e}), (\ref{pr-02f}), (\ref{pr-02g}), (\ref{pr-02h}) and a similar argument to Type B-case.

Next, we suppose $k=n-2$. In the following three cases,  the condition (\ref{cl-1}) clearly holds by Definition \ref{tab-def-lam}:
\begin{itemize}
\item 
$\iota^{(n-2)}<\iota^{(n-3)},\ \iota^{(n-2)}<\iota^{(n-1)}$ and
$\iota^{(n-2)}<\iota^{(n)}$,\\ \vspace{-4mm}
\item
$\iota^{(n-2)}<\iota^{(n-3)},\ \iota^{(n-2)}>\iota^{(n-1)}$ and
$\iota^{(n-2)}<\iota^{(n)}$,\\ \vspace{-4mm}
\item
$\iota^{(n-2)}<\iota^{(n-3)},\ \iota^{(n-2)}<\iota^{(n-1)}$ and
$\iota^{(n-2)}>\iota^{(n)}$.
\end{itemize}
In the following two cases, the condition (\ref{cl-1}) follows from
Definition \ref{tab-def-lam} and
a similar argument to Case 2, 3 of type B:
\begin{itemize}
\item$\iota^{(n-2)}>\iota^{(n-3)}$, $\iota^{(n-2)}<\iota^{(n-1)}$,
$\iota^{(n-2)}<\iota^{(n)}$,\\ \vspace{-4mm}
\item
$\iota^{(n-2)}<\iota^{(n-3)}$, $\iota^{(n-2)}>\iota^{(n-1)}$,
$\iota^{(n-2)}>\iota^{(n)}$.
\end{itemize}
In the cases
\begin{itemize}
\item$\iota^{(n-2)}>\iota^{(n-3)}$, $\iota^{(n-2)}>\iota^{(n-1)}$,
$\iota^{(n-2)}<\iota^{(n)}$,\\ \vspace{-4mm}
\item
$\iota^{(n-2)}>\iota^{(n-3)}$, $\iota^{(n-2)}<\iota^{(n-1)}$,
$\iota^{(n-2)}>\iota^{(n)}$,
\end{itemize}
by Definition \ref{tab-def-lam},
each element in ${\rm Tab}_{{\rm D,\iota,n-2}}[\lambda]$ is written as
$[\ovl{n+1},j_2,\cdots,j_s]_{-1-P(n-2)}+\lan h_{n-2},\lambda\ran$ with some positive integer $s$
and $j_2,\cdots,j_s\in J_{\rm D}$.
Considering the explicit forms of boxes (\ref{D-box-ag-2}), (\ref{D-box-ag-3}) and
the argument after (\ref{D-box-ag-7}) in Case 2.6, we see that the condition (\ref{cl-1}) holds.

In the case
\begin{itemize}
\item$\iota^{(n-2)}>\iota^{(n-3)}$, $\iota^{(n-2)}>\iota^{(n-1)}$,
$\iota^{(n-2)}>\iota^{(n)}$,
\end{itemize}
by Definition \ref{tab-def-lam},
each element $T$ in ${\rm Tab}_{{\rm D,\iota,n-2}}[\lambda]\setminus \{\lambda^{(n-2)}\}$ is written as
$T=[j_1,\cdots,j_{n-2}]^{\rm D}_{-P(n-3)}+\lan h_{n-2},\lambda\ran
=\sum_{i=1}^{n-2}\fbox{$j_i$}^{\rm D}_{-P(n-3)+n-2-i}+\lan h_{n-2},\lambda\ran$.
It follows from the explicit form (\ref{D-box}) of the boxes and
the inequalities (\ref{pr-02e})-(\ref{pr-02h}) for $k=n-2$ in the proof of Proposition \ref{closednessD} (i), (ii) that
the condition (\ref{cl-1}) holds.

Finally, we suppose $k=n-1$ or $n$. 
In the case $\iota^{(k)}<\iota^{(n-2)}$, the condition (\ref{cl-1}) clearly holds by
Definition \ref{tab-def-lam}. Thus, we consider the case
$\iota^{(k)}>\iota^{(n-2)}$.
Each element $T$ in ${\rm Tab}_{{\rm D,\iota,k}}[\lambda]\setminus \{\lambda^{(k)}\}$
is written as
\begin{eqnarray*}
T-\lan h_k,\lambda\ran
&=&[\ovl{n+1},j_2,\cdots,j_{s}]^{\rm D}_{-P(n-2)}-\lan h_k,\lambda\ran\\
&=&\fbox{$\ovl{n+1}$}^{\rm D}_{-P(n-2)+s-1} + \sum^{s}_{i=2} \fbox{$j_i$}^{\rm D}_{-P(n-2)+s-i}.
\end{eqnarray*}
with some integer $s$ and $j_2,\cdots,j_s\in \{\ovl{n},\cdots,\ovl{1}\}$.
Replacing $k$ with $s$ in the definition (\ref{D-box-ag}) of the boxes and the inequalities
(\ref{D-pr-03-b})-(\ref{D-pr-03-d}), we see that the condition (\ref{cl-1}) holds. \qed

\end{document}